\documentclass[10pt,reqno]{extarticle}

\usepackage{geometry}
\usepackage{graphicx}

\usepackage{tocloft}
\usepackage[utf8]{inputenc}
\usepackage[english]{babel}

\usepackage{amssymb}
\usepackage{amsmath}
\usepackage{amsthm}
\usepackage{mathrsfs}
\usepackage{mathtools}

\usepackage[all,arc]{xy}
\usepackage{tikz}
\usepackage{tikz-cd}
\usepackage{float}
\usepackage{subcaption}

\usepackage{soul}
\usepackage{cancel}
\usepackage{slashed}
\usepackage{tensor}
\usepackage{scalerel}

\usepackage{color}
\usepackage{enumitem}
\usepackage{hyperref}
\hypersetup{
  colorlinks   = true, 
  urlcolor     = blue, 
  linkcolor    = blue, 
  citecolor   = red 
}

\usepackage[final]{pdfpages}
\usepackage{showlabels}

\numberwithin{equation}{section}
\allowdisplaybreaks

\newtheorem{theorem}{Theorem}[section]
\newtheorem{definition}[theorem]{Definition}

\newtheorem{remark}[theorem]{Remark}
\newtheorem{lemma}[theorem]{Lemma}
\newtheorem{proposition}[theorem]{Proposition}

\newcommand{\uL}{\underline{L} \makebox[0ex]{}}
\newcommand{\uchi}{\underline{\chi} \makebox[0ex]{}}
\newcommand{\ueta}{\underline{\eta} \makebox[0ex]{}}

\newcommand{\slashg}{g \mkern-8.2mu \scaleto{\boldsymbol{\slash}}{1.6ex} \mkern+1mu \makebox[0ex]{}}  
\newcommand{\subslashg}{g \mkern-7.7mu \scaleto{\boldsymbol{\slash}}{1ex} \mkern+1mu \makebox[0ex]{}}  
\newcommand{\tr}{\mathrm{tr}}              
\newcommand{\hatchi}{\hat{\chi} \makebox[0ex]{}}
\newcommand{\hatuchi}{\underline{\hat{\chi}} \makebox[0ex]{}}

\newcommand{\slashdiv}{\mathrm{div} \mkern-16.5mu \raisebox{+0.2ex}[0ex][0ex]{$\scaleto{\boldsymbol{\slash}}{2ex}$} \mkern+9mu \makebox[0ex]{}}        

\newcommand{\dvol}{\mathrm{dvol} \makebox[0ex]{}}

\newcommand{\ucalH}{\underline{\mathcal{H}} \makebox[0ex]{}}
\newcommand{\uomega}{\underline{\omega} \makebox[0ex]{}}

\newcommand{\ed}{\mathrm{d} \makebox[0ex]{}}

\newcommand{\slashDelta}{\Delta \mkern-10mu \raisebox{0.4ex}[0ex][0ex]{$\scaleto{\boldsymbol{\slash}}{1.7ex}$} \mkern+2mu \makebox[0ex]{}}

\newcommand{\slashgl}[1]{{}^{#1} \mkern-2.5mu g \mkern-8.2mu \scaleto{\boldsymbol{\slash}}{1.6ex} \mkern+1mu \makebox[0ex]{}}

\newcommand{\chil}[1]{{}^{#1}\mkern-1.5mu \chi \makebox[0ex]{}}
\newcommand{\uchil}[1]{{}^{#1}\mkern-1.5mu \underline{\chi} \makebox[0ex]{}}
\newcommand{\etal}[1]{{}^{#1}\mkern-1.5mu \eta \makebox[0ex]{}}
\newcommand{\uomegal}[1]{{}^{#1}\mkern-1.5mu \underline{\omega} \makebox[0ex]{}}

\newcommand{\hatchil}[1]{{}^{#1}\mkern-1.5mu \hat{\chi} \makebox[0ex]{}}
\newcommand{\hatuchil}[1]{{}^{#1}\mkern-1.5mu \underline{\hat{\chi}} \makebox[0ex]{}}

\newcommand{\Kl}[1]{{}^{#1} \mkern-2.5mu K \makebox[0ex]{}}

\newcommand{\mul}[1]{{}^{#1} \mkern-1.5mu \mu \makebox[0ex]{}}

\newcommand{\bSigma}{\bar{\Sigma} \makebox[0ex]{}}

\newcommand{\dL}{\dot{L} \makebox[0ex]{}}
\newcommand{\duL}{\underline{\dot{L}} \makebox[0ex]{}}

\newcommand{\dmul}[1]{{}^{#1} \mkern-1.5mu \dot{\mu} \makebox[0ex]{}}

\newcommand{\dchil}[1]{{}^{#1}\mkern-1.5mu \dot{\chi} \makebox[0ex]{}}
\newcommand{\duchil}[1]{{}^{#1}\mkern-1.5mu \underline{\dot{\chi}} \makebox[0ex]{}}
\newcommand{\detal}[1]{{}^{#1}\mkern-1.5mu \dot{\eta} \makebox[0ex]{}}

\newcommand{\bslashdiv}{\bar{\mathrm{div}} \mkern-16.5mu \raisebox{+0.2ex}[0ex][0ex]{$\scaleto{\boldsymbol{\slash}}{2ex}$} \mkern+9mu \makebox[0ex]{}} 

\newcommand{\bslashd}{\bar{\mathrm{d}} \mkern-8.9mu \raisebox{+0.2ex}[0ex][0ex]{$\scaleto{\boldsymbol{\slash}}{2ex}$} \mkern+2.3mu \makebox[0ex]{}}

\newcommand{\bslashDelta}{\bar{\Delta} \mkern-10mu \raisebox{0.4ex}[0ex][0ex]{$\scaleto{\boldsymbol{\slash}}{1.7ex}$} \mkern+2mu \makebox[0ex]{}}

\newcommand{\bL}{\bar{L} \makebox[0ex]{}}
\newcommand{\buL}{\underline{\bar{L}} \makebox[0ex]{}}

\newcommand{\br}{\bar{r} \makebox[0ex]{}}

\newcommand{\bslashgl}[1]{{}^{#1} \mkern-1.5mu \bar{g} \mkern-8.2mu \scaleto{\boldsymbol{\slash}}{1.6ex} \mkern+1mu \makebox[0ex]{}}
\newcommand{\subbslashgl}[1]{{}^{#1} \mkern-1.5mu \bar{g} \mkern-7.7mu \scaleto{\boldsymbol{\slash}}{1ex} \mkern+1mu \makebox[0ex]{}}

\newcommand{\bchil}[1]{{}^{#1}\mkern-1.5mu \bar{\chi} \makebox[0ex]{}}
\newcommand{\buchil}[1]{{}^{#1}\mkern-1.5mu \underline{\bar{\chi}} \makebox[0ex]{}}
\newcommand{\btal}[1]{{}^{#1}\mkern-1.5mu \bar{\eta} \makebox[0ex]{}}
\newcommand{\buomegal}[1]{{}^{#1}\mkern-1.5mu \underline{\bar{\omega}} \makebox[0ex]{}}

\newcommand{\bhatchil}[1]{{}^{#1}\mkern-1.5mu \widehat{\bar{\chi}} \makebox[0ex]{}}
\newcommand{\bhatuchil}[1]{{}^{#1}\mkern-1.5mu \underline{\widehat{\bar{\chi}}} \makebox[0ex]{}}

\newcommand{\btr}{\bar{\mathrm{tr}} \makebox[0ex]{}}

\newcommand{\bKl}[1]{{}^{#1} \mkern-2.5mu \bar{K} \makebox[0ex]{}}

\newcommand{\bmul}[1]{{}^{#1} \mkern-1.5mu \bar{\mu} \makebox[0ex]{}}

\newcommand{\bal}[1]{{}^{#1} \mkern-1.5mu \bar{a} \makebox[0ex]{}}

\newcommand{\rhol}[1]{{}^{#1} \mkern-1.5mu \rho \makebox[0ex]{}}
\newcommand{\subslashgl}[1]{{}^{#1} \mkern-1.5mu g \mkern-7.7mu \scaleto{\boldsymbol{\slash}}{1ex} \mkern+1mu \makebox[0ex]{}}

\newcommand{\slashepsilon}{\epsilon \mkern-5.8mu  \raisebox{0.12ex}[0ex][0ex]{$\scaleto{\boldsymbol{\slash}}{1.4ex}$} \mkern+0.5mu \makebox[0ex]{}}

\newcommand{\ualpha}{\underline{\alpha} \makebox[0ex]{}}
\newcommand{\ubeta}{\underline{\beta} \makebox[0ex]{}}

\newcommand{\slashd}{\mathrm{d} \mkern-8.9mu \raisebox{+0.2ex}[0ex][0ex]{$\scaleto{\boldsymbol{\slash}}{2ex}$} \mkern+2.3mu \makebox[0ex]{}}

\newcommand{\slashcurl}{\mathrm{curl} \mkern-16.5mu \raisebox{+0.2ex}[0ex][0ex]{$\scaleto{\boldsymbol{\slash}}{2ex}$} \mkern+9mu \makebox[0ex]{}}

\newcommand{\slashnabla}{\nabla \mkern-13mu \raisebox{0.3ex}[0ex][0ex]{$\scaleto{\boldsymbol{\slash}}{1.7ex}$} \mkern+5mu \makebox[0ex]{}}

\newcommand{\lie}{\mathcal{L} \makebox[0ex]{}}

\newcommand{\slashGamma}{\Gamma \mkern-9.5mu \raisebox{0.4ex}[0ex][0ex]{$\scaleto{\boldsymbol{\slash}}{1.7ex}$} \mkern+1mu \makebox[0ex]{}}

\newcommand{\us}{\underline{s} \makebox[0ex]{}}
\newcommand{\uC}{\underline{C} \makebox[0ex]{}}

\newcommand{\circg}{\overset{\raisebox{-0.6ex}[0.6ex][0ex]{$\scaleto{\circ}{0.6ex}$}}{g} \makebox[0ex]{}}
\newcommand{\subcircg}{\overset{\raisebox{-0.5ex}[0.5ex][0ex]{$\scaleto{\circ}{0.5ex}$}}{g} \makebox[0ex]{}}

\newcommand{\uf}{\underline{f} \makebox[0ex]{}}
\newcommand{\ufl}[1]{{}^{#1} \mkern-4mu \underline{f} \makebox[0ex]{}}
\newcommand{\uh}{\underline{h} \makebox[0ex]{}}

\newcommand{\dpartial}{\dot{\partial} \makebox[0ex]{}}
\newcommand{\fl}[1]{{}^{#1} \mkern-4mu f \makebox[0ex]{}}

\newcommand{\uvarepsilon}{\underline{\varepsilon} \makebox[0ex]{}}

\newcommand{\dslashd}{\dot{\raisebox{0ex}[1.5ex][0ex]{$\mathrm{d}$}} \mkern-8.9mu \raisebox{+0.2ex}[0ex][0ex]{$\scaleto{\boldsymbol{\slash}}{2ex}$} \mkern+2.3mu \makebox[0ex]{}}
\newcommand{\db}{\dot{b} \makebox[0ex]{}}
\newcommand{\dslashg}{\dot{g} \mkern-8.2mu \scaleto{\boldsymbol{\slash}}{1.6ex} \mkern+1mu \makebox[0ex]{}}
\newcommand{\circepsilon}{\overset{\raisebox{-0.6ex}[0.6ex][0ex]{$\scaleto{\circ}{0.6ex}$}}{\epsilon} \makebox[0ex]{}}

\newcommand{\dchi}{\dot{\chi} \makebox[0ex]{}}
\newcommand{\dtr}{\dot{\raisebox{0ex}[1.2ex][0ex]{$\mathrm{tr}$}} \makebox[0ex]{}}
\newcommand{\duchi}{\underline{\dot{\chi}} \makebox[0ex]{}}
\newcommand{\deta}{\dot{\eta} \makebox[0ex]{}}
\newcommand{\duomega}{\underline{\dot{\omega}} \makebox[0ex]{}}

\newcommand{\sym}{\mathrm{sym} \makebox[0ex]{}}

\newcommand{\circnabla}{\overset{\raisebox{0ex}[0.6ex][0ex]{$\scaleto{\circ}{0.6ex}$}}{\raisebox{0ex}[1.2ex][0ex]{$\nabla$}} \makebox[0ex]{}}
\newcommand{\circGamma}{\overset{\raisebox{0ex}[0.6ex][0ex]{$\scaleto{\circ}{0.6ex}$}}{\raisebox{0ex}[1.2ex][0ex]{$\Gamma$}} \makebox[0ex]{}}

\newcommand{\dalpha}{\dot{\alpha} \makebox[0ex]{}}
\newcommand{\dbeta}{\dot{\beta} \makebox[0ex]{}}
\newcommand{\drho}{\dot{\rho} \makebox[0ex]{}}
\newcommand{\dsigma}{\dot{\sigma} \makebox[0ex]{}}
\newcommand{\dubeta}{\underline{\dot{\beta}} \makebox[0ex]{}}
\newcommand{\dualpha}{\underline{\dot{\alpha}} \makebox[0ex]{}}

\newcommand{\dB}{\dot{B} \makebox[0ex]{}}

\newcommand{\ddslashg}{\ddot{g} \mkern-8.2mu \scaleto{\boldsymbol{\slash}}{1.6ex} \mkern+1mu \makebox[0ex]{}}
\newcommand{\subddslashg}{\ddot{g} \mkern-7.7mu \scaleto{\boldsymbol{\slash}}{1ex} \mkern+1mu \makebox[0ex]{}}  

\newcommand{\ddL}{\ddot{L} \makebox[0ex]{}}
\newcommand{\dduL}{\underline{\ddot{L}} \makebox[0ex]{}}
\newcommand{\dvarepsilon}{\dot{\varepsilon} \makebox[0ex]{}}
\newcommand{\ddslashd}{\ddot{\raisebox{0ex}[1.5ex][0ex]{$\mathrm{d}$}} \mkern-8.9mu \raisebox{+0.2ex}[0ex][0ex]{$\scaleto{\boldsymbol{\slash}}{2ex}$} \mkern+2.3mu \makebox[0ex]{}}

\newcommand{\dslashnabla}{\dot{\nabla} \mkern-13mu \raisebox{0.3ex}[0ex][0ex]{$\scaleto{\boldsymbol{\slash}}{1.7ex}$} \mkern+5mu \makebox[0ex]{}}
\newcommand{\dslashGamma}{\dot{\Gamma} \mkern-9.5mu \raisebox{0.4ex}[0ex][0ex]{$\scaleto{\boldsymbol{\slash}}{1.7ex}$} \mkern+1mu \makebox[0ex]{}}
\newcommand{\ddpartial}{\ddot{\partial} \makebox[0ex]{}}

\newcommand{\ddchi}{\ddot{\chi} \makebox[0ex]{}}
\newcommand{\dduchi}{\underline{\ddot{\chi}} \makebox[0ex]{}}
\newcommand{\ddeta}{\ddot{\eta} \makebox[0ex]{}}

\newcommand{\circdot}{\makebox[2ex]{$\circ\mkern-7.2mu\cdot$}}

\newcommand{\ddtr}{\ddot{\raisebox{0ex}[1.2ex][0ex]{$\mathrm{tr}$}} \makebox[0ex]{}}

\newcommand{\ddalpha}{\ddot{\alpha} \makebox[0ex]{}}
\newcommand{\ddbeta}{\ddot{\beta} \makebox[0ex]{}}
\newcommand{\ddrho}{\ddot{\rho} \makebox[0ex]{}}
\newcommand{\ddsigma}{\ddot{\sigma} \makebox[0ex]{}}
\newcommand{\ddubeta}{\underline{\ddot{\beta}} \makebox[0ex]{}}
\newcommand{\ddualpha}{\underline{\ddot{\alpha}} \makebox[0ex]{}}

\newcommand{\dslashepsilon}{\dot{\epsilon} \mkern-5.8mu  \raisebox{0.12ex}[0ex][0ex]{$\scaleto{\boldsymbol{\slash}}{1.4ex}$} \mkern+0.5mu \makebox[0ex]{}}
\newcommand{\ddslashepsilon}{\ddot{\epsilon} \mkern-5.8mu  \raisebox{0.12ex}[0ex][0ex]{$\scaleto{\boldsymbol{\slash}}{1.4ex}$} \mkern+0.5mu \makebox[0ex]{}}

\newcommand{\dd}[1]{\mathfrak{d}\{ #1 \} \makebox[0ex]{}}
\newcommand{\bdd}[1]{\delta\{#1\} \makebox[0ex]{}}

\newcommand{\circDelta}{\overset{\raisebox{0ex}[0.6ex][0ex]{$\scaleto{\circ}{0.6ex}$}}{\raisebox{0ex}[1.2ex][0ex]{$\Delta$}} \makebox[0ex]{}}

\newcommand{\balphal}[1]{{}^{#1}\mkern-1.5mu \bar{\alpha} \makebox[0ex]{}}
\newcommand{\bbetal}[1]{{}^{#1}\mkern-1.5mu \bar{\beta} \makebox[0ex]{}}
\newcommand{\bsigmal}[1]{{}^{#1}\mkern-1.5mu \bar{\sigma} \makebox[0ex]{}}
\newcommand{\brhol}[1]{{}^{#1}\mkern-1.5mu \bar{\rho} \makebox[0ex]{}}
\newcommand{\bubetal}[1]{{}^{#1}\mkern-1.5mu \bar{\underline{\beta}} \makebox[0ex]{}}
\newcommand{\bualphal}[1]{{}^{#1}\mkern-1.5mu \bar{\underline{\alpha}} \makebox[0ex]{}}

\newcommand{\circdiv}{\overset{\raisebox{0ex}[0.6ex][0ex]{$\scaleto{\circ}{0.6ex}$}}{\raisebox{0ex}[1.2ex][0ex]{$\mathrm{div}$}} \makebox[0ex]{}}

\newcommand{\bddsub}[2]{\delta_{#1}\{#2\}  \makebox[0ex]{}}

\newcommand{\deltasub}[1]{\delta_{#1} \makebox[0ex]{}}

\newcommand{\circcurl}{\overset{\raisebox{0ex}[0.6ex][0ex]{$\scaleto{\circ}{0.6ex}$}}{\raisebox{0ex}[1.2ex][0ex]{$\mathrm{curl}$}} \makebox[0ex]{}}

\newcommand{\calM}{\mathcal{M} \makebox[0ex]{}}

\newcommand{\bfg}{\mathbf{g} \makebox[0ex]{}}
\newcommand{\bfk}{\mathbf{k} \makebox[0ex]{}}

\newcommand{\Ima}{\mathrm{im} \makebox[0ex]{}}

\newcommand{\calS}{\mathcal{S} \makebox[0ex]{}}
\newcommand{\bfs}{\mathbf{s} \makebox[0ex]{}}
\newcommand{\udelta}{\underline{\delta} \makebox[0ex]{}}
\newcommand{\rmW}{\mathrm{W} \makebox[0ex]{}}
\newcommand{\Vvert}{\vert \mkern-2mu \vert \mkern-2mu \vert \makebox[0ex]{}}

\title{\textsc{Linearised Perturbation of Constant Mass Aspect Function Foliation in Schwarzschild Black Hole Spacetime}}
\author{Pengyu Le}
\newcommand{\Address}{{
  \bigskip
  \footnotesize
  \textsc{Yanqi Lake Beijing Institute of Mathematical Sciences and Applications, Beijing, China}
  
  \textit{E-mail address}: \texttt{pengyu.le@bimsa.cn}
}}
\date{}

\begin{document}

\maketitle

\begin{abstract}
We study the linearised perturbation of the constant mass aspect function foliation in a Schwarzschild black hole spacetime. In particular, we investigate the linearised perturbation of the asymptotic geometry of the foliation at null infinity. We show that there is a $4$-dimensional linear space for the linearised perturbation of the initial leaf inside the event horizon, corresponding to which the linearised perturbation of the asymptotic geometry of the foliation at null infinity is preserved to be round. For such a linearised perturbation of the initial leaf in this $4$-dimensional linear space, we calculate the corresponding linearised perturbation of the energy-momentum vector and the Bondi mass at null infinity. We show that the linearised perturbations of the Bondi energy and the Bondi mass both vanish, and every possible linearised perturbation of the linear momentum can be achieved by a linearised perturbation of the initial leaf in the $4$-dimensional linear space.
\end{abstract}

\tableofcontents

\section{Introduction}\label{sec 1}
The theme of this paper is the study of the linearised perturbation of the constant mass aspect function foliation in a Schwarzschild black hole spacetime. 

The concept of a constant mass aspect function foliation (see definitions \ref{def 2.3}, \ref{def 2.6}) was introduced in \cite{C1991} \cite{CK1993} for a spacelike hypersurface, and in \cite{C2003} \cite{S2008} for a null hypersurface. It plays an important role in the monumental proof of the global nonlinear stability of the Minkowski spacetime in \cite{CK1993}. It was also found that the Hawking mass satisfies an elegant variation formula (see formula \eqref{eqn 2.1}) along a constant mass aspect function foliation in \cite{C2003} \cite{S2008}, similar to the Geroch monotonicity formula of the Hawking mass along an inverse mean curvature flow in \cite{G1973} \cite{HI2001}.

Lead by the variation formula of the Hawking mass along a constant mass aspect function foliation, \cite{S2008} initiated the project using the constant mass aspect function foliation to prove the null Penrose inequality on a null hypersurface close to the spherically symmetric null hypersurface in a Schwarzschild spacetime. The famous Penrose inequality was proposed by Penrose in \cite{P1973} which conjectures the relation between the area of an outmost marginally trapped surface and the total mass of the spacetime. The null Penrose inequality is such a type of inequality on a null hypersurface, relating the area of an outmost marginally trapped surface with the Bondi mass of the null hypersurface at null infinity.

Although it didnot prove the null Penrose inequality in this case, \cite{S2008} pointed out that the obstacle towards a proof lies in the asymptotic geometry of the foliation at null infinity. The asymptotic geometry of the foliation at null infinity plays an important role in the relation between the limit of the Hawking mass along the foliation and the Bondi mass at null infinity, see subsections \ref{subsec 2.6}, \ref{subsec 2.7}.

Christodoulou and Sauter proposed the project to investigate the change of the asymptotic geometry of the constant mass aspect function foliation at null infinity by varying the foliation, and find a null hypersurface, in which the constant mass aspect function foliation starting from a marginally trapped surface is asymptotically round at null infinity. Carrying out their proposal will resolve the above mentioned difficulty arisen from the asymptotic geometry of the foliation at null infinity. 

As the first step to accomplish the proposal of Christodoulou and Sauter, we investigate the linearised perturbation of the constant mass aspect function foliation in a Schwarzschild spacetime. We sketch the main results obtained in our investigation.
\begin{theorem}[Sketch of theorem \ref{thm 11.1}]\label{thm 1.1}
For the spherically symmetric constant mass aspect function foliation in a spherically symmetric null hypersurface in a Schwarzschild spacetime, let $\Sigma_{0,0}$ be its leaf in the event horizon. There exists a $4$-dimensional linear space $V$ of the linearised perturbation of $\Sigma_{0,0}$ inside the event horizon, such that the resulting constant mass aspect function foliation is preserved to be asymptotically round at null infinity on the linearised level. This $4$-dimensional linear space $V$ can be characterised by the spherical harmonics of degrees $0$ and $1$.
\end{theorem}

Moreover, we investigate the change of the energy-momentum vector relative to the linearised perturbation of the constant mass aspect function foliation.
\begin{theorem}[Sketch of theorem \ref{thm 12.1}]\label{thm 1.2}
The linearised perturbation of the Bondi energy relative to the linearised perturbation of $\Sigma_{0,0}$ in the $4$-dimensional linear space $V$ in theorem \ref{thm 1.1} vanishes. The linear map from the $4$-dimensional linear space $V$ to the $3$-dimensional space of the linearised perturbation of the linear momentum vector is surjective, with the $1$-dimensional kernel consisting of constant functions.
\end{theorem}

By theorem \ref{thm 1.2}, we obtain that given any linearised perturbation of the linear momentum, there exists some linearised perturbation of $\Sigma_{0,0}$ in the $4$-dimensional space $V$ such that its corresponding linearised perturbation of the linear momentum is the given one.

\section{Constant mass aspect function foliation}
In this section, we define the mass aspect function on a spacelike surface, then introduce the constant mass aspect function foliation on a null hypersurface. We introduce a system of equations to study the geometry of a constant mass aspect function foliation.

\subsection{Mass aspect function}
\begin{figure}[H]
\begin{center}
\begin{tikzpicture}[scale=1]
\draw (-1.5,0)
to [out=80, in=180] (0.1, 0.9)
to [out=0, in=90] (1.4, 0) node[right] {$\Sigma$}
to [out=-90, in=10] (0, -1.1) 
to [out=-170, in=-100] (-1.5,0);
\draw[blue,->] (0.1, 0.9) -- (0.6,1.4) node[above]{$L$};
\draw[blue,->] (0.1, 0.9) -- (0.6,0.4) node[below]{$\uL$};
\node[right] at (2,0.9) {$g(L,\uL)=2$};
\end{tikzpicture}
\caption{A conjugate null frame $\{L, \uL\}$.}
\label{fig 1}
\end{center}
\end{figure}
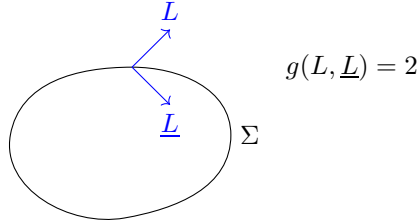

Let $(M,g)$ be a time-oriented $4$-dimensional Lorentzian manifold and $\Sigma$ be a spacelike surface in $(M,g)$. We introduce the notion of a conjugate null frame of $\Sigma$: $\{L,\uL\}$ is called a conjugate null frame of $\Sigma$ if $L, \uL$ are null normal vectors of $\Sigma$ with the condition
\begin{align*}
g(L,\uL)= 2.
\footnotemark
\end{align*}
\footnotetext{The sign here is opposite to the ones in many literatures, for example \cite{C2009}. It is obvious that $g(L, \uL) =2$ here implies that $L, \uL$ have different causal directions.}See figure \ref{fig 1}. Without loss of generality, we assume that $L$ is future directed throughout this paper.
For any positive number $a$, the transformation $\{L, \uL\} \mapsto \{ a L, a^{-1} \uL \}$ transforms a conjugate null frame to another one. We can define the following geometric quantities on $\Sigma$ relative to a conjugate null frame $\{L,\uL\}$.
\begin{definition}[Connection coefficients]\label{def 2.1}
We define the second fundamental forms of $\Sigma$ relative to $L, \uL$: let $X,Y$ be tangential to $\Sigma$,
\begin{align*}
\chi(X,Y)= g(\nabla_X L, Y), 
\quad
\uchi(X,Y) = g(\nabla_X \uL,Y),
\end{align*}
where $\nabla$ is the covariant derivative of $(M,g)$. We also define the torsion of the conjugate null frame $\{ L, \uL \}$ of $\Sigma$,
\begin{align*}
\eta(X) = \frac{1}{2} g ( \nabla_X \uL, L).
\end{align*}
Let $\slashg$ be the intrinsic metric of $\Sigma$. We decompose $\chi$, $\uchi$ into their traces and trace-free parts: let $\tr$ be the trace operator relative to $\slashg$,
\begin{align*}
\chi = \hatchi + \frac{1}{2} \tr \chi\, \slashg,
\quad
\uchi = \hatuchi + \frac{1}{2} \tr \uchi\, \slashg.
\end{align*}
$\tr \chi$, $\tr \uchi$ are called null expansions and $\hatchi$, $\hatuchi$ are called shears of $\Sigma$ relative to the conjugate null frame $\{L, \uL\}$.
\end{definition}
With the help of the concept of null expansions, we give the definitions of a trapped surface and a marginally trapped surface here.
\begin{definition}\label{def 2.2}
We call a spacelike surface $\Sigma$ in $(M,g)$ a trapped surface if
\begin{align*}
\tr \chi <0,
\quad
- \tr \uchi <0.\footnotemark
\end{align*}
\footnotetext{The negative sign here is due to the fact that $\uL$ is past-directed.}We call $\Sigma$ a margianally trapped surface if one of $\tr \chi$ and $-\tr \uchi$ vanishes identically while the other one is non-positive.
\end{definition}

Now we define the mass aspect function.
\begin{definition}[Mass aspect function]\label{def 2.3}
The mass aspect function $\mu$ of a spacelike surface $\Sigma$ relative to a conjugate null frame $\{L ,\uL\}$ is
\begin{align*}
\mu = K -  \frac{1}{4} \tr \chi \tr \uchi - \slashdiv \eta,
\end{align*}
where $K$ is the Gauss curvature of $(\Sigma, \slashg)$ and $\slashdiv$ is the divergence operator on $(\Sigma,\slashg)$.
\end{definition}
Note that the dependence of $\{L, \uL\}$ in the mass aspect function solely lies in the last term $\slashdiv \eta$. The mass aspect function function is closely related to the Hawking mass. If the surface $\Sigma$ has the topology of a sphere, then we have that by the Gauss-Bonnet theorem,
\begin{align*}
m_H ( \Sigma) 
= 
\frac{r(\Sigma)}{8\pi} \int_\Sigma \mu\, \dvol_{\subslashg} 
= 
\frac{r(\Sigma)}{2} \Big( 1 - \frac{1}{16\pi} \int_\Sigma \tr \chi \tr \uchi\, \dvol_{\subslashg} \Big),
\end{align*}
where $m_H(\Sigma)$ is the Hawking mass of $\Sigma$, $r(\Sigma)$ is the area radius of $\Sigma$ defined by $|\Sigma| = 4\pi r^2(\Sigma)$, and $\dvol_{\subslashg}$ is the area element of $(\Sigma,\slashg)$. We see that if $\Sigma$ is marginally trapped, then $m_H(\Sigma) = \frac{r(\Sigma)}{2} = \sqrt{\frac{|\Sigma|}{16\pi}}$. 

The concept of the Hawking mass was introduce by Hawking in \cite{H1968}. The mass aspect function was introduced in \cite{C1991} by Christodoulou and applied in the monumental proof of global nonlinear stability of Minkowski spacetime by Christodoulou and Klainerman in \cite{CK1993}.

\subsection{Variation of Hawking mass along a foliation on a null hypersurface}\label{subsec 2.2}

\begin{figure}[H]
\begin{center}
\begin{tikzpicture}[scale=1]
\draw[dashed] (-1.5,0)
node[below right]{\footnotesize $\Sigma_u$}
to [out=50, in=180] (0.1, 0.6)
to [out=0, in=130] (1.4, 0); 
\draw (1.4,0) 
to [out=-50, in=10] (0, -1.1) 
to [out=-170, in=-130] (-1.5,0); 
\draw[dashed] (-2,-1) to [out=70, in=180] (0.1, -0.8) to [out=0, in=150] (2, -1); 
\draw (2,-1) to [out=-30, in=10] (-0.2, -2.2)to [out=-170, in=-110] (-2,-1);  
\draw[blue] (-0.9,0.7) to [out=-110, in=50] (-1.5,0);
\draw[blue] (-1.5,0) to [out=-130, in=70] (-2,-1) to [out=-110,in=65] (-3,-2);
\draw[blue] (1.2,0.5) to [out=-60,in=130] (1.4,0) to [out=-50, in=150] (2,-1) to [out=-30,in=115] (3,-2) node[above right] {\footnotesize $\ucalH$}; 
\draw[blue] (0.1,0) to [out=-80, in=80] (0,-1.1);
\draw[blue] (0,-1.1) to [out=-100, in=60] (-0.2,-3); 
\node[above right,blue] at (0,-2.8) {\footnotesize a null geodesic in $\ucalH$}; 
\draw[blue,->] (0, -1.1)  -- (0.5,-0.6) node[above]{\footnotesize $L'^u$};
\draw[blue,->] (0, -1.1)  -- (-0.21,-2.2) node[below]{\footnotesize $\uL^u$};
\end{tikzpicture}
\end{center}
\caption{The conjugate null frame $\{L'^u, \uL^u\}$ relative to $\{ \Sigma_u \}$.}
\label{fig 2}
\end{figure}
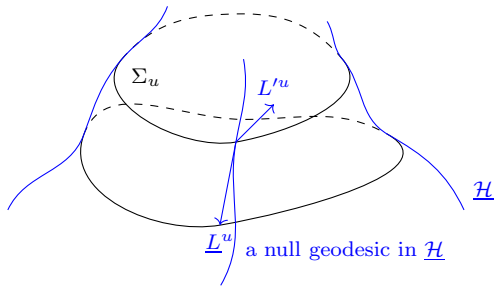

Let $\ucalH$ be a null hypersurface in $(M,g)$, and $\{\Sigma_u\}$ be a foliation by spacelike surfaces of $\ucalH$. Assume that the parameter $u$ increases toward the past direction. We can define the conjugate null frame $\{L'^u, \uL^u\}$\footnote{We add a prime in $L'^u$ of the conjugate null frame, since the simpler notation $L$ will be reserved for the outgoing null vector field in the double null coordinate system in a Schwarzschild spacetime.} on each $\Sigma_u$ relative to the foliation $\{\Sigma_u\}$ as follows. Let $\uL^u$ be the tangent null vector of $\ucalH$ with the condition
\begin{align*}
\uL^u u=1,
\end{align*}
then define $L'^u$ to be the other normal null vector of $\Sigma_u$ such that $g(L'^u, \uL^u)=2$, thus $\{L'^u, \uL^u\}$ is a conjugate null frame of $\Sigma_u$. See figure \ref{fig 2}.
Then we can define the corresponding connection coefficients and the mass aspect function $\mul{u}$ of $\Sigma_u$ relative to $\{ L'^u, \uL^u\}$ by the constructions in definitions \ref{def 2.1} and \ref{def 2.3}.\footnote{We use either a subscript or a superscript $u$ (on the left or right) to indicate the corresponding notation is defined for $\Sigma_u$. In the following, we omit the subscript or superscript $u$ for the sake of brevity when no ambiguity arises.\label{footnote 4}} Furthermore, we can define the acceleration of the tangent null vector field $\uL^u$ on $\ucalH$ as follows.
\begin{definition}\label{def 2.4}
Define the acceleration $\uomegal{u}$ of the tangent null vector field $\uL^u$ on $\ucalH$ by
\begin{align*}
\uomegal{u} = \frac{1}{4} g( \nabla_{\uL^u} \uL^u, L'^u).
\end{align*}
Equivalently, the acceleration $\uomegal{u}$ is given by
\begin{align*}
\nabla_{\uL^u} \uL^u = 2\, \uomegal{u} \uL^u.
\end{align*}
\end{definition}
We have the following nice proposition for the variation of the Hawking mass $m_H(\Sigma_u)$ in a vacuum spacetime $(M,g)$.
\begin{proposition}\label{pro 2.5}
Suppose that $(M,g)$ is a vacuum spacetime, i.e. $\mathrm{Ric}=0$, then the variation of the Hawking mass along the foliation $\{\Sigma_u\}$ on $\ucalH$ is
\begin{align}
\begin{aligned}
\frac{\ed}{\ed s} m_H(\Sigma_u) 
&= 
\frac{r_{s}}{32\pi} \int_{\Sigma_u}
\big( \tr \chil{u}'\, \vert \hatuchil{u} \vert^2 + 2 \tr \uchil{u}\, \vert \etal{u} \vert^2 \big) \dvol_u
\\
&\phantom{=}
- 
\frac{r_{s}}{16\pi} \int_{\Sigma_u}
\big( \tr \uchil{u} - \overline{\tr \uchil{u}} \big) \big( \mul{u} - \overline{\mul{u}} \big) \dvol_u,
\end{aligned}
\label{eqn 2.1}
\end{align}
where $r_u$ is the area radius of $\Sigma_u$, $\vert \cdot \vert$ is the norm relative to the metric $\slashgl{u}$ on $\Sigma_u$, $\dvol_u$ is the area element of $(\Sigma_u, \slashgl{u})$, and the overline $\overline{\makebox[1.5ex]{$\cdot$}}$ is the mean value on $(\Sigma_u, \slashgl{u})$ that
$\overline{f} = \frac{\int_{\Sigma_u} f\, \dvol_u}{|\Sigma_u|}$.
\end{proposition}
We refer to \cite{S2008} for the derivation of the above proposition. In \cite{H1968}, Hawking derived the variation of the Hawking mass employing the Newman-Penrose formalism in a general spacetime (not necessary being vacuum). In particular he considered the special foliation where $u$ is a luminosity parameter, which is equivalent to $\tr\uchil{u} = \overline{\tr \uchil{u}} = \frac{2}{r_u}$. This neat formula \eqref{eqn 2.1} in a vacuum spacetime is due to Christodoulou in \cite{C2008} and his unpublished lecture note \cite{C2003}. Note that the mass aspect function appears naturally in formula \eqref{eqn 2.1}. 

An important corollary of formula \eqref{eqn 2.1} is the monotonicity of the Hawking mass. If $\tr \chil{u}'$ and $\tr \uchil{u}$ are both non-negative, then $\tr \uchil{u} = \overline{\tr \uchil{u}}$ or $\mul{u} = \overline{\mul{u}}$ implies that 
\begin{align*}
\frac{\ed}{\ed s} m_H(\Sigma_u)  \geq 0,
\end{align*}
i.e. the Hawking mass is monotonically nondecreasing along the foliation. 

Such kind of monotonicity proposition of the Hawking mass was observed by Geroch in \cite{G1973} for an inverse mean curvature flow in a $3$-dimensional Riemannian manifold with non-negative scalar curvature, and he proposed using the inverse mean curvature flow to prove the positive energy theorem, then still a conjecture until fully proved by Schoen and Yau in \cite{SY1979}. Geroch's proposal was modified by Jang and Wald in \cite{JW1977} to prove the Penrose inequality in the above context. This proposal was fully accomplished in the remarkable proof of the Riemannian Penrose inequality by Huisken and Ilmanen in \cite{HI2001} which overcomes the difficulty that the inverse mean curvature flow would develop singularities. Another different proof of the Riemannian Penrose inequality by the conformal flow method was given by Bray in \cite{Br2001} independently at the same time.

\subsection{Constant mass aspect function foliation}\label{subsec 2.3}
The variation formula \eqref{eqn 2.1} of the Hawking mass gives rise to two interesting kinds of foliations:
\begin{enumerate}[label=\arabic*.]
\item constant null expansion foliation (foliation by a luminosity parameter in the terminology of \cite{H1968}): $\tr \uchil{u}= \overline{\tr \uchil{u}} = \frac{2}{r_u}$. Note that $\frac{1}{r_u^2}\dvol_u = \mathrm{const.}$ along the flow generated by the tangent null vector field $\uL^u$. This property is shared by the inverse mean curvature flow, thus we can view the constant null expansion foliation as the analogue of the inverse mean curvature flow in a null hypersurface.

\item constant mass aspect function foliation: $\mul{u} = \overline{\mul{u}}$. We shall elaborate on this foliation later.
\end{enumerate}
Since the constant mass aspect function foliation is the main topic of this paper, we give its definition again here.
\begin{definition}\label{def 2.6}
Let $\ucalH$ be a null hypersurface in $(M,g)$, and $\{\Sigma_u\}$ be a foliation by spacelike surfaces in $\ucalH$. We call $\{\Sigma_u\}$ a constant mass aspect function foliation if the mass aspect function $\mul{u}$ relative to the conjugate null frame $\{ L'^u, \uL^u \}$ is constant along each leaf $\Sigma_u$, i.e. $ \mul{u} = \overline{\mul{u}} = \frac{2 m_H (\Sigma_u)}{r_u^3}$ for every $u$.
\end{definition}
Note that relabelling the foliation preserves the fact that the mass aspect function being constant. We often choose the area radius to remove this freedom of relabelling by requiring that $r_{u_0+u} = r_{u_0}+u$, where the freedom left in the labelling of the foliation is just the difference by a constant. Throughout this paper, we always assume that the foliation has the labelling $r_{u_0+u} = r_{u_0}+u$ unless specifying other kinds of labelling.

Like the constant null expansion foliation has the inverse mean curvature flow as its analogue in a spacelike hypersurface, the constant mass aspect function foliation also has such an analogue. We sketch its analogue before proceeding the discussion of the constant mass aspect function foliation.
\begin{enumerate}[label=\ref{subsec 2.3}.\alph*.,leftmargin=.45in]
\item
Let $H$ be a spacelike hypersurface in $(M,g)$, and $\{\Sigma_u\}$ be a foliation by surfaces (being spacelike automatically) in $H$. 
\item
Define $N_u$ as the normal vector field of $\Sigma_u$ in $H$ such that $N_u u=1$. Let $n_u$ be the unit normal vector field of $\Sigma_u$ in $H$ such that $n_u u>0$. 
\item
There exists a positive function $a_u$ on $\Sigma_u$ such that $N_u= a_u n_u$. We call $a_u$ the lapse function of the foliation $\Sigma_u$. See figure \ref{fig 3}.
\item
Let $K_u$ be the Gauss curvature of $(\Sigma_u, \slashg)$, $\theta$ be the second fundamental form of $\Sigma_u$ relative to $n_u$ in $H$, and $\tr \theta$ be the corresponding mean curvature.
\item
Define the function $\nu_u$ on each leaf $\Sigma_u$ by $\nu_u = K_u -\frac{1}{4} ( \tr \theta)^2 $, then we require the lapse function satisfies the equation
\begin{align*}
\slashDelta \log a_u = \nu_u - \overline{\nu_u}.
\end{align*} 
Then the above construction determines the foliation up to the freedom in relabelling the leaves.
\end{enumerate}
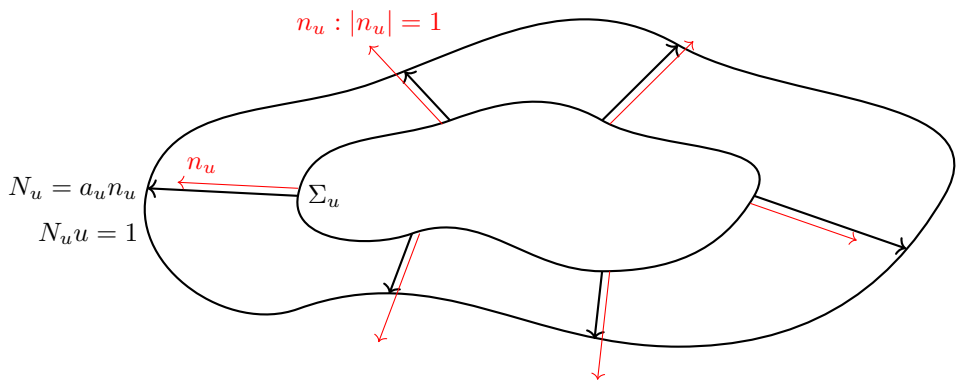
\begin{figure}[H]
\begin{center}
\begin{tikzpicture}
\draw[thick] (-3,0) 
to [out=80,in=-160]
(-1,1)
to [out=20,in=150]
(1,1)
to [out=-30,in=60]
(3,0)
to [out=-120,in=0]
(1,-1)
to [out=-180,in=20]
(-1.5,-0.5)
to [out=-160,in=-100]
(-3,0)
node[right] {$\Sigma_u$};
\draw[thick] (-5,0) 
to [out=80,in=-160]
(-2,1.5)
to [out=20,in=150]
(2,2)
to [out=-30,in=60]
(5.5,0)
to [out=-120,in=0]
(2,-2)
to [out=-180,in=20]
(-3,-1.5)
to [out=-160,in=-100]
(-5,0);
\draw[thick,->]
(-3,0) -- (-4.98,0.1)
node[left] {$N_u=a_u n_u$} ;
\node[left] at (-4.98,-0.5) {$N_u u=1$};
\draw[thick,->]
(-1,1) -- (-1.6,1.65);
\draw[thick,->]
(1,1) -- (2,2);
\draw[thick,->]
(3,0) -- (5,-0.7);
\draw[thick,->]
(1,-1) -- (0.9,-1.9);
\draw[thick,->]
(-1.5,-0.5) -- (-1.8,-1.3);
\draw[->,red]
(-3,0+0.1) -- (-4.98*0.8+3*0.8 -3, 0.1*0.8-0*0.8 +0 +0.1)
node[above right] {$n_u$};
\draw[->,red]
(-1-0.1,1-0.05) -- (-1.6*1.6+1*1.6 -1 -0.1, 1.65*1.6-1*1.6 +1 -0.05)
node[above] {$n_u: |n_u|=1$};
\draw[->,red]
(1+0.1,1-0.05) -- (2*1.1-1*1.1 +1 +0.1, 2*1.1-1*1.1 +1 -0.05);
\draw[->,red]
(3-0.05,0-0.1) -- (5*0.7-3*0.7 +3 -0.05, -0.7*0.7-0*0.7 +0 -0.1);
\draw[->,red]
(1+0.1,-1) -- (0.9*1.6-1*1.6 +1+0.1, -1.9*1.6+1*1.6 -1);
\draw[->,red]
(-1.5+0.1,-0.5) -- ( -1.8*1.8+1.5*1.8 -1.5 +0.1, -1.3*1.8+0.5*1.8 -0.5);
\end{tikzpicture}
\end{center}
\caption{Lapse function $a_u$: $N_u=a_u n_u$.}
\label{fig 3}
\end{figure}
The above analogue of the constant mass aspect function foliation was introduced in \cite{CK1993}, and plays an important role in the proof of global nonlinear stability of the Minkowski spacetime. We recommend the exposition on this topic in \cite{C2008} to interested readers.

We return to the constant mass aspect function foliation. Because of the monotonicity property of the Hawking mass along the two foliations, they can be potentially applied to prove the Penrose inequality on a null hypersurface, like the application of the inverse mean curvature flow to the Riemannian Penrose inequality.

Motived by this, Sauter in his thesis \cite{S2008} developed the local existence theory of two foliations, and proved the global existence theorems for both foliations on a null hypersurface close to the spherically symmetric null hypersurface in a Schwarzschild black hole spacetime. His work pointed out the importance of the asymptotic geometry of the foliation at null infinity if one want to apply the two foliations to prove the Penrose inequality in this scenario. We shall discuss this further in subsections \ref{subsec 2.6} and \ref{subsec 2.7}.

Because of the above, Christodoulou and Sauter proposed the project to study the change of the asymptotic geometry of two foliations when varying the null hypersurface. In particular, they asked whether one can find a null hypersurface in which at least one of the two foliations starting from a marginally trapped surface is asymptotically round at null infinity. 

Naturally as the first step to carry out their proposal, one shall study the linearised perturbation of the asymptotic geometry of the foliation in a Schwarzschild spacetime. In this paper, we take this task for the constant mass aspect function foliation.

\subsection{Formulation of construction of the foliation as an inverse lapse problem}\label{subsec 2.4}
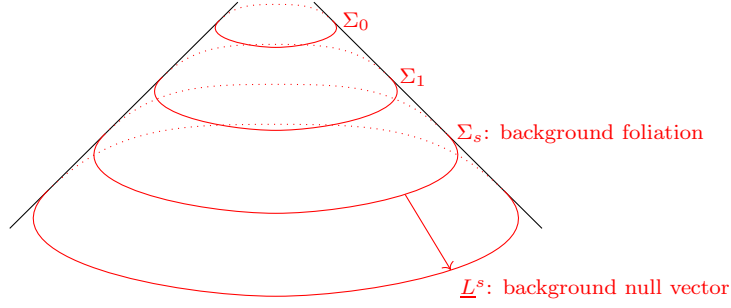
\begin{figure}[H]
\begin{center}
\begin{tikzpicture}
\draw (-0.5,-0.5) -- (-3.5,-3.5)  (0.5,-0.5) -- (3.5,-3.5); 
\draw[opacity=1,dotted,red] (-0.75, -0.75) 
to [out=45, in=180]  (0,-0.53) 
to [out=0, in=135]  (0.75,-0.75) 
node[right] { \footnotesize $\Sigma_0$};
\draw[opacity=1, red] (0.75,-0.75)
to [out=-45, in=0]  (0,-1.1)
to [out=180, in=-135] (-0.75,-0.75); 
\draw[opacity=1,dotted,red] (-1.5, -1.5) 
to [out=45, in=180]   (0,-1.06) 
to [out=0, in=135]  (1.5,-1.5)
node[right] { \footnotesize $\Sigma_1$};
\draw[opacity=1,red] (1.5,-1.5)
to [out=-45, in=0]  (0,-2.2)
to [out=180, in=-135] (-1.5,-1.5); 
\draw[opacity=1,dotted,red] (-2.25, -2.25) 
to [out=45, in=180]   (0,-1.59) 
to [out=0, in=135]  (2.25,-2.25)
node[right] { \footnotesize  $\Sigma_s$: background foliation};
\draw[opacity=1,red] (2.25,-2.25)
to [out=-45, in=0]  (0,-3.3)
to [out=180, in=-135] (-2.25,-2.25); 
\draw[opacity=1,dotted,red] (-3, -3) 
to [out=45, in=180]   (0,-2.12) 
to [out=0, in=135]  (3,-3);
\draw[opacity=1,red] (3,-3)
to [out=-45, in=0]  (0,-4.4)
to [out=180, in=-135] (-3,-3); 
\draw[->,opacity=1, red] (1.7,-3.05) -- (2.3,-4.05)
node[below right] {\footnotesize $\uL^s$: background null vector};
\end{tikzpicture}
\end{center}
\caption{Background foliation $\Sigma_s$ and background null vector $\uL^s$.}
\label{fig 4}
\end{figure}

In this subsection, we review the formulation of the construction of a constant mass aspect function foliation as an inverse lapse problem, following \cite{S2008}. Suppose that $\{\Sigma_s\}$ is a background foliation on $\ucalH$. Let $\uL^s$ be the corresponding tangent null vector on $\ucalH$ that $\uL^s s=1$. See figure \ref{fig 4}. Introduce a coordinate system $\{ s, \theta^1, \theta^2 \}$ of $\ucalH$ that $\uL^s \theta^1 =\uL^s \theta^2 =0$, which is equivalent to that $\uL^s = \partial_s$. For the sake of brevity, we use $\vartheta$ to denote the coordinate $(\theta^1, \theta^2)$.

We formulate the problem of constructing a constant mass aspect function foliation as follows. Given a spacelike surface $\bSigma_{u_0}$\footnote{We use the overbar to distiguish $\{\bSigma_u\}$ with the background foliation $\{\Sigma_s\}$ more clearly. Later, we will use the overbar to indicate the quantities associated with $\{\bSigma_u\}$.} in $\ucalH$, find the constant mass aspect function foliation $\{\bSigma_u\}$ of $\ucalH$ containing $\Sigma_{u_0}$. Moreover we require that $\br_{u_0 + u} = \br_{u_0} + u$.

The above construction problem can be reformulated as an inverse lapse problem as follows.
\begin{enumerate}[label=\ref{subsec 2.4}.\alph*.,leftmargin=.45in]
\item Parameterise a spacelike surface $\bSigma_u$ as the graph of a function $f(u, \vartheta)$ of $s$ over the $\vartheta$ domain in the coordinate system $\{s, \vartheta\}$.

\item Denote the metric on $\bSigma_u$ by $\bslashgl{u}$. Let $\btr$ be the trace operator relative to $\bslashgl{u}$. Let $\bslashdiv$ and $\bslashDelta$ be the divergence and Laplace operators on $(\bSigma_u, \bslashgl{u})$.

\item Construct the conjugate null frame $\{ \dL'^u, \duL^u\}$ of $\Sigma_u$, where we choose that $\duL^u = \uL^s$. Note that $\dL'^{u} \neq L'^s$ in general.\footnote{Similar to the overbar notation, we use an overdot to indicate that the quantity being associated with $\bSigma_u$ instead of the background foliation $\Sigma_s$.} This conjugate null frame is determined by one leaf $\bSigma_u$ only, not depending on the whole foliation $\{\bSigma_u\}$.

\item We can obtain the mass aspect function $\dmul{u}$ of $\bSigma_u$ relative to $\{\dL'^u, \duL^u\}$. Let $\btr \dchil{u}'$, $\btr \duchil{u}$ be the null expansions of $\bSigma_u$ relative to $\{\dL'^u, \duL^u\}$, and $\detal{u}$ be the torsion of $\{\dL'^u, \duL^u\}$, then
\begin{align*}
\dmul{u} = \bKl{u} - \frac{1}{4} \btr \dchil{u}'\, \btr \duchil{u} - \bslashdiv \detal{u}.
\end{align*}
Since $\{\dL'^u, \duL^u\}$ depends only on $\bSigma_u$, so is the mass aspect function $\dmul{u}$.

\item Construct the corresponding conjugate null frame $\{\bL'^u, \buL^u\}$ relative to the whole foliation $\{\bSigma_u\}$. Since $\buL^u$ and $\duL^u$ are collinear, we introduce the following function $\bal{u}$ by
\begin{align*}
\buL^u = \bal{u} \duL^u = \bal{u} \uL^s.
\end{align*}
$\bal{u}$ is called the lapse function. See figure \ref{fig 5}.
\begin{figure}[H]
\begin{center}
\begin{tikzpicture}
\draw (-0.5,-0.5) -- (-3.5,-3.5)  (0.5,-0.5) -- (3.5,-3.5); 
\draw[opacity=0.8,dotted,red] (-0.75, -0.75) 
to [out=45, in=180]  (0,-0.53) 
to [out=0, in=135]  (0.75,-0.75) 
node[right] { \footnotesize $\Sigma_0$};
\draw[opacity=0.8,red] (0.75,-0.75)
to [out=-45, in=0]  (0,-1.1)
to [out=180, in=-135] (-0.75,-0.75); 
\draw[opacity=0.8,dotted,red] (-1.5, -1.5) 
to [out=45, in=180]   (0,-1.06) 
to [out=0, in=135]  (1.5,-1.5)
node[right] { \footnotesize $\Sigma_1$};
\draw[opacity=0.8,red] (1.5,-1.5)
to [out=-45, in=0]  (0,-2.2)
to [out=180, in=-135] (-1.5,-1.5); 
\draw[opacity=0.8,dotted,red] (-2.25, -2.25) 
to [out=45, in=180]   (0,-1.59) 
to [out=0, in=135]  (2.25,-2.25)
node[right] { \footnotesize  $\Sigma_s$: background foliation};
\draw[opacity=0.8,red] (2.25,-2.25)
to [out=-45, in=0]  (0,-3.3)
to [out=180, in=-135] (-2.25,-2.25); 
\draw[opacity=0.8,dotted,red] (-3, -3) 
to [out=45, in=180]   (0,-2.12) 
to [out=0, in=135]  (3,-3);
\draw[opacity=0.8,red] (3,-3)
to [out=-45, in=0]  (0,-4.4)
to [out=180, in=-135] (-3,-3); 
\draw[thick,dashed] (-2.1, -2.1) node[above left] { \footnotesize $\bSigma_u$}
to [out=45, in=-165]  (-0.2,-2) 
to [out=15, in=135]  (2.4,-2.4);
\draw[thick] (2.4,-2.4)
to [out=-45, in=-20]  (0.3,-3)
to [out=160, in=-135] (-2.1,-2.1); 
\draw[thick,dashed] (-3, -3)
to [out=45, in=-165]  (-0.4,-2.8) 
to [out=15, in=135]  (3.15,-3.15);
\draw[thick] (3.15,-3.15)
to [out=-45, in=-30]  (-1.5,-4)
to [out=150, in=-135] (-3,-3); 
\draw[->,opacity=0.8,red] (1.7+0.1,-3.05+0.02) -- (2.3+0.1,-4.05+0.02)
node[above right] {\footnotesize $\duL^u = \uL^s$: background null vector};
\draw[thick,->] (1.7,-3.05) -- (2.3*1.25-1.7*1.25 +1.7,-4.05*1.25+3.05*1.25 -3.05)
node[below right] {\footnotesize $\buL^u = \bal{u} \uL^s$, $\buL^u u=1$};
\end{tikzpicture}
\end{center}
\caption{Lapse function $\bal{u}$: $\buL^u = \bal{u} \uL^s$.}
\label{fig 5}
\end{figure}

The parameterisation $f(u,\vartheta)$ of $\{ \bSigma_u \}$ satisfies the equation
\begin{align*}
\frac{\ed}{\ed u} f(u, \vartheta) = \bal{u} (\vartheta).
\end{align*}

$\{\bL'^u, \buL^u\}$ and $\{ \dL'^u, \duL^u\}$ are related by the following transformation
\begin{align*}
\dL'^u \rightarrow \bL'^u = \bal{u}^{-1} \dL'^u,
\quad
\duL^u \rightarrow \buL^u = \bal{u} \duL^u.
\end{align*}

\item Consider the mass aspect function $\bmul{u}$ of $\bSigma_u$ relative to $\{\bL'^u, \buL^u\}$. Introduce the null expansions $\btr \bchil{u}'$, $\btr \buchil{u}$ of $\bSigma_u$ relative to $\{ \bL'^u, \buL^u \}$, and the torsion of $\btal{u}$ of $\{ \bL'^u, \buL^u\}$. Then
\begin{align*}
\bmul{u} = \bKl{u} - \frac{1}{4} \btr \bchil{u}'\, \btr \buchil{u} - \bslashdiv \btal{u}.
\end{align*}

\item 
The pairs of null expansions $\btr \bchil{u}'$, $\btr \buchil{u}$ and $\btr \dchil{u}'$, $\btr \duchil{u}$ are related by
\begin{align*}
\btr \bchil{u}' = \bal{u}^{-1} \btr \dchil{u}',
\quad
\btr \buchil{u} = \bal{u} \btr \buchil{u}.
\end{align*}
The two torsions $\btal{u}$ and $\detal{u}$ are related by
\begin{align*}
\btal{u} = \detal{u} + \bslashd \log \bal{u}.
\end{align*}
Then the mass aspect functions $\bmul{u}$ and $\dmul{u}$ are related by
\begin{align*}
\bmul{u} = \dmul{u} - \bslashDelta \log \bal{u}.
\end{align*}

\item The equations of the lapse function $\bal{u}$ of a constant mass aspect function foliation are
\begin{align}
\left\{
\begin{aligned}
&
\textstyle
\bslashDelta \log \bal{u} = \dmul{u} - \bmul{u} = \dmul{u} - \frac{2m_{H}(\bSigma_u)}{r_u^3},
\\
&
\textstyle 
\overline{\bal{u}\, \btr \duchil{u}}
=\overline{\btr \buchil{u}} 
= \frac{2}{r_u}.
\end{aligned}
\right.
\label{eqn 2.2}
\end{align}
Note that $\bal{u}$ depends only on $\bSigma_u$, not on the whole foliation.
The above system gives rise to a map from a spacelike surface $\bSigma_u$ to its lapse function $\bal{u}$ of a constant mass aspect function foliation. We denote this map simply by $a$:
\begin{align*}
\bSigma_u \mapsto a[\bSigma_u] = \bal{u}. 
\end{align*}

\item The inverse lapse problem associated with a constant mass aspect function foliation in the background coordinate system is solving the following equation
\begin{align}
\frac{\ed}{\ed u} f(u, \vartheta) = \bal{u} = a[\bSigma_u](\vartheta).
\label{eqn 2.3}
\end{align}

\end{enumerate}

\subsection{Basic equations of constant mass aspect function foliation}\label{subsec 2.5}
In this subsection, we collect the equations on the geometry of a constant mass aspect function foliation. We shall assume that $(M,g)$ is vacuum from now on.

Let $\{\Sigma_u\}$ be a constant mass aspect function foliation of a null hypersurface $\ucalH$ in $(M,g)$. Let $\{ L', \uL \}$\footnote{We omit the superscript $u$ for the sake of brevity as remarked in footnote \ref{footnote 4}.} be the corresponding conjugate null frame on $\Sigma_u$. We introduce the decomposition of the curvature tensor relative to $\{L', \uL\}$.
\begin{definition}\label{def 2.7}
Suppose that $X,Y$ are tangent vectors of $\Sigma_u$. Define the curvature components relative to $\{ L', \uL\}$ on $\Sigma_u$ as follows,\footnote{Convention of the curvature tensor: $\mathrm{R}(X,Y)Z = \nabla_X \nabla_Y Z - \nabla_Y \nabla_X Z - \nabla_{[X,Y]} Z$, $\mathrm{R}(X,Y,Z,W) = g(\mathrm{R}(X,Y)Z, W)$.}
\begin{align*}
&
\ualpha(X,Y)=\mathrm{R}(\uL,X,Y,\uL),
&&
\alpha(X,Y) = \mathrm{R}(L',X,Y,L'),
\\
&
\ubeta(X) =\frac{1}{2} \mathrm{R}(X, \uL,L',\uL),
&&
\beta(X) = \frac{1}{2} \mathrm{R}(X,L',\uL,L'),
\\
&
\rho= \frac{1}{4} \mathrm{R}(L',\uL,\uL,L'),
&&
\sigma  \slashepsilon(X,Y) = \frac{1}{2} \mathrm{R}(X,Y,\uL,L'),
\end{align*}
where $\slashepsilon$ is the volume form of $(\Sigma_u, \slashg)$. In a local coordinate system $\{\theta^1, \theta^2\}$ of $\Sigma_u$, the above is equivalent to
\begin{align*}
&
\ualpha_{ab} = \mathrm{R}_{\uL a b \uL},
&&
\alpha_{ab}= \mathrm{R}_{L' a b L'},
\\
&
\ubeta_a = \frac{1}{2} \mathrm{R}_{a \uL L' \uL},
&&
\beta_a = \frac{1}{2} \mathrm{R}_{a L' \uL L'},
\\
&
\rho = \frac{1}{4} \mathrm{R}_{L' \uL \uL L'},
&&
\sigma \slashepsilon_{ab} = \frac{1}{2} \mathrm{R}_{a b \uL L'},
\end{align*}
where subscripts $a,b$ denote indices in $\{1,2\}$.
\end{definition}
The above curvature components contain all the information of the curvature tensor, because we have that
\begin{align*}
&
\mathrm{R}_{a b c \uL} = \slashepsilon_{ab} \slashepsilon_{cd} \ubeta^d,
&&
\mathrm{R}_{a b c L'} = \slashepsilon_{ab} \slashepsilon_{cd} \beta^d,
\\
&
\mathrm{R}_{a b c d} = \rho \slashepsilon_{ab} \slashepsilon_{cd},
&&
\mathrm{R}_{\uL a b L'} = \rho \slashg_{ab} - \sigma \slashepsilon_{ab}.
\end{align*}

We list the system of equations on the geometry of a constant mass aspect function foliation in the following. Equations \eqref{eqn 2.4}-\eqref{eqn 2.12} are called \emph{the basic equations for a constant mass aspect function foliation}.
\begin{align}
&
\phantom{\big\{ \}}
\mu = \overline{\mu},
\label{eqn 2.4}
\\
&
\phantom{\big\{ \}}
\uL \overline{\mu} 
=
- \frac{3}{2} \overline{\mu}\, \overline{\tr \uchi}
+ \frac{1}{4} \overline{\tr \chi' | \hatuchi |^2}
+ \frac{1}{2} \overline{\tr \uchi |\eta|^2},
\label{eqn 2.5}
\\
&
\phantom{\big\{ \}}
\uL \tr \uchi
= 
2 \uomega\, \tr \uchi 
- | \hatuchi |^2 
- \frac{1}{2} ( \tr \uchi )^2,
\label{eqn 2.6}
\\
&
\phantom{\big\{ \}}
\uL \tr \chi' 
= 
-2 \uomega\, \tr \chi' 
- \frac{1}{2} \tr \uchi\, \tr \chi' 
- 2 | \eta |^2
+ 2 \mu,
\label{eqn 2.7}
\\
&
\phantom{\big\{ \}}
\slashdiv \hatuchi
- \frac{1}{2} \slashd \tr \uchi
- \hatuchi \cdot \eta 
+ \frac{1}{2} \tr \uchi \, \eta
=
- \ubeta,
\label{eqn 2.8}
\\
&
\phantom{\big\{ \}}
\slashdiv \hatchi'
- \frac{1}{2} \slashd \tr \chi'
+ \hatchi' \cdot \eta 
- \frac{1}{2} \tr \chi' \, \eta
=
- \beta,
\label{eqn 2.9}
\\
&
\left\{
\begin{aligned}
&
\slashcurl \eta
=
\frac{1}{2} \hatchi' \wedge \hatuchi 
+ \sigma,
\\
&
\slashdiv \eta 
= 
- \rho 
- \frac{1}{2} ( \hatuchi, \hatchi' ) 
- \mu,
\label{eqn 2.10}
\end{aligned}
\right.
\\
&
\phantom{\big\{ \}}
\begin{aligned}
&
2 \slashDelta \uomega
=
-\frac{3}{2} ( \mu\, \tr \uchi - \overline{\mu\, \tr \uchi} )
+ \frac{1}{2} ( \tr \uchi | \eta |^2 - \overline{\tr \uchi | \eta |^2} )
+ \frac{1}{4} ( \tr \chi' | \hatuchi |^2 - \overline{ \tr \chi' | \hatuchi|^2} )
\\
&
\phantom{2 \slashDelta \uomega =}
+ 4 (\slashdiv \hatuchi, \eta )
+ 4 ( \hatuchi, \slashnabla \eta )
- 2 \slashdiv \ubeta.
\end{aligned}
\label{eqn 2.11}
\end{align}
In the above equations, the dot $\cdot$ and the parenthesis $(\cdot, \cdot)$ are both the inner product relative to the metric $\slashg$ that
\begin{align*}
&
(\hatuchi \cdot \eta)_a = \hatuchi_{ab}\, \eta_c\, (\slashg^{-1})^{bc},
\quad \quad
(\hatchi' \cdot \eta)_a = \hatchi'_{ab}\, \eta_c\, (\slashg^{-1})^{bc},
\\
&
( \hatchi', \hatuchi ) = \hatchi'_{ab}\, \hatuchi_{cd}\, ( \slashg^{-1} )^{ac}\, ( \slashg^{-1} )^{bd},
\quad \quad
( \hatuchi, \slashnabla \eta ) = \hatchi'_{ab}\, \slashnabla_c \eta_d\, ( \slashg^{-1} )^{ac}\, ( \slashg^{-1} )^{bd},
\\
&
(\slashdiv \hatuchi, \eta )= (\slashdiv \hatuchi)_a\, \eta_b\, (\slashg^{-1})^{ab},
\end{align*}
the wedge operator $\wedge$ is defined by
\begin{align*}
\hatchi' \wedge \hatuchi = \hatchi'_{ab}\, \hatuchi_{cd}\, (\slashg^{-1})^{ac}\, \slashepsilon^{bd},
\end{align*}
and $\slashcurl$ is the curl operator on $(\Sigma_u, \slashg)$ that $(\slashcurl \eta)_a = \slashepsilon^{ab} \slashnabla_a \eta_b$. 

Note that equation \eqref{eqn 2.11} determines the acceleration $\uomega$ up to a constant, which gives essentially the freedom of relabelling the leaves of the foliation. In order to remove this freedom, we require that $r_{u_0+u}= r_{u_0} +u$, which is equivalent to $\overline{\tr \uchi} = \frac{2}{r}$ on every $\Sigma_u$. For such a choice of the labelling, we can determine $\uomega$ completely with an additional equation
\begin{align}
\overline{\uomega}
=
-\frac{r}{2} \overline{(\uomega - \overline{\uomega})(\tr\uchi - \overline{\tr\uchi})}
- \frac{r}{8} \overline{(\tr \uchi - \overline{\tr \uchi})^2}
+ \frac{r}{4} \overline{|\hatuchi|^2}.
\label{eqn 2.12}
\end{align}

We elaborate on the structures and geometric meanings of equations \eqref{eqn 2.4}-\eqref{eqn 2.12}.
\begin{enumerate}[label=\alph*.]
\item Propagation equations: equations \eqref{eqn 2.4}-\eqref{eqn 2.7}. From these equations, we can solve $\mu, \tr \uchi, \tr \chi'$ by integrations from their values at an arbitrary leaf. It is worth to point out that equation \eqref{eqn 2.7} is equivalent to
\begin{align*}
\uL \tr \chi' 
= 
-2 \uomega\, \tr \chi' 
- \frac{1}{2} \tr \uchi\, \tr \chi' 
- (\hatchi', \hatuchi) 
- 2( \slashdiv \eta + | \eta |^2) 
-2 \rho,
\end{align*}
as we have that
\begin{align*}
\mu =  - \rho - \frac{1}{2} ( \hatchi', \hatuchi ) - \slashdiv \eta,
\end{align*}
which follows from the Gauss equation of $\Sigma_u$
\begin{align}
K - \frac{1}{4} \tr \chi' \, \tr \uchi + \frac{1}{2} ( \hatchi', \hatuchi ) = - \rho.
\label{eqn 2.13}
\end{align}

\item Elliptic equations: equations \eqref{eqn 2.8}-\eqref{eqn 2.12}. If we assume that $\mu, \tr \uchi, \tr \chi'$ and the curvature components $\ubeta, \beta, \rho, \sigma$ are known, then we can obtain $\hatuchi, \hatchi', \eta, \uomega$. 

\item Above propagation equations and elliptic equations form a coupled system. If the curvature components $\ubeta, \beta, \rho, \sigma$ are known, then we can solve the mass aspect function $\mu$ and the connection coefficients $ \tr \uchi$, $\tr \chi'$, $\hatuchi$, $\hatchi'$, $\eta$, $\uomega$.

\item Geometric meanings: equations \eqref{eqn 2.7} and \eqref{eqn 2.8} are the Codazzi equations of the null second fundamental forms $\uchi$ and $\chi'$, and the first equation of \eqref{eqn 2.10} is the Gauss equation for the curvature of the normal bundle of $\Sigma_u$.

\end{enumerate}

We conclude this section with a comparison of the constant mass aspect function foliation with a general foliation. Note that equations \eqref{eqn 2.6}-\eqref{eqn 2.9} are valid for any foliation, and equation \eqref{eqn 2.12} holds for any foliation with the condition $\overline{\tr \uchi} =\frac{2}{r}$. While equations \eqref{eqn 2.5} and \eqref{eqn 2.11} are no longer valid for a general foliation. Instead, the following equations hold,
\begin{align*}
&
\uL \overline{\mu} 
= 
- \overline{\mu}\, \overline{\tr \uchi}
- \frac{1}{2} \overline{\mu\, \tr \uchi}
+ \frac{1}{4} \overline{\tr \chi'\, | \hatuchi|^2}
+ \frac{1}{2} \overline{\tr \uchi\, | \eta |^2},
\tag{\ref{eqn 2.5}$'$}
\label{eqn 2.5'}
\\
&
\begin{aligned}
&
\uL \mu
=
- 2 \slashDelta \uomega
- \frac{3}{2} \mu\, \tr \uchi 
+ \frac{1}{2}  \tr \uchi\, | \eta |^2 
+ \frac{1}{4}  \tr \chi'\, | \hatuchi |^2 
+ 4 (\slashdiv \hatuchi, \eta )
\\
&
\phantom{ 2 \slashDelta \uomega =}
+ 4 ( \hatuchi, \slashnabla \eta )
- 2 \slashdiv \ubeta.
\end{aligned}
\tag{\ref{eqn 2.11}$'$}
\label{eqn 2.11'}
\end{align*}
We refer to references \cite{C2009} \cite{S2008} \cite{L2018} for the proofs of equations \eqref{eqn 2.5}-\eqref{eqn 2.12}, \eqref{eqn 2.5'} and \eqref{eqn 2.11'}. However for the sake of self-containedness, we shall derive equations \eqref{eqn 2.5'} and \eqref{eqn 2.11'} in Appendix \ref{appen A}.

\subsection{Asymptotic geometry of foliation at null infinity}\label{subsec 2.6}
As mentioned in subsection \ref{subsec 2.3}, we are interested in the asymptotic geometry of the constant mass aspect function foliation at null infinity, for its importance in the application to the Penrose inequality pointed out in \cite{S2008}. Thus we introduce the geometric quantities describing the asymptotic geometry of a constant mass aspect function foliation at null infinity.

Adopting the same notations in subsection \ref{subsec 2.5}, $\{\Sigma_u\}$ is a constant mass aspect function foliation of a null hypersurface $\ucalH$ in $(M,g)$. In order to describe the asymptotic geometry of $\{\Sigma_u\}$ at null infinity, we consider the renormalised geometry of $\Sigma_u$: define the renormalised metric $\slashgl{u,r}$ on $\Sigma_u$\footnote{We shall use the superscript $r$ on the left to indicate renormalised quantities.}
\begin{align}
\slashgl{u,r} = r_u^{-2} \cdot \slashgl{u},
\label{eqn 2.14}
\end{align}
and define the renormalised Gauss curvature $\Kl{u,r}$ of $\Sigma_u$, which is the Gauss curvature of the renormalised metric $\slashgl{u,r}$,
\begin{align}
\Kl{u,r} = r_u^2 \cdot \Kl{u}.
\label{eqn 2.15}
\end{align}
When $u$ approaches $+\infty$,  the leaf $\Sigma_u$ approaches null infinity. We take the limit of $\slashgl{u,r}$ and $\Kl{u,r}$ as $u \rightarrow +\infty$. If the limits exist, we define them as the limit renormalised metric and Gauss curvature of the foliation at null infinity
\begin{align}
\slashgl{\infty, r} = \lim_{u \rightarrow +\infty} \slashgl{u,r},
\quad
\Kl{\infty, r} = \lim_{u \rightarrow +\infty} \Kl{u,r}.
\label{eqn 2.16}
\end{align}
When the limit renormalised metric $\slashgl{\infty,r}$ is round, or equivalently the limit renormalised Gauss curvature $\Kl{\infty,r}$ is constant one, we call the foliation $\{ \Sigma_u \}$ asymptotically round.

\subsection{Asymptotic reference frame and energy-momentum vector at null infinity}\label{subsec 2.7}
When the foliation $\{ \Sigma_u \}$ is asymptotically round, we call that $\{\Sigma_u\}$ defines an asymptotic reference frame, as in \cite{S2008}. Assume so for the foliation $\{ \Sigma_u \}$ and denote the corresponding reference frame by $\gamma_{\infty}$. Following \cite{S2008}, we define the Bondi energy relative to the reference frame $\gamma_{\infty}$.
\begin{definition}[Bondi energy]\label{def 2.8}
Assume that the foliation $\{ \Sigma_u \}$ of a null hypersurface $\ucalH$ defines the asymptotic reference frame $\gamma_{\infty}$ at null infinity. The Bondi energy of $\ucalH$ relative to $\gamma_{\infty}$ is defined by
\begin{align}
E^{\gamma_{\infty}}(\ucalH) 
=
\lim_{u \rightarrow \infty} m_H ( \Sigma_u ).
\label{eqn 2.17}
\end{align}
\end{definition}

We introduce the following asymptotic quantities as in \cite{C1991} \cite{CK1993}:
\begin{align*}
P 
&=
\lim_{u \rightarrow \infty} r^3 \cdot \rhol{u},
\\
\varSigma
&=
\lim_{u \rightarrow \infty} \hatuchil{u},
\\
\Xi
&=
\lim_{u \rightarrow \infty} r^{-1} \cdot \hatchil{u},
\\
N
&=
-P - \frac{1}{2} ( \varSigma, \Xi )_{\subslashgl{\infty,r}}.
\end{align*}
Since the Hawking mass can be calculated by
\begin{align*}
m_H(\Sigma_u) = \frac{r}{8\pi} \int_{\Sigma_u} \big[ -\rhol{u} - \frac{1}{2} (\hatuchil{u}, \hatchil{u})_{\subslashgl{u}} \big] \dvol_{\subslashgl{u}},
\end{align*}
then formula \eqref{eqn 2.17} is equivalent to
\begin{align}
E^{\gamma_{\infty}}(\ucalH) 
=
\frac{1}{8\pi} \int N \dvol_{\subslashgl{\infty,r}}.
\label{eqn 2.18}
\end{align}

Since $\slashgl{\infty,r}$ is round, we choose an isometric embedding to the standard round sphere of radius $1$ centring at the origin. Then pull back the Cartesian coordinate functions $\{x^1, x^2, x^3\}$ to null infinity. Relative to this set of functions $\{x^1, x^2, x^3\}$, we can define the linear momentum of $\ucalH$ at null infinity with respect to the reference frame $\gamma_{\infty}$.
\begin{definition}\label{def 2.9}
Assume that the foliation $\{ \Sigma_u \}$ of a null hypersurface $\ucalH$ defines the asymptotic reference frame $\gamma_{\infty}$ at null infinity. Let $\{x^1, x^2, x^3\}$ be a set of functions at null infinity obtained from pulling back the Cartesian coordinate functions via an isometric embedding of $\slashgl{\infty,r}$ into $3$-dim Euclidean space. 

Relative to $\{x^1, x^2, x^3\}$, the linear momentum $\vec{P}^{\gamma_{\infty}}=(P^{\gamma_{\infty},1}, P^{\gamma_{\infty},2}, P^{\gamma_{\infty},3})$ of $\ucalH$ at null infinity with respect to the reference frame $\gamma_{\infty}$ is defined by
\begin{align}
P^{\gamma_{\infty},i}(\ucalH) 
=
\frac{1}{8\pi} \int x^i\cdot N \dvol_{\subslashgl{\infty,r}}.
\label{eqn 2.19}
\end{align}
Together with the Bondi energy, the four dimensional vector $(E^{\gamma_{\infty}}, \vec{P}^{\gamma_{\infty}})$ is the energy-momentum vector of $\ucalH$ at null infinity with respect to the reference frame $\gamma_{\infty}$.

If the linear momentum $\vec{P}^{\gamma_{\infty}}$ vanishes, then we call $\gamma_{\infty}$ an asymptotic center-of-mass reference frame.
\end{definition}

Following \cite{S2008}, we give the definition of the Bondi mass of $\ucalH$ at null infinity.
\begin{definition}[Bondi mass]\label{def 2.10}
Assume that the foliation $\{ \Sigma_u \}$ of a null hypersurface $\ucalH$ defines the asymptotic reference frame $\gamma_{\infty}$ at null infinity. The energy-momentum vector of $\ucalH$ with respect to the reference frame $\gamma_{\infty}$ is $(E^{\gamma_{\infty}}, \vec{P}^{\gamma_{\infty}})$. Then the Bondi mass of $\ucalH$ is defined by
\begin{align}
m_B(\ucalH) = \sqrt{(E^{\gamma_{\infty}})^2 - | \vec{P}^{\gamma_{\infty}} |^2}.
\label{eqn 2.20}
\end{align}
If $\gamma_{\infty}$ is an asymptotic center-of-mass reference frame, i.e. $\vec{P}^{\gamma_{\infty}} =0$, then
\begin{align*}
m_B(\ucalH) = E^{\gamma_{\infty}}.
\end{align*}
\end{definition}

\section{Schwarzschild black hole spacetime}
In this section, we review the Schwarzschild black hole spacetime. We introduce a double null coordinate system and give the geometric information of the spacetime relative to this coordinate system.

\subsection{Background double null coordinate system}\label{subsec 3.1}
In the coordinate system $\{ t, r, \theta, \phi\}$, the Schwarzschild metric $g$ is
\begin{align*}
\textstyle
g=-\left( 1-\frac{2m}{r} \right) \ed t^2 + \left( 1-\frac{2m}{r} \right)^{-1} \ed r ^2 + r^{2} \left( \ed \theta^2 + \sin^2 \theta \ed \phi^2 \right).
\end{align*}
Let $r_0=2m$, where $m$ is the mass, and $r_0$ is the area radius of the event horizon of the black hole. Consider the following coordinate transformation $\{ t, r\} \leftrightarrow \{ \us, s\}$ in \cite{L2018} \cite{L2022}
\begin{align*}
\left\{
\begin{aligned}
&
(r-r_0)^{\frac{1}{2}} \exp \frac{t+r}{2r_0}
=
\exp\frac{\us}{r_0},
\\
&
(r-r_0)^{\frac{1}{2}} \exp \frac{-t+r}{2r_0}
=
s \exp \frac{s+r_0}{r_0},
\end{aligned}
\right.
\end{align*}
then the Schwarzschild metric $g$ takes the form
\begin{align*}
g
= 
2 \Omega^2 \left( \ed s\otimes \ed \us + \ed \us \otimes \ed s  \right) + r^2 \left( \ed \theta^2 + \sin^2 \theta \ed \phi^2 \right),
\end{align*}
where
\begin{align*}
\Omega^2 = \frac{s+r_0}{r}\exp\frac{\us+s+r_0-r}{r_0},
\quad
(r-r_0)\exp\frac{r}{r_0} = s\exp\frac{\us+s+r_0}{r_0}.
\end{align*}
We shall use $\circg$ to denote the standard metric on the sphere of radius $1$ that $\circg = \ed \theta^2 + \sin^2 \theta \ed \phi^2$.

The coordinate system $\{ s, \us, \theta, \phi \}$ is called a double null coordinate system, as the level sets of $\us,s$ are both null hypersurfaces. Denote the level set of $\us$ by $\uC_{\us}$ and the level set of $s$ by $C_s$. Let $\Sigma_{s,\us}$ be the intersection of $\uC_{\us}$ and $C_s$. We denote the metric on $\Sigma_{s,\us}$ by $\slashg$, which is in fact $\slashg = r^2 \circg$.

In the above double null coordinate system, coordinates $\us$, $s$ are essential, while $\theta, \phi$ can be replaced by other coordinates on the sphere. For example, let $\{\theta^1, \theta^2\}$ be a coordinate system on the sphere, then we can construct a double null coordinate system $\{\us, s, \theta^1, \theta^2\}$ for the Schwarzschild spacetime where the metric takes the form
\begin{align*}
g
= 
2 \Omega^2 \left( \ed s\otimes \ed \us + \ed \us \otimes \ed s  \right) + r^2 \circg_{ab} \ed \theta^a \otimes \ed \theta^b.
\end{align*}
Sometimes for the sake of brevity, we simply use the notation $\vartheta$ to denote the coordinates $(\theta^1, \theta^2)$ on the sphere.

\subsection{Geometry of background coordinate system}\label{subsec 3.2}

We construct a conjugate null frame relative to a double null coordinate system $\{\us, s, \vartheta\}$. Construct the tangent null vector field $\uL$ on each null hypersurface $\uC_{\us}$ such that
\begin{align*}
\uL s =1.
\end{align*}
Then define $L'$ as the conjugate normal null vector of $\Sigma_{s,\us}$ relative to $\uL$ such that $g(L', \uL)=2$. Then $\{L', \uL\}$ is a conjugate null frame of $\Sigma_{s,\us}$. This construction is the same as in subsection \ref{subsec 2.2}.

We can express $L', \uL$ in terms of the coordinate frame, that
\begin{align*}
\uL = \partial_s,
\quad
L' = \Omega^{-2} \partial_{\us}.
\end{align*}
Symmetrically, we can also construct another conjugate null frame $\{ L, \uL'\}$ of $\Sigma_{s,\us}$ that
\begin{align*}
L = \partial_{\us},
\quad
\uL' = \Omega^{-2} \partial_s.
\end{align*}
In this paper, we shall mainly work with the conjugate null frame $\{ L', \uL\}$. We define the background connection coefficients relative to the double null coordinate system $\{s, \us\}$, by adopting definitions \ref{def 2.1} and \ref{def 2.4}.
\begin{definition}[Background connection coefficients]\label{def 3.1}
Let $X,Y$ be vector fields tangential to $\Sigma_{s,\us}$. We define that
\begin{align*}
&
\uchi(X,Y) = g(\nabla_X \uL, Y),
&&
\chi(X,Y) = g(\nabla_X L, Y),
\\
&
\uchi'(X,Y) = g(\nabla_X \uL, Y),
&&
\chi'(X,Y) = g(\nabla_X L, Y),
\\
&
\eta(X) = \frac{1}{2} g(\nabla_X \uL, L'),
&&
\ueta(X) = \frac{1}{2} g(\nabla_X L, \uL'),
\\
&
\omega = \frac{1}{4} g(\nabla_L L, \uL') =  L \log \Omega
&&
\uomega = \frac{1}{4} g(\nabla_{\uL}\uL, L') =  \uL \log \Omega.
\end{align*}
We decompose $\chi, \uchi$ into their traces and trace-free parts relative to $\slashg$:
\begin{align*}
\chi = \hatchi + \frac{1}{2} \tr \chi \slashg,
\quad
\uchi = \hatuchi + \frac{1}{2} \tr \uchi \slashg,
\end{align*}
and similarly for $\chi', \uchi'$.
\end{definition}

We also adopt definition \ref{def 2.7} to define the background curvature components relative to $\{L', \uL\}$.
\begin{definition}[Background curvature components]\label{def 3.2}
Suppose that $X,Y$ are tangent vectors of $\Sigma_{s,\us}$. Define the background curvature components relative to $\{ L', \uL\}$ on $\Sigma_{s,\us}$ as follows,
\begin{align*}
&
\ualpha(X,Y)=\mathrm{R}(\uL,X,Y,\uL),
&&
\alpha(X,Y) = \mathrm{R}(L',X,Y,L'),
\\
&
\ubeta(X) =\frac{1}{2} \mathrm{R}(X, \uL,L',\uL),
&&
\beta(X) = \frac{1}{2} \mathrm{R}(X,L',\uL,L'),
\\
&
\rho= \frac{1}{4} \mathrm{R}(L',\uL,\uL,L'),
&&
\sigma  \slashepsilon(X,Y) = \frac{1}{2} \mathrm{R}(X,Y,\uL,L'),
\end{align*}
where $\slashepsilon$ is the volume form of $(\Sigma_{s,\us}, \slashg)$. In a double null coordinate system $\{\us, s, \theta^1, \theta^2\}$, the above is equivalent to
\begin{align*}
&
\ualpha_{ab} = \mathrm{R}_{\uL a b \uL},
&&
\alpha_{ab}= \mathrm{R}_{L' a b L'},
\\
&
\ubeta_a = \frac{1}{2} \mathrm{R}_{a \uL L' \uL},
&&
\beta_a = \frac{1}{2} \mathrm{R}_{a L' \uL L'},
\\
&
\rho = \frac{1}{4} \mathrm{R}_{L' \uL \uL L'},
&&
\sigma \slashepsilon_{ab} = \frac{1}{2} \mathrm{R}_{a b \uL L'},
\end{align*}
where subscripts $a,b$ denote indices in $\{1,2\}$.
\end{definition}

Above definitions \ref{def 3.1}, \ref{def 3.2} are not confined to the Schwarzschild spacetime, but apply to a double null coordinate system in any vacuum spacetime.\footnote{Definition \ref{def 3.1} also applies to non-vacuum spacetimes.} We can compute the background connection coefficients and curvature components in the Schwarzschild spacetime explicitly. We list the results here:
\begin{align}
\begin{aligned}
&
\partial_{\us} r= \frac{r-r_0}{r}, 
\quad
\partial_s r= \frac{s+r_0}{r} \cdot \frac{r-r_0}{s},
\\
&
\tr \chi'= \frac{2s}{r(s + r_0)},
\quad
\tr\uchi = \frac{2(s+r_0)}{r^2} \cdot \frac{r-r_0}{s} = \frac{2(s+r_0)}{r^2} \exp\frac{\us+s+r_0-r}{r_0},
\\
&
\hatchi' = \hatuchi =0,
\quad
\eta = \ueta= 0,
\\
&
\omega =\frac{r_0}{2r^2}, 
\quad 
\uomega =\frac{1}{2(s+r_0)} +\frac{1}{2r_0} - \Big( \frac{1}{2r} + \frac{1}{2r_0} \Big)  \frac{s+r_0}{r} \exp \frac{\us+s+r_0-r}{r_0},
\end{aligned}
\label{eqn 3.1}
\end{align}
and
\begin{align}
\alpha = \ualpha =0,
\quad
\beta = \ubeta = 0,
\quad
\sigma=0,
\quad
\rho = - \frac{ r_0}{r^3}.
\label{eqn 3.2}
\end{align}
We can also compute the mass aspect function $\mu$ relative to $\{L', \uL\}$ which is
\begin{align*}
\mu = K - \frac{1}{4} \tr \chi'\, \tr \uchi - \slashdiv \eta = \frac{r_0}{r^3}.
\end{align*}
Thus the foliation $\{\Sigma_{s,\us}\}_{s}$ is a constant mass aspect function foliation on $\uC_{\us}$.

\section{Parameterisation and geometry of a spacelike surface}
In this section, we introduce the method to parameterise a spacelike surface with the help of a double null coordinate system $\{s, \us, \vartheta\}$ in a Schwarzschild spacetime. Then we give the formulae of the connection coefficients and curvature components on a spacelike surface in terms of its parameterisation and the background geometric quantities.

\subsection{Parameterisation of a spacelike surface}\label{subsec 4.1}
We adopt the methods to parameterise a spacelike surface in \cite{L2018}, \cite{L2020}, which introduced two kinds of parameterisations. Let $\Sigma$ be a spacelike surface in the Schwarzschild spacetime, which is close to some $\Sigma_{s,\us}$ in the double null coordinate system.
\begin{enumerate}[label=\ref{subsec 4.1}.\alph*.,leftmargin=.45in]
\item First kind of parameterisation $(f,\uf)$. Parameterise $\Sigma$ by a pair of two functions $(f,\uf)$ as its graph of $(s,\us)$ over $\vartheta$ domain in the double null coordinate system, i.e.
\begin{align*}
\varphi_{(f,\uf)}:
\quad
\vartheta
\mapsto
(s,\us,\vartheta) = ( f(\vartheta), \uf(\vartheta), \vartheta) \in \Sigma.
\end{align*}
See figure \ref{fig 6}.
\begin{figure}[H]
\begin{center}
\begin{tikzpicture}
\draw[dashed] (-1,0)
to [out=70, in=180] (0,0.5)
to [out=0,in=110] (1,0);
\draw (1,0)
to [out=-70,in=0] (0,-0.8) node[below]{\tiny $\Sigma_{0,0}$}
to [out=180,in=-110] (-1,0); 
\draw[dashed] (-1,0) to [out=70,in=-110] (-0.7,0.9);
\draw (-1,0) to [out=-110,in=70] (-1.9,-2.5);
\draw[dashed] (1,0) to [out=110,in=-70] (0.7,0.9);
\draw[->] (1,0) to [out=-70,in=110] (1.9,-2.5) node[right] {\small $s$}; 
\draw[dashed] (-1,0) to [out=-45,in=135] (-0.5,-0.5);
\draw (-1,0) to [out=135, in= -45] (-2,1);
\draw[dashed] (1,0) to [out=-135,in=45] (0.5,-0.5);
\draw[->] (1,0) to [out=45, in= -135] (2,1) node[right] {\small $\us$}; 
\node[below right] at (1.4,0.6) {\tiny $C_{s=0}$};
\node[above right] at (1.55,-2) {\tiny $\uC_{\us=0}$};
\draw[dashed] (-2.7,-2.2+1) to [out=70,in=180] (-0.5,-2.5+1) to [out=0,in=110] (2.7,-2.2+1);
\draw (2.7,-2.2+1) node[right] {\scriptsize $\Sigma=\{(s,\us,\vartheta): s= f(\vartheta), \us=\uf(\vartheta)\}$} to [out=-70,in=0] (1,-2+0.8) to [out=180,in=-110] (-2.7,-2.2+1); 
\end{tikzpicture}
\end{center}
\caption{First kind of parameterisation of $\Sigma$.}
\label{fig 6}
\end{figure}
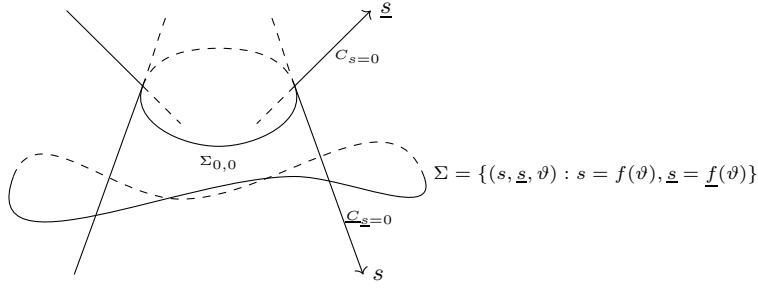

\item Second kind of parameterisation $(f, \ufl{s=0})$. Let $\ucalH$ be the incoming null hypersurface containing $\Sigma$. Denote the intersection of $\ucalH$ and $C_{s=0}$ by $\Sigma_0$. Parameterise $\Sigma_0$ by a function $\ufl{s=0}$ as its graph of $\us$ over $\vartheta$ domain in the coordinate system $\{\us, \vartheta\}$ in $C_{s=0}$. Restrict the coordinate functions $s,\vartheta$ on $\ucalH$ to obtain a coordinate system $\{ s, \vartheta\}$ on $\ucalH$. Then parametrise $\Sigma$ in $\ucalH$ by a function $f$ as its graph of $s$ over $\vartheta$ domain in the coordinate system $\{s,\vartheta\}$. See figure \ref{fig 7}. 
\begin{figure}[H]
\begin{center}
\begin{tikzpicture}
\draw[dashed] (-1,0)
to [out=70, in=180] (0,0.5)
to [out=0,in=110] (1,0);
\draw (1,0)
to [out=-70,in=0] (0,-0.8) node[below]{\tiny $\Sigma_{0,0}$}
to [out=180,in=-110] (-1,0); 
\draw[dashed] (-1,0) to [out=70,in=-110] (-0.7,0.9);
\draw (-1,0) to [out=-110,in=70] (-1.85,-2.4);
\draw[dashed] (1,0) to [out=110,in=-70] (0.7,0.9);
\draw[->] (1,0) to [out=-70,in=110] (1.85,-2.4) node[right] {\small $s$}; 
\node[above right] at (1.5,-1.9) {\tiny $\uC_{\us=0}$};
\draw[dashed] (-1,0) to [out=-45,in=135] (-0.5,-0.5);
\draw (-1,0) to [out=135, in= -45] (-2,1);
\draw[dashed] (1,0) to [out=-135,in=45] (0.5,-0.5);
\draw[->] (1,0) to [out=45, in= -135] (1.8,0.8) node[right] {\tiny $C_{s=0}$}to [out=45,in=-135] (2.3,1.3) node[right] {\small $\us$}; 
\draw[dashed] (-1.5,0.5) to [out=135,in=45] (1.5,0.5);
\draw (1.5,0.5) to [out=-135,in=0] (0,0) node[below]{\scriptsize $\Sigma_0$} to [out=180,in=-45] (-1.5,0.5); 
\draw[dashed] (-1.4,1.3) to [out=-110,in=70] (-1.6,0.7);
\draw (-1.6,0.7) to [out=-110,in=70] (-2.75,-2.4);
\draw[dashed] (1.4,1.3) to [out=-70,in=110] (1.6,0.7);
\draw[->] (1.6,0.7) to [out=-70,in=110] (2.75,-2.4) node[right]{\small $s$}; 
\node[right] at (2.4,-1.5) {\scriptsize $\ucalH$}; 
\draw[dashed] (-2.3,-2.2+1) to [out=70,in=180] (-0.5,-2.5+1) node[below] {\scriptsize $\Sigma$} to [out=0,in=110] (2.3,-2.2+1);
\draw (2.3,-2.2+1) to [out=-70,in=0] (1,-2+0.9) to [out=180,in=-110] (-2.3,-2.2+1); 
\draw[->] (-1.73,0.73) to [out=-110,in=70] (-2.45,-1.2);
\node[left] at (-2.3,-0.8) {$f$};
\draw[->] (-1.05,-0.1) to [out=135,in=-45] (-1.65,0.5);
\node[below] at (-1.55,0.1) {$\ufl{s=0}$};
\end{tikzpicture}
\end{center}
\caption{Second kind of parameterisation of $\Sigma$.}
\label{fig 7}
\end{figure}
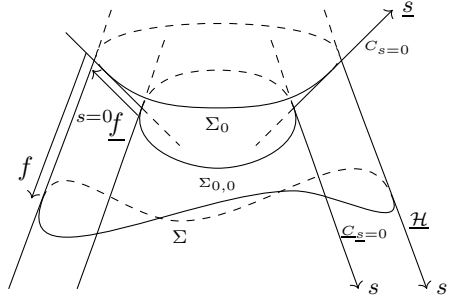
\end{enumerate}

Two methods to transform from the second parameterisation $(f,\ufl{s=0})$ to the first paramterisation $(f,\uf)$ were introduced in \cite{L2018} \cite{L2020}. We summarise them here. Firstly, two parameterisations share the same parameterisation function $f$.
\begin{enumerate}[label=\ref{subsec 4.1}.\Roman*.,leftmargin=.45in]
\item Method I. Parameterise $\ucalH$ by a function $\uh$ as its graph of $\us$ over $(s,\vartheta)$ domain in the double null coordinate system $\{ s, \us, \vartheta \}$, i.e.
\begin{align*}
\varphi_{\uh}:
\quad
(s,\vartheta)
\mapsto
(s, \us, \vartheta) = (s, \uh(s,\vartheta), \vartheta) \in \ucalH.
\end{align*}
Then the parameterisation function $\uf$ is given by
\begin{align*}
\uf(\vartheta) = \uh(f(\vartheta), \vartheta).
\end{align*}
In order to obtain $\uf$, it is sufficient to determine $\uh$. $\uh$ satisfies the following first order fully nonlinear equation: let $\dpartial_s, \dpartial_{i=1,2}$ be the coordinate derivatives of the coordinate system $\{s, \vartheta\}$ on $\ucalH$,
\begin{align}
\dpartial_s \uh = \Omega^2 (\slashg^{-1})^{ij} \dpartial_i \uh\, \dpartial_j \uh,
\label{eqn 4.1}
\end{align}
with the initial data $\uh(s=0, \vartheta) = \ufl{s=0}(\vartheta)$.

\item Method II. Introduce functions $\fl{t}= t f$, $t\in [0,1]$. Introduce a family of surfaces $\{S_t\}$ in $\ucalH$, where $S_t$ is the graph of $\fl{t}$ of $s$ over $\vartheta$ domain in the coordinate system $\{ s, \vartheta\}$ on $\ucalH$. Thus $S_t$ has the second parameterisation $(\fl{t}, \ufl{s=0})$. See figure \ref{fig 8}.
\begin{figure}[H]
\begin{center}
\begin{tikzpicture}
\draw (-1,0) to [out=135, in= -45] (-2.3,1.3);
\draw[->] (1,0) to [out=45, in= -135] (1.8,0.8) node[right] {\tiny $C_{s=0}$}to [out=45,in=-135] (2.3,1.3) node[right] {\small $\us$};
\draw[dashed] (-1.5,0.5) to [out=135,in=45] (1.5,0.5);
\draw (1.5,0.5) to [out=-135,in=0] (0,0) node[above]{\scriptsize $S_{t=0}:\Sigma_0$} to [out=180,in=-45] (-1.5,0.5); 
\draw[dashed] (-1.4,1.3) to [out=-110,in=70] (-1.6,0.7);
\draw (-1.6,0.7) to [out=-110,in=70] (-3,-3);
\draw[dashed] (1.4,1.3) to [out=-70,in=110] (1.6,0.7);
\draw[->] (1.6,0.7) to [out=-70,in=110] (3,-3) node[right]{\small $s$}; 
\node[right] at (2.4,-1.5) {\scriptsize $\ucalH$};
\draw[dashed] (-2.3+0.23,-2.2+1.6) to [out=70,in=180] (-0.4,-0.7) node[below] {\scriptsize $S_t$} to [out=0,in=110] (2.3-0.23,-2.2+1.6);
\draw (2.3-0.23,-2.2+1.6) to [out=-70,in=0] (0.7,-0.5) to [out=180,in=-110] (-2.3+0.23,-2.2+1.6); 
\draw[->] (-1.71,0.71) to [out=-110,in=70] (-2.45+0.25,-1.2+0.6);
\node[left] at (-2,-0.3) {$\fl{t}$}; 
\draw[dashed] (-2.3-0.4,-2.2) to [out=70,in=180] (-0.5,-2.5) node[below] {\scriptsize $S_{t=1}: \Sigma$} to [out=0,in=110] (2.3+0.4,-2.2);
\draw (2.3+0.4,-2.2) to [out=-70,in=0] (1,-2) to [out=180,in=-110] (-2.3-0.4,-2.2); 
\draw[->] (-1.7-0.3,0.7+0.3) to [out=-110,in=70] (-2.3-0.8,-2.2+0.4);
\draw (-2.3-0.4,-2.2) to [out=135, in= -45] (-2.3-0.4-0.7,-2.2+0.7);
\node[left] at (-2.3-0.55,-2.2+1-0.2) {$f$}; 
\end{tikzpicture}
\end{center}
\caption{The family of surfaces $\{ S_t \}$.}
\label{fig 8}
\end{figure}
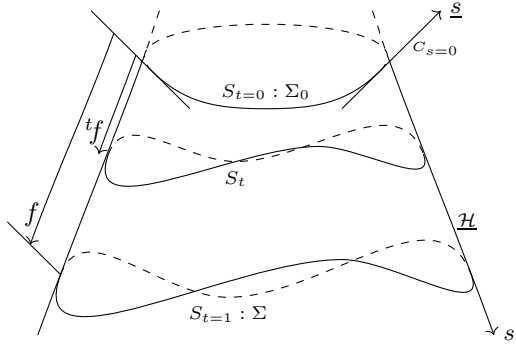
Suppose that the first parameterisation of $S_t$ is $(\fl{t}, \ufl{t})$. We see that $S_{t=0}= \Sigma_0$ and $S_{t=1} = \Sigma$. Then $\ufl{t=0} = \ufl{s=0}$ and $\ufl{t=1} = \uf$. Thus in order to determine $\uf$, it is sufficient to obtain $\ufl{t=1}$. We have that $\ufl{t}$ satisfies the following first order fully nonlinear equation:
\begin{align}
\partial_t \ufl{t} = f \cdot [ 1- t\underline{e}^i f_i - t e^i f_i \cdot \uvarepsilon ]^{-1} \cdot \left[ \uvarepsilon - \underline{e}^i \, \ufl{t}_i - e^i \, \ufl{t}_i \cdot \uvarepsilon \right],
\label{eqn 4.2}
\end{align}
where
\begin{align*}
\begin{aligned}
&
\underline{\varepsilon} = \frac{ -|\underline{e}|^2}{(2\Omega^2 + e\cdot \underline{e}) + \sqrt{(2\Omega^2 + e\cdot \underline{e})^2 -|e|^2 |\underline{e}|^2}},
\\
&
|e|^2 = \slashg_{ij}e^ie^j,
\quad
|\underline{e}|^2 = \slashg_{ij} \underline{e}^i \underline{e}^j, 
\quad
e\cdot \underline{e} =\slashg_{ij} e^i \underline{e}^j,
\\
&
e^k =-2\Omega^2 \cdot tf_j \left(\slashg^{-1}\right)^{jk}, 
\quad
\underline{e}^k = -2\Omega^2 \cdot \ufl{t}_j \left( \slashg^{-1} \right)^{jk},
\end{aligned}
\end{align*}
and $f_i = \partial_i f$, $\ufl{t}_i = \partial_i \ufl{t}$.
\end{enumerate}

\subsection{Geometry of a null hypersurface with induced coordinate system}\label{subsec 4.2}
In the second kind of parameterisation of a spacelike surface $\Sigma$, we consider an intermediate step which is the embedding of $\Sigma$ in an incoming null hypersurface $\ucalH$.

When considering the geometry of the spacelike surface $\Sigma$, we can make use of this intermediate step, i.e. we first study the geometry of the incoming null hypersurface $\ucalH$, then study the geometry of $\Sigma$ by viewing it as an embedded surface in $\ucalH$. This is the strategy adopted in \cite{L2018}, which is also the same strategy calculating the outgoing null expansion of $\Sigma$ as in \cite{L2020}.

In the following, we present the geometry of $\ucalH$ in the coordinate system $\{s, \vartheta\}$ induced from the double null coordinate system $\{s, \us, \vartheta\}$ in the Schwarzschild spacetime. Let $\Sigma_s$ be the intersection of $\ucalH$ with the outgoing null hypersurface $C_s$. $\{\Sigma_s\}$ is a foliation of $\ucalH$, which is viewed as a background foliation.

As introduced in method I in subsection \ref{subsec 4.1}, assume that $\ucalH$ is parameterised by a function $\uh$ as its graph of $\us$ over $(s, \vartheta)$ domain in the double null coordinate system $\{s, \us, \vartheta\}$. Let $\dpartial_s, \dpartial_{i=1,2}$ be the coordinate frame vector of $\{s, \vartheta\}$ on $\ucalH$.

\begin{enumerate}[label=\ref{subsec 4.2}.\alph*.,leftmargin=.45in]
\item
The coordinate frame vectors of $\{s, \vartheta\}$ coordinate system on $\ucalH$ are
\begin{align*}
\dpartial_s = \partial_s + \dpartial_s \uh\cdot \partial_{\us},
\quad
\dpartial_i = \partial_i + \dpartial_i \uh\cdot \partial_{\us}.
\end{align*}
We use $\cdot$ on the top to indicate the corresponding notation being associated with $\Sigma_s$ or the $\{s, \vartheta\}$ coordinate system on $\ucalH$.

\item
Introduce the conjugate null frame $\{ \duL, \dL'\}$ associated with the background foliation $\{\Sigma_s\}$,
\begin{align*}
\left\{
\begin{aligned}
&
\dL'=L',
\\
&
\duL=\uL+\uvarepsilon L + \uvarepsilon^i \partial_i,
\end{aligned}
\right.
\end{align*}
where
\begin{align*}
\uvarepsilon=-\Omega^2 (\slashg^{-1} )^{ij} \uh_i \uh_j = -\Omega^2 \vert \dslashd \uh \vert_{\subslashg}^2,
\quad
\uvarepsilon^i = -2\Omega^2 ( \slashg^{-1} )^{ij} \uh_j.
\end{align*}

\item
The shifting vector $\db$ between $\dpartial_s$ and $\duL$ is given by
\begin{align*}
\duL = \dpartial_s + \db^i \dpartial_i
\quad
\Rightarrow
\quad
\db^i 
= 
- 2 \Omega^2 ( \slashg^{-1} )^{ij} \uh_j.
\end{align*}

Let $\dslashg$ be the intrinsic metric on $\Sigma_s$ and $\dslashepsilon$ be the volume form of $\dslashg$, then 
\begin{align*}
\dslashg_{ij} = \slashg_{ij} = r^2 \circg_{ij},
\quad
( \dslashg^{-1} )^{ij} = ( \slashg^{-1} )^{ij} = r^{-2} (\circg^{-1})^{ij},
\quad
\dslashepsilon_{ij} = \slashepsilon_{ij} = r^2 \circepsilon_{ij}.
\end{align*}
The degenerated metric on $\ucalH$ in $\{s, \vartheta\}$ coordinate system is
\begin{align*}
g\vert_{\ucalH} =  \dslashg_{ij} \big( \dot{\ed} \theta^i - \db^i \dot{\ed} s  \big) \otimes \big( \dot{\ed} \theta^j - \db^j \dot{\ed} s \big).
\end{align*}

\item
The connection coefficients on $\ucalH$ relative to the background foliation $\{\Sigma_s\}$ and $\{\duL, \dL' \}$ are given by the following formulae:
\begin{align}
\begin{aligned}
\dchi'_{ij} & = \chi'_{ij} = \frac{1}{2} \tr \chi' \slashg_{ij},  \quad \dtr \dchi' = \tr \chi',
\\
\duchi_{ij}
& =
\uchi_{ij} 
-\Omega^2 \vert \dslashd \uh \vert^2_{\subslashg} \chi_{ij}  
-2\Omega^2 \slashnabla_{ij}^2\uh
- 4\omega\Omega^2 (\dslashd \uh \otimes \dslashd \uh )_{ij} 
+ 2 \tr \chi \Omega^2 (\dslashd \uh \otimes \dslashd \uh )_{ij} ,
\\
\dtr \duchi & = \left( \dslashg^{-1} \right)^{ij} \duchi_{ij}
=  
\tr \uchi - 2\Omega^2 \slashDelta \uh + \Omega^2 \tr \chi  \vert \dslashd \uh \vert_{\subslashg}^2
- 4 \Omega^2 \omega \vert \dslashd \uh \vert^2_{\subslashg} ,
\\
\deta_i & = \frac{1}{2} \tr \chi \, ( \dslashd \uh )_i,
\\
\duomega 
& =
\uomega - \frac{1}{2} \Omega^2 \tr \chi \vert \dslashd \uh \vert^2_{\subslashg}.
\end{aligned}
\label{eqn 4.3}
\end{align}
In the above formulae, we use $\cdot$ to denote the inner product relative to $\slashg$, and $\dtr$ to denote the trace relative to $\dslashg= \slashg$. $\dslashd$ is differential operator on $\Sigma_s$. $\slashnabla$ in $\slashnabla \uh, \slashnabla^2_{ij} \uh$ is the pull back of the covariant derivative $\slashnabla$ of $(\Sigma_{\us,s}, \slashg)$ to $\Sigma_s$.\footnote{Here we abuse the notation $\slashnabla$ to denote both the covariant derivative of $(\Sigma_{\us,s}, \slashg)$ and its pull back to $\Sigma_s$. Which meaning $\slashnabla$ represents in a concrete formula depends on where the differentiated object is defined. For example, if a vector field $V$ is defined on $\Sigma_s$, then $\slashnabla$ in $\slashnabla V$ is the pull back to $\Sigma_s$. If $\slashnabla$ can be interpreted in both way in a formula, we will state the precise meaning of $\slashnabla$ to avoid ambiguity.  \label{footnote 10}}
$\slashDelta$ in $\slashDelta \uh$ is the operator $\left( \slashg^{-1} \right)^{ij} \slashnabla^2_{ij}$.

The precise meaning of the pull back $\slashnabla$ to $\Sigma_s$ is as follows: let $\slashGamma_{ij}^k$ be the Christoffel symbol of the covariant connection $\slashnabla$ of $(\Sigma_{\us,s}, \slashg)$, then
\begin{align*}
&
\text{$\phi$: a function on $\Sigma_s$,}
&&
\slashnabla_i \phi = \dpartial_i \phi, \quad \slashnabla^i \phi= ( \slashg^{-1} )^{ij} \dpartial_j \phi,
\\
&
\text{$V$: a vector field on $\Sigma_s$,}
&&
\slashnabla_i V^k =  \dpartial_i V^k + \slashGamma_{ij}^k V^j,
\\
&
\text{$T$: a tensor field on $\Sigma_s$,}
&&
\slashnabla_{i} T_{i_1 \cdots i_k}^{j_1 \cdots j_l}
=
\dpartial_i T_{i_1 \cdots i_k}^{j_1 \cdots j_l} 
-  \slashGamma_{i i_m}^{r}  T_{i_1 \cdots \underset{\hat{i_m}}{r}\cdots i_k}^{j_1 \cdots j_l}
+ \slashGamma_{i s}^{j_n} T_{i_1 \cdots i_k}^{j_1 \cdots \overset{\hat{j_n}}{s}\cdots  j_l}.
\end{align*}

Note that since $\slashg = r^2 \circg$, the covariant derivative $\slashnabla$ is simply the covariant derivative $\circnabla$ of the standard round sphere metric $\circg$. Their Christoffel symbols are equal, $\slashGamma= \circGamma$.

\item The curvature components on $\ucalH$ relative to the background foliation $\{\Sigma_s\}$ and $\{ \duL, \dL' \}$ are given by the following formulae:

\begin{align} 
\begin{aligned}
\dualpha_{ij} 
&=
8\Omega^4\uh_i \uh_j  \rho 
+2\Omega^2 \uvarepsilon \slashg_{ij}  \rho
+\uvarepsilon^k \uvarepsilon^l \slashepsilon_{ki} \slashepsilon_{jl} \rho,
\\
\dubeta_i
&=
-3 \Omega^2 \uh_i \rho,
\\
\dsigma \cdot \slashepsilon_{ij} 
&=
0,
\\
\drho 
&=
\rho,
\\ 
\dbeta_i 
&=
0,
\\
\dalpha_{ij} 
&=
0
\end{aligned}
\label{eqn 4.4}
\end{align}

\item
Denote the covariant derivative of $(\Sigma_s, \dslashg)$ by $\dslashnabla$. Let $\dslashGamma_{ij}^k$ be the Christoffel symbol of $\dslashnabla$. It is given by the following formula,
\begin{align*}
\dslashGamma_{ij}^k 
& = 
\slashGamma_{ij}^k 
+ ( \slashg^{-1} )^{kl} \big( \dpartial_i \uh \cdot \chi_{jl} + \dpartial_j \uh \cdot \chi_{il} - \dpartial_l \uh \cdot \chi_{ij}  \big)
\\
&=
\circGamma_{ij}^k 
+ \frac{1}{2} \tr \chi ( \slashg^{-1} )^{kl} \big( \dpartial_i \uh \cdot \slashg_{jl} + \dpartial_j \uh \cdot \slashg_{il} - \dpartial_l \uh \cdot \slashg_{ij}  \big).
\end{align*}
Introduce the tensor $\triangle_{ij}^k$ to denote the difference of $\dslashGamma_{ij}^k$ with $\slashGamma_{ij}^k$,
\begin{align*}
\triangle_{ij}^k = \frac{1}{2} \tr \chi ( \slashg^{-1} )^{kl} \big( \dpartial_i \uh \cdot \slashg_{jl} + \dpartial_j \uh \cdot \slashg_{il} - \dpartial_l \uh \cdot \slashg_{ij}  \big).
\end{align*}

\end{enumerate}

\subsection{Geometry of a spacelike surface}\label{subsec 4.3}
Let $\Sigma$ be a spacelike surface. In subsection \ref{subsec 4.2}, we already obtain the geometry of the incoming null hypersurface $\ucalH$ where $\Sigma$ is embedded in. Suppose $\Sigma$ is parameterised by a function $f$ as its graph of $s$ over $\vartheta$ domain in the coordinate system $\{s, \vartheta\}$ on $\ucalH$. In the following, we present the geometry of $\Sigma$ in terms of its parameterisation function $f$ and the geometry of $\ucalH$.

\begin{enumerate}[label=\ref{subsec 4.3}.\alph*.,leftmargin=.45in]
\item
Let $\{\vartheta\}$ be the coordinate system on $\Sigma$ which is induced from the coordinate system $\{s, \vartheta\}$ on $\ucalH$. The coordinate frame vector of $\{\vartheta\}$ on $\Sigma$ is given by
\begin{align*}
\ddpartial_i = \dpartial_i + f_i \dpartial_s = \dB_i^j \dpartial_j + f_i \duL,
\quad
\dB_i^j = \delta_i^j - f_i \db^j.
\end{align*}
We use $\cdot\cdot$ on the top to indicate the corresponding notation being associated with $\Sigma$.

\item
Let $\ddslashg$ be the intrinsic metric on $\Sigma$,
\begin{align*}
\ddslashg_{ij} 
=
\dB_i^k \dB_i^l \dslashg_{kl}
=
\slashg_{ij} - \big( \slashg_{ik} \db^i f_j + \slashg_{jl} \db^l f_i \big) + f_i f_j \slashg_{kl} \db^k \db^l.
\end{align*}
Let $\ddslashepsilon$ be the volume form of $\ddslashg$.

\item
Introduce the conjugate null frame $\{ \dduL, \ddL' \}$ on $\Sigma$
\begin{align*}
\left\{
\begin{aligned}
&
\ddL' = \dL' + \dvarepsilon' \duL + \dvarepsilon'^i \dpartial_i,
\\
&
\dduL = \duL,
\end{aligned}
\right.
\end{align*}
where
\begin{align*}
\dvarepsilon'^i
=
-2 ( \dslashg^{-1} )^{ik} \big( \dB^{-1} \big)_k^j f_j,
\quad
\dvarepsilon'
=
-\vert \ddslashd f \vert_{\subddslashg}^2
=
- ( \ddslashg^{-1} )^{ij} f_i f_j.
\end{align*}
$\ddslashd$ is the differential operator on $\Sigma$.

\item 
$\dslashnabla$ is the covariant derivative of $(\Sigma_s, \dslashg)$, and we also use it to denote the pull back of $\dslashnabla$ to $\Sigma$. The precise interpretation of $\dslashnabla$ shall be understood in the context.\footnote{This is similar to the pull back of $\slashnabla$. See footnote \ref{footnote 10}. Which meaning $\dslashnabla$ is interpreted as depends on where the differentiated object is defined. For example, if $V$ is a vector field defined on $\Sigma$, then $\dslashnabla$ in $\dslashnabla V$ should be interpreted as the pull back to $\Sigma$. If it can be interpreted in both way, we will point out the precise meaning of $\slashnabla$ to avoid ambiguity.} 
We have the following formulae for the pull back of $\dslashnabla$ to $\Sigma$
\begin{align*}
&
\text{$\phi$: a function on $\Sigma$,}
&&
\dslashnabla_i \phi = \ddpartial_i \phi, \quad \dslashnabla^i \phi= ( \dslashg^{-1} )^{ij} \ddpartial_j \phi = ( \slashg^{-1} )^{ij} \ddpartial_j \phi ,
\\
&
\text{$V$: a vector field on $\Sigma$,}
&&
\dslashnabla_i V^k =  \ddpartial_i V^k + \dslashGamma_{ij}^k V^j,
\\
&
\text{$T$: a tensor field on $\Sigma$,}
&&
\dslashnabla_{i} T_{i_1 \cdots i_k}^{j_1 \cdots j_l}
=
\ddpartial_i T_{i_1 \cdots i_k}^{j_1 \cdots j_l} 
-  \dslashGamma_{i i_m}^{r}  T_{i_1 \cdots \underset{\hat{i_m}}{r}\cdots i_k}^{j_1 \cdots j_l}
+ \dslashGamma_{i s}^{j_n} T_{i_1 \cdots i_k}^{j_1 \cdots \overset{\hat{j_n}}{s}\cdots  j_l}.
\end{align*}

\item
The connection coefficients on $\Sigma$ relative to the conjugate null frame $\{\dduL, \ddL'\}$ are given by the following formulae:
\begin{align}
\begin{aligned}
\dduchi_{ij}
&=
\duchi_{ij} + 2 \sym  \big\{   [-\duchi (\db )]  \otimes \ddslashd f  \big\}_{ij} + \duchi (\db,\db ) f_i f_j,
\quad
\ddtr \dduchi
=
\big(\ddslashg^{-1} \big)^{ij} \dduchi_{ij},
\\
\ddchi'_{ij} 
&=
\dchi'_{ij} + \dvarepsilon' \duchi_{ij} +  \big( \db \circdot \vec{\dvarepsilon}' -2 \big) \dslashnabla^2_{ij} f
\\
&\phantom{=}
+2\sym  \Big\{  
\ddslashd f \otimes  \big[ \dslashnabla \db \circdot \vec{\dvarepsilon}' - \duchi (\vec{\dvarepsilon}' ) -\dvarepsilon' \duchi (\db ) - \dchi' (\db ) -2 \deta  \big]
\Big\}_{ij}
\\
&\phantom{=}
+  \Big[
2\duchi (\db,\vec{\dvarepsilon}' ) +\dvarepsilon' \duchi (\db,\db ) +\dchi' (\db,\db ) +4\deta (\db ) 
\\
&
\phantom{= +  \Big[2\duchi (\db,\vec{\dvarepsilon}' ) +\dvarepsilon' \duchi (\db,\db ) }
-\dslashnabla_{\db} \db \circdot \vec{\dvarepsilon}' -\dpartial_s \db \circdot \vec{\dvarepsilon}' -4\duomega
\Big] f_i f_j,
\\
\ddtr \ddchi'
&=
\big(\ddslashg^{-1} \big)^{ij} \ddchi'_{ij},
\\
\ddeta_i &= \deta_i + \frac{1}{2} \duchi (\vec{\dvarepsilon}' )_i + \big[ 2\duomega-\deta (\db ) -\frac{1}{2} \duchi (\db,\vec{\dvarepsilon}' )  \big] f_i.
\end{aligned}
\label{eqn 4.5}
\end{align}
In the above formulae, we use $\circdot$ to denote the inner product relative to $\dslashg$, where
\begin{align*}
 ( \dslashnabla \db \circdot \vec{\dvarepsilon}'  )_i = \dslashg_{kl} \cdot \dslashnabla_i \db^k \cdot \dvarepsilon'^l,
\quad
\dslashnabla_{\db} \db \circdot \vec{\dvarepsilon}' = \db^i \cdot \dslashnabla_i \db^k \cdot \dvarepsilon'^l \cdot \dslashg_{kl},
\quad
\dpartial_s \db \circdot \vec{\dvarepsilon}' = \dslashg_{ij} \cdot \dpartial_s \db^i \cdot \dvarepsilon'^j.
\end{align*}
Numerically, $\circdot$ is the same as the inner product relative to $\slashg$, thus we will simply write $\cdot$ for $\circdot$ later. We use $\ddtr$ to denote the trace relative to $\ddslashg$.

\item The curvature components on $\Sigma$ relative to $\{\dduL, \ddL'\}$ are given by the following formulae:
\begin{align}
\begin{aligned}
\ddualpha_{ij}&=\dB_i^k \dB_j^l \dualpha_{kl},
\\
\ddubeta_i &= \dB_i^j \dubeta_j -\frac{1}{2} \dB_i^j \dvarepsilon'^k \dualpha_{jk},
\\ 
\ddrho &= \drho - \dvarepsilon'^j\dubeta_j + \frac{1}{4} \dvarepsilon'^i \dvarepsilon'^j \dualpha_{ij},
\\ 
\ddsigma \cdot \ddslashepsilon_{ij}
&=
( \dB_j^k f_i -\dB_i^k f_j) \dubeta_k  
+\frac{1}{2}\dB_i^k \dB_j^l \dvarepsilon'^m \dslashepsilon_{lk} \dslashepsilon_{mn} \ubeta^n 
+ \frac{1}{2} (f_j \dB_i^k - f_i \dB_j^k) \dvarepsilon'^m \dualpha_{km},
\\
\ddbeta_i
&=
- 3 f_i \drho 
- \dB_i^j \dvarepsilon' \dubeta_j 
+ 2 f_i \dvarepsilon'^l \dubeta_l 
+ \frac{1}{2} \dB_i^j \dvarepsilon'^k \dvarepsilon'^l \dslashepsilon_{kj} \dslashepsilon_{lm} \dubeta^m
\\
&\phantom{=}
+ \frac{1}{2} \dB_i^j \dvarepsilon'^l  \dvarepsilon' \dualpha_{jl}
-\frac{1}{2} f_i \dvarepsilon'^k \dvarepsilon'^l \dualpha_{kl},
\\
\ddalpha_{ij}
&=
-  (f_j \dB_i^k\dvarepsilon'^n + f_i \dB_j^k \dvarepsilon'^n) \dslashg_{nk} \drho
+4 f_i f_j \drho 
+ 2\dvarepsilon' \dB_i^k \dB_j^l \dslashg_{kl} \drho 
+ \dB_i^k \dB_j^l \dvarepsilon'^m \dvarepsilon'^n \dslashepsilon_{mk} \dslashepsilon_{ln} \drho
\\
&\phantom{=}
- 4 f_i f_j \dvarepsilon'^n \dubeta_n
+ 2 (\dB_j^k f_i +\dB_i^k  f_j ) \dvarepsilon' \dubeta_k 
+ (\dB_i^k \dB_j^l + \dB_j^k \dB_i^l) \dvarepsilon' \dvarepsilon'^n \dslashepsilon_{nl} \dslashepsilon_{km} \dubeta^m
\\
&\phantom{=}
+ (\dB_i^k f_j + \dB_j^k f_i ) \dvarepsilon'^m \dvarepsilon'^n \dslashepsilon_{km} \dslashepsilon_{nl} \dubeta^l
+ \dB_i^k \dB_j^l \dvarepsilon'^2 \dualpha_{kl}
+ f_i f_j \dvarepsilon'^m \dvarepsilon'^n \dualpha_{mn}
\\
&\phantom{=}
- (\dB_i^k f_j + \dB_j^k  f_i) \dvarepsilon' \dvarepsilon'^n \dualpha_{nk}.
\end{aligned}
\label{eqn 4.6}
\end{align}

\end{enumerate}

\section{Perturbation and linearised perturbation of parameterisation of spacelike surface}
In this section, we explain how to use the parameterisations in subsection \ref{subsec 4.1} to describe a perturbation of the coordinate surface $\Sigma_{s, \us=0}$. Then we give an appropriate linearised perturbation of parameterisations for a perturbation of $\Sigma_{s,\us=0}$. It should be pointed out that the constructions in this section are consistent with the ones in \cite{L2018} \cite{L2020}.

\subsection{Perturbation of parameterisation}
Let $\Sigma$ be a perturbation of the coordinate surface $\Sigma_{s,\us=0}$. We use the parameterisations in subsection \ref{subsec 4.1} to describe this perturbation.

Suppose that the first parameterisation of $\Sigma$ is $(f,\uf)$, and its second parameterisation is $(f, \ufl{s=0})$. Then we define the following perturbation functions of the parameterisations between $\Sigma_{s,\us=0}$ and $\Sigma$:
\begin{align*}
\dd{f} = f - s,
\quad
\dd{\uf}= \uf -0 = \uf,
\quad
\dd{\ufl{s=0}} = \ufl{s=0} - 0 = \ufl{s=0}.
\end{align*}

In subsection \ref{subsec 4.1}, we described how to transform the second parameterisation from the first. We can also use the above mentioned transformation to obtain $\dd{\uf}$ from $\dd{\ufl{s=0}}$. In fact, given method II in subsection \ref{subsec 4.1}, the procedure to obtain $\dd{\uf}$ from $\dd{\ufl{s=0}}$ is trivial, since $\dd{\uf}= \uf$ and $\dd{\ufl{s=0}} = \ufl{s=0}$, then equation \eqref{eqn 4.2} gives the result. 

\subsection{Linearised perturbation of parameterisation}\label{subsec 5.2}
We consider the linearised perturbation of parameterisation at $\Sigma_{s,\us=0}$. Denote the linearised perturbations of parameterisation functions by
\begin{align*}
\bdd{f},
\quad
\bdd{\ufl{s=0}},
\quad
\bdd{\uf}.
\end{align*}
In previous subsection, we learnt how to transform the second parameterisation from the first. In this subsection, we want to achieve a similar transformation for $\bdd{\ufl{s=0}}$ and $\bdd{\uf}$. Moreover, we require such a transformation being linear.

First we define linearised perturbations of the first parameterisation functions simply by
\begin{align*}
\bdd{f} = \dd{f},
\quad
\bdd{\ufl{s=0}} = \dd{\ufl{s=0}}.
\end{align*}
In order to obtain a linear transformation from $\bdd{\ufl{s=0}}$ to $\bdd{\uf}$, we consider the linearisation of equation \eqref{eqn 4.2}
\begin{align}
\partial_t \ufl{t} 
=
f \cdot [ 1- t\underline{e}^i f_i - t e^i f_i \cdot \uvarepsilon ]^{-1} \cdot \left[ \uvarepsilon - \underline{e}^i \, \ufl{t}_i - e^i \, \ufl{t}_i \cdot \uvarepsilon \right].
\tag{\ref{eqn 4.2}}
\end{align}
Note that the right hand side of above equation is quadratic in $\ufl{t}$, thus the linearisation of the right hand side at $\ufl{t}=0$ shall be zero. Hence the linearisation of $\ufl{t}$ is constant for all $t$, thus the transformation from $\bdd{\ufl{s=0}}$ to $\bdd{\uf}$ is simply the identity map:
\begin{align*}
\bdd{\uf} = \bdd{\ufl{s=0}} = \dd{\ufl{s=0}}.
\end{align*}

\section{Linearised perturbation of geometry of spacelike surface at $\Sigma_{s,\us=0}$}

In this section, we construct a linearised perturbation for the geometry of spacelike surfaces at $\Sigma_{s,\us}$. We shall use the formulae of geometric quantities obtained in subsections \ref{subsec 4.2} and \ref{subsec 4.3} as important tools in the construction.

\subsection{Sketch of construction of linearised perturbation}
Recalling the procedure to obtain the geometry of a spacelike surface $\Sigma$ in subsections \ref{subsec 4.2}, \ref{subsec 4.3}, we first study the geometry of the incoming null hypersurface $\ucalH$, then obtain the geometry of $\Sigma$ by viewing it as an embedded surface in $\ucalH$.

In the construction of the linearised perturbation of the geometry of spacelike surfaces, we adopt a similar procedure.
\begin{enumerate}[label=\alph*.]
\item Construct the linearised perturbations of connection coefficients $\dchi'$, $\duchi$, $\deta$, $\duomega$ (formulae \eqref{eqn 4.3}) and curvature components $\dalpha$, $\dualpha$, $\dbeta$, $\dubeta$, $\drho$, $\dsigma$ (formulae \eqref{eqn 4.4}) on incoming null hypersurfaces. 

\item Construct the linearised perturbation of connection coefficients $\ddchi'$, $\dduchi$, $\ddeta$ (formulae \eqref{eqn 4.5}) and curvature components $\ddalpha$, $\ddualpha$, $\ddubeta$, $\ddbeta$, $\ddrho$, $\ddsigma$ (formulae \eqref{eqn 4.6}) on spacelike surfaces.
\end{enumerate}

Note that in order to achieve step a, we need to construct the linearised perturbation of the differential $\dslashd \uh$ and Hessian $\slashnabla^2 \uh$ of the parameterisation function $\uh$ on spacelike surfaces. Thus we will first construct them, then turn to step a and b mentioned above.

\subsection{Linearised perturbations of differential and Hessian of $\uh$ at $\Sigma_{s,\us=0}$}\label{subsec 6.2}

In order to construct the linearised perturbations of $\dslashd \uh$, $\slashnabla^2 \uh$, we consider equation \eqref{eqn 4.1} of $\uh$,
\begin{align}
\dpartial_s \uh = \Omega^2 ( \slashg^{-1} )^{ij} \dpartial_i \uh\, \dpartial_j \uh.
\tag{\ref{eqn 4.1}}
\end{align}
$\Sigma_{s, \us=0}$ is embedded in the incoming surface $\uC_{\us=0}$, whose parameterisation function $\uh_{\uC_{\us=0}} = 0$. Since the right hand side of equation \eqref{eqn 4.1} is quadratic in $\uh$, thus the linearisation of the right hand side of equation \eqref{eqn 4.1} at $\uh_{\uC_{\us=0}}$ is simply zero. We denote the linearisation of the parameterisation function $\uh$ by $\bdd{\uh}$, then the linearisation $\bdd{\uh}$ at $\uh_{\uC_{\us=0}}$ satisfies the equation
\begin{align*}
\dpartial_s \bdd{\uh} = 0,
\end{align*}
with the initial condition $\bdd{\uh}(s=0,\vartheta) = \bdd{\ufl{s=0}}(\vartheta)$. Hence the linearisation $\bdd{\uh}$ is simply
\begin{align*}
\bdd{\uh}(s, \vartheta) = \bdd{\ufl{s=0}}(\vartheta).
\end{align*}
Then we define the linearisations of $\dslashd \uh$ and $\slashnabla^2 \uh$ by
\begin{align}
\bdd{\dslashd \uh} = \slashd \bdd{\uh} = \slashd \bdd{\ufl{s=0}},
\quad
\bdd{\slashnabla^2 \uh} = \circnabla^2 \bdd{\uh} = \circnabla^2 \bdd{\ufl{s=0}}.
\label{eqn 6.1}
\end{align}
Written in components, the above can be rewritten as
\begin{align}
\bdd{(\dslashd \uh)_i} = (\slashd \bdd{\uh})_i = (\slashd \bdd{\ufl{s=0}})_i,
\quad
\bdd{\slashnabla^2_{ij} \uh} = \circnabla^2_{ij} \bdd{\uh} = \circnabla^2_{ij} \bdd{\ufl{s=0}}.
\tag{\ref{eqn 6.1}$'$}
\label{eqn 6.1'}
\end{align}

\subsection{Linearised perturbation of geometry of null hypersurface at $\Sigma_{s,\us=0}$}\label{subsec 6.3}

In this subsection, we shall construct linearised perturbation of the geometry of null hypersurfaces at $\Sigma_{s,\us=0}$ by using the formulae in subsection \ref{subsec 4.2}, the linearised perturbation $\bdd{f}$ in subsection \ref{subsec 5.2}, and the linearised perturbations $\bdd{\dslashd \uh}$, $\bdd{\slashnabla^2 \uh}$ in subsection \ref{subsec 6.2}. More precisely, let $\Sigma$ be a spacelike surface embedded in an incoming null hypersurface $\ucalH$, then we are interested in the linearised perturbation from the geometry of $\uC_{\us=0}$ restricted on $\Sigma_{s, \us=0}$ to the geometry of $\ucalH$ restricted on $\Sigma$.

In the following, we use $\bdd{a}$ to denote the linearised perturbation of some quantity $a$ from $\Sigma_{s,\us=0}$ to $\Sigma$.

\begin{enumerate}[label=\ref{subsec 6.3}.\alph*.,leftmargin=.45in]
\item Suppose the linearised perturbations of the second parameterisation functions $f$, $\ufl{s=0}$ are $\bdd{f}$, $\bdd{\ufl{s=0}}$. Following subsections \ref{subsec 5.2}, \ref{subsec 6.2}, we have the linearised perturbation of the first parameterisation function $\uf$ is
\begin{align*}
\bdd{\uf} = \bdd{\ufl{s=0}},
\end{align*}
and the linearised perturbations of $\dslashd \uh$, $\slashnabla^2 \uh$ are
\begin{align*}
\bdd{\dslashd \uh} = \slashd \bdd{\ufl{s=0}},
\quad
\bdd{\slashnabla^2 \uh} = \circnabla^2 \bdd{\ufl{s=0}}.
\end{align*}

\item The linearised perturbations of $\uvarepsilon$ and $\vec{\uvarepsilon}$ in subsection \ref{subsec 4.2}.b. are
\begin{align*}
\bdd{\uvarepsilon} = 0,
\quad
\bdd{\uvarepsilon} = - 2 r^{-2} (\circg^{-1})^{ij} (\bdd{\uh})_j
\end{align*}

\item The linearised perturbation of the shifting vector $\db$ is
\begin{align*}
\bdd{\db^i} = - 2 r^{-2} (\circg^{-1})^{ij} (\bdd{\uh})_j.
\end{align*}
The linearised perturbation of the intrinsic metric $\dslashg$ is given by
\begin{align*}
&
\bdd{r}
=
\partial_{s} r \cdot \bdd{f} + \partial_{\us} r \cdot \bdd{\uf}
=
\bdd{f} + \frac{s}{r} \bdd{\uf},
\\
&
\bdd{\dslashg}
=
\bdd{\slashg}
=
2r\bdd{r} \circg
=
2(r\bdd{f} + s \bdd{\uf}) \circg.
\end{align*}

\item The linearised perturbations of connection coefficients $\dchi', \duchi, \deta, \duomega$ are given by
\begin{align*}
\begin{aligned}
&
\left\{
\begin{aligned}
\bdd{\dtr \dchi'}
& = 
\bdd{\tr \chi'} = \partial_s \tr \chi' \cdot \bdd{f} + \partial_{\us} \tr \chi' \cdot \bdd{\uf}
= \frac{2 (r_0- s)}{r^3} \bdd{f} - \frac{2s^2}{r^4} \bdd{\uf},
\\
\bdd{\dchi'} 
& =
\frac{1}{2} \bdd{\tr \chi'} r^2 \circg + \frac{1}{2} \tr \chi' \bdd{\slashg}
=
\bdd{f} \circg+ \frac{s^2}{r^2} \bdd{\uf} \circg,
\end{aligned}
\right.
\\
&
\left\{
\begin{aligned}
\bdd{\tr \uchi} 
& =
\partial_s \tr \uchi \cdot \bdd{f} + \partial_{\us} \tr \uchi \cdot \bdd{\uf}
=
-\frac{2}{r^2} \bdd{f} -\frac{2r-4r_0}{r^3} \bdd{\uf},
\\
\bdd{\uchi} 
& =
\frac{1}{2} \bdd{\tr \uchi} r^2 \circg + \frac{1}{2} \tr \uchi \bdd{\slashg}
=
\bdd{f} \circg
+
\bdd{\uf} \circg,
\\
\bdd{\dtr \duchi}
& =
\bdd{\tr \uchi} - 2\bdd{\slashDelta \uh}
=
-\frac{2}{r^2} \bdd{f}
- \frac{2r-4r_0}{r^3} \bdd{\uf}
- \frac{2}{r^2} \circDelta \bdd{\uh},
\\
\bdd{\duchi}
& =
\bdd{\uchi} - 2\bdd{\slashnabla^2 \uh}
=
\bdd{f} \circg
+
\bdd{\uf} \circg
- 
2\circnabla^2 \bdd{\uh},
\end{aligned}
\right.
\\
&\phantom{\Big\{}
\begin{aligned}
\bdd{\deta} 
&= 
\frac{1}{2} \tr \chi \bdd{ \dslashd \uh}
=
\frac{s}{r^2} \slashd \bdd{\uh}, 
\end{aligned}
\\
&
\phantom{\Big\{}
\begin{aligned}
\bdd{\duomega}
& =
\bdd{\uomega}
= 
\partial_s \uomega \cdot \bdd{f} + \partial_{\us} \uomega \cdot \bdd{\uf}
=
-\frac{r_0}{r^3} \bdd{\uf}.
\end{aligned}
\end{aligned}
\end{align*}

\item The linearised perturbations of curvature components $\dualpha, \dubeta, \dsigma, \drho, \dbeta, \dalpha$ are given by
\begin{align*}
\bdd{\dualpha} &= 0,
\\
\bdd{\dubeta} &= - 3\Omega^2 \rho \cdot \slashd \bdd{\uh} = \frac{3 r_0}{r^3} \slashd \bdd{\uh},
\\
\bdd{\dsigma} &= 0,
\\
\bdd{\drho} &= \bdd{\rho} 
= \partial_s \rho \cdot \bdd{f} + \partial_{\us} \rho \cdot \bdd{\uf} 
= \frac{3r_0}{r^4} \bdd{f} + \frac{3r_0 s}{r^5} \bdd{\uf},
\\
\bdd{\dbeta} &= 0,
\\
\bdd{\dalpha} &= 0.
\end{align*}

\item The linearised perturbation of the Christoffel symbol $\dslashGamma_{ij}^k$ is given by
\begin{align*}
\bdd{\dslashGamma_{ij}^k} = \bdd{\triangle_{ij}^k}
&=
 \frac{1}{2} \tr \chi ( \slashg^{-1} )^{kl} \big( \partial_i \bdd{\uh} \cdot \slashg_{jl} + \partial_j \bdd{\uh} \cdot \slashg_{il} - \partial_l \bdd{\uh} \cdot \slashg_{ij}  \big)
\\
&=
\frac{s}{r^2} \big( \partial_i \bdd{\uh} \cdot \delta_j^k + \partial_j \bdd{\uh} \cdot \delta_i^k - \partial_l \bdd{\uh} \cdot (\circg^{-1})^{kl}  \circg_{ij}  \big)
\end{align*}

\end{enumerate}

\subsection{Linearised perturbation of geometry of spacelike surface at $\Sigma_{s,\us=0}$}\label{subsec 6.4}
In this subsection, we shall construct the linearised perturbation of the geometry of spacelike surfaces at $\Sigma_{s,\us=0}$ by using the formulae of the geometry of spacelike surfaces in subsection \ref{subsec 4.3} and the linearised perturbations constructed in previous subsections.

\begin{enumerate}[label=\ref{subsec 6.4}.\alph*.,leftmargin=.45in]
\item Following subsection \ref{subsec 6.3}, assume that the linearised perturbations of the second parameterisation functions $f$, $\ufl{s=0}$ are $\bdd{f}$, $\bdd{\ufl{s=0}}$, then all the constructions in subsection \ref{subsec 6.3} apply and
\begin{align*}
\bdd{\dB_i^j} = \bdd{\delta_i^j} - \bdd{f_i \db^j} = 0
\end{align*}

\item The linearised perturbation of the intrinsic metric $\ddslashg$ is
\begin{align*}
\bdd{\ddslashg} = \bdd{\slashg} = 2 r \bdd{f} \circg + 2s \bdd{\uf} \circg.
\end{align*}

\item The linearised perturbations of $\dvarepsilon'$, $\dvarepsilon'^i$ in subsection \ref{subsec 4.3} are
\begin{align*}
\bdd{\dvarepsilon'^i} = - 2 r^{-2} (\circg^{-1})^{ij} (\bdd{f})_j,
\quad
\bdd{\dvarepsilon'} =0.
\end{align*}

\item The linearised perturbation of the Christoffel symbol $\dslashGamma_{ij}^k$ is presented previously in subsection \ref{subsec 6.3}.f.. 

\item The linearised perturbations of connection coefficients $\dduchi, \ddchi',\ddeta$ are given by
\begin{align*}
&
\left\{
\begin{aligned}
\bdd{\dduchi}
&=
\bdd{\duchi}
=
\bdd{f} \circg
+
\bdd{\uf} \circg
- 
2\circnabla^2 \bdd{\uh},
\\
\bdd{\ddtr \dduchi}
&=
\bdd{\dtr \duchi}
=
-\frac{2}{r^2} \bdd{f}
- \frac{2r-4r_0}{r^3} \bdd{\uf}
- \frac{2}{r^2} \circDelta \bdd{\uh},
\end{aligned}
\right.
\\
&
\left\{
\begin{aligned}
\bdd{\ddchi'}
&=
\bdd{\dchi'} - 2 \circnabla^2 \bdd{f}
=
\bdd{f} \circg+ \frac{s^2}{r^2} \bdd{\uf} \circg  - 2 \circnabla^2 \bdd{f},
\\
\bdd{\ddtr \ddchi'}
&=
\bdd{\dtr \dchi'}
-
\frac{2}{r^2} \circDelta \bdd{f}
=
\frac{2(r_0-s)}{r^3} \bdd{f} 
- \frac{2s^2}{r^4} \bdd{\uf} 
- \frac{2}{r^2} \circDelta \bdd{f}
\end{aligned}
\right.
\\
&\phantom{\Big\{}
\begin{aligned}
\bdd{\ddeta_i} 
&= 
\bdd{\deta_i}  + \frac{1}{2} \uchi_{ij} \bdd{\dvarepsilon'^j}
=
\frac{s}{r^2} (\bdd{\uh})_i - \frac{1}{r} (\bdd{f})_i.
\end{aligned}
\end{align*}

\item The linearised perturbations of curvature components $\ddualpha, \ddubeta, \ddrho, \ddsigma, \ddbeta, \ddalpha$ are given by
\begin{align*}
\bdd{\ddualpha}
&=
0,
\\
\bdd{\ddubeta}
&=
\bdd{\dubeta} 
=
\frac{3 r_0}{r^3} \slashd \bdd{\uh},
\\
\bdd{\ddrho}
&=
\bdd{\drho}
=
\frac{3r_0}{r^4} \bdd{f} + \frac{3 r_0 s}{r^5} \bdd{\uf},
\\
\bdd{\ddsigma}
&
=0,
\\
\bdd{\ddbeta}
&=
-3 \rho \slashd \bdd{f}
=
\frac{3r_0}{r^3} \slashd \bdd{f},
\\
\bdd{\ddalpha}
&=
0.
\end{align*}
\end{enumerate}

\section{Strategies to construct linearised perturbation of constant mass aspect function foliation at $\{\Sigma_{s,\us=0}\}_{s\geq 0}$}\label{sec 7}

Recall that the goal of this paper is to study the linearised perturbation of constant mass aspect function foliations in a Schwarzschild spacetime. In this section, we shall explain two strategies to construct the linearised perturbation of the foliations at the spherically symmetric constant mass aspect function foliation $\{\Sigma_{s, \us=0}\}_{s\geq 0}$ of $\uC_{\us=0}$.

\subsection{Constant mass aspect function foliation in Schwarzschild spacetime}\label{subsec 7.1}
Let $\ucalH$ be an incoming null hypersurface in the Schwarzschild spacetime, and suppose that $\ucalH$ is parameterised by $\uh$ as its graph of $\us$ over $(s,\vartheta)$ domain in $\{\us, s, \vartheta\}$ coordinate system. We studied the geometry of $\ucalH$ in subsection \ref{subsec 4.2}.

Assume that $\{\bSigma_u\}_{u\geq 0}$ is a constant mass aspect function foliation of $\ucalH$ as in definition \ref{def 2.6} and satisfies the parameterisation condition that $\br_{u_0+u} = \br_{u_0}+u$, where $\br_u$ is the area radius of $\bSigma_u$. $\{ \bL'^u, \buL^u\}$ is the conjugate null frame relative to $\{\bSigma_u\}$, that $\buL^u u=1$. 

Associated with this foliation $\{\bSigma_u\}_{u \geq 0}$, we have the following quantities on its parameterisation and geometry.
\begin{enumerate}[label=\ref{subsec 7.1}.\alph*.,leftmargin=.45in]

\item Parameterisation functions of $\bSigma_u$. Denote the parameterisation functions of the first parameterisation of $\bSigma_u$ by $\fl{u}$, $\ufl{u}$, and the parameterisation functions of the second parameterisation by $\fl{u}, \ufl{s=0}$.\footnote{Note that the first and second parameterisation share the same parameterisation function $\fl{u}$. Also note that $\{\bSigma_u\}$ share the same parameterisation function $\ufl{s=0}$ for the second parameterisation, since all $\bSigma_u$ are embedded in one null hypersurface $\ucalH$.} By the results in subsection \ref{subsec 4.3}, we can obtain the geometry of $\bSigma_u$ relative to the frame $\{\dduL, \ddL'\}$.

\item The lapse function $\bal{u}$ of $\{\bSigma_u\}$ relative to the induced background coordinate system $\{s, \vartheta\}$ on $\ucalH$.

\item The intrinsic metric $\bslashgl{u}$ on $\bSigma_u$.

\item The connection coefficients $\buchil{u}$, $\bchil{u}'$, $\btal{u}$, $\buomegal{u}$ of $\{ \bSigma_u\}$ relative to the frame $\{\bL'^u, \buL^u\}$.

\item The Gauss curvature $\bKl{u}$ of $\bSigma_u$.

\item The mass aspect function $\bmul{u}$ of $\{\bSigma_u\}$.

\item The curvature components $\bualphal{u}$, $\bubetal{u}$, $\brhol{u}$, $\bsigmal{u}$, $\bbetal{u}$, $\balphal{u}$ relative to the frame $\{ \bL'^u, \buL^u \}$.
\end{enumerate}

Recall that in subsections \ref{subsec 2.4}, \ref{subsec 2.5}, we introduced several equations useful for studying a constant mass aspect function foliation, which are equation \eqref{eqn 2.3} and the basic equations \eqref{eqn 2.4}-\eqref{eqn 2.12}. Additional to above equations, we have two more for the lapse function $\bal{u}$ which follow from $\nabla_{\buL^u} \buL^u = 2 \buomegal{u} \buL^u$ and $\overline{\btr \buchil{u}} = \frac{2}{\br_u}$,
\begin{align}
&
\buL^u \log \bal{u} = 2\, \buomegal{u} - 2\, \bal{u} \duomega,
\label{eqn 7.1}
\\
&
\overline{\bal{u}\, \ddtr \dduchi} = \frac{2}{\br_u},
\tag{\ref{eqn 2.2}}
\end{align}
where $\duomega$ is the background acceleration of the tangent null vector field $\duL$ on $\ucalH$ and $\ddtr \dduchi$ is the null expansion of $\bSigma_u$ relative to $\dduL$, see equation \eqref{eqn 4.3} in subsection \ref{subsec 4.2} and equation \eqref{eqn 4.5} in subsection \ref{subsec 4.3}.

\subsection{Two strategies to construct linearised perturbation of foliation}\label{subsec 7.2}
For the linearised perturbation of the constant mass aspect function foliation, we want to construct the linearised perturbations for not only the parameterisation functions $\bdd{\fl{u}}$, $\bdd{\ufl{u}}$, $\bdd{\ufl{s=0}}$, but also the quantities associated with the foliation listed in \ref{subsec 7.1}.a.-g.. Namely, we want to construct the following linearised perturbations.
\begin{enumerate}[label=\ref{subsec 7.2}.\alph*.,leftmargin=.45in]
\item Linearised perturbations of parameterisation functions: $\bdd{\fl{u}}$, $\bdd{\ufl{u}}$, $\bdd{\ufl{s=0}}$, $\bdd{\uh}$.

\item Linearised perturbation of the lapse function: $\bdd{\bal{u}}$.

\item Linearised perturbations of the intrinsic metric and the area radius: $\bdd{\bslashgl{u}}$, $\bdd{\br_u}$.

\item Linearised perturbations of the connection coefficients: $\bdd{\buchil{u}}$, $\bdd{\bchil{u}'}$, $\bdd{\btal{u}}$, $\bdd{\buomegal{u}}$.

\item Linearised perturbation of the Gauss curvature $\bdd{\bKl{u}}$ of $\bSigma_u$.

\item Linearised perturbation of the mass aspect function: $\bdd{\bmul{u}}$.

\item Linearised perturbations of the curvature components: $\bdd{\bualphal{u}}$, $\bdd{\bubetal{u}}$, $\bdd{\brhol{u}}$, $\bdd{\bsigmal{u}}$, $\bdd{\bbetal{u}}$, $\bdd{\balphal{u}}$.

\end{enumerate}

Among all the linearised perturbations, the most basic ones are the linearised perturbations of the parameterisation functions $\bdd{\fl{u=0}}$, $\bdd{\ufl{s=0}}$ of the initial leaf from $\Sigma_{s=0,\us=0}$ to $\bSigma_{u=0}$. The reason is the following: the constant mass aspect function foliation $\{\bSigma_u\}$ is completely determined by its initial leaf $\bSigma_{u=0}$. Thus all the linearised perturbations listed in \ref{subsec 7.2}.a.-f. should be derived from $\bdd{\fl{u=0}}$, $\bdd{\ufl{s=0}}$. The goal of this section is to achieve it.

There are two strategies to construct the above linearised perturbations:
\begin{enumerate}[label=\arabic*.]
\item In subsection \ref{subsec 2.4}, we formulate the construction of a constant mass aspect function foliation as an inverse lapse problem. Then one can use this formulation to construct the linearised perturbation of the foliation by linearising equations \eqref{eqn 2.2}, \eqref{eqn 2.3}.

\item In subsection \ref{subsec 2.5}, we introduce the basic equations for a constant mass aspect function foliation on its geometry. Then one can also construct the linearised perturbation of the foliation by linearising the basic equations.

\end{enumerate}

In the following, we shall employ both strategies and show that they result in the same linearised perturbation. We shall use $\bddsub{1}{\cdot}$ to denote the linearised perturbation constructed by the first strategy and $\bddsub{2}{\cdot}$ by the second strategy.

\section{Linearised perturbation of foliation at initial leaf $\Sigma_{s=0,\us=0}$}\label{sec 8}

In both strategies, we first need to construct the linearised perturbation of the initial leaf from $\Sigma_{s=0,\us=0}$ to $\bSigma_{u=0}$. We will construct it in this section.

\subsection{Construction of linearised perturbation of foliation at initial leaf $\Sigma_{s=0,\us=0}$}\label{subsec 8.1}

As already mentioned, the most basic linearised perturbations are $\bdd{\fl{u=0}}$ and $\bdd{\ufl{s=0}}$. Thus we shall assume that these two are known, then derive all the other linearised perturbations. 

We present the procedure to obtain all the other linearised perturbations from $\bdd{\fl{u=0}}$, $\bdd{\ufl{s=0}}$ in the following.

\begin{enumerate}[label=\ref{subsec 8.1}.\alph*.,leftmargin=0.45in]
\item $\bdd{\ufl{u=0}} = \bdd{\ufl{s=0}}$, $\bdd{\uh}= \bdd{\ufl{s=0}}$. See subsection \ref{subsec 5.2}.

\item $\bdd{\bal{u=0}}$ satisfies the following elliptic equation
\begin{align}
\left\{
\begin{aligned}
\frac{1}{r_0^2} \circDelta \bdd{\bal{u=0}}
&=
- \bdd{\ddrho}|_{\Sigma_{0,0}} 
+ \overline{\bdd{\ddrho}|_{\Sigma_{0,0}}} 
- \frac{1}{r_0^2} \circdiv \bdd{\ddeta}|_{\Sigma_{0,0}}
\\
&=
- \frac{3}{r_0^3} \big( \bdd{\fl{u=0}} - \overline{\bdd{\fl{u=0}}} \big)
+ \frac{1}{r_0^3} \circDelta \bdd{\fl{u=0}},
\\
\overline{\bdd{\bal{u=0}}}
&=
-\frac{1}{r_0} \bdd{\br_{u=0}} - \frac{r_0}{2} \overline{\bdd{\ddtr \dduchi}|_{\Sigma_{0,0}}}
=
-\frac{1}{r_0} \overline{\bdd{\ufl{s=0}}},
\end{aligned}
\right.
\label{eqn 8.1}
\end{align}
which follows from linearising equation \eqref{eqn 2.2}.

\item $\bdd{\bslashgl{u=0}} = \bdd{\ddslashg}|_{\Sigma_{0,0}} = 2r_0 \bdd{\fl{u=0}} \circg$, $\bdd{\br_{u=0}} = \overline{\bdd{\fl{u=0}}}$.

\item Linearised perturbations of connection coefficients at the initial leaf,
\begin{align}
\begin{aligned}
&
\left\{
\begin{aligned}
\bdd{\buchil{u=0}}
&= 
\bdd{\bal{u=0}\, \dduchi}|_{\Sigma_{0,0}}
=
\bdd{\dduchi}|_{\Sigma_{0,0}} + \bdd{\bal{u=0}} \cdot \uchi|_{\Sigma_{0,0}}
\\
&=
\bdd{\fl{u=0}} \circg  
+ \bdd{\bal{u=0}} r_0 \circg
+ \big( \bdd{\ufl{s=0}} 
- \overline{\bdd{\ufl{s=0}}} \big) \circg 
- 2 \circnabla^2 \bdd{\ufl{s=0}},
\\
\bdd{\btr \buchil{u=0}}
&=
\bdd{\bal{u=0}\, \ddtr \dduchi}|_{\Sigma_{0,0}} 
\\
&=
-\frac{2}{r_0^2} \bdd{\fl{u=0}} 
+ \frac{2}{r_0} \bdd{\bal{u=0}}
+ \frac{2}{r_0^2} \big( \bdd{\ufl{s=0}} 
- \overline{\bdd{\ufl{s=0}}} \big) 
- \frac{2}{r_0^2} \circDelta \bdd{\ufl{s=0}},
\\
\bdd{\bhatuchil{u=0}}
&=
\bdd{\buchil{u=0} - \frac{1}{2} \btr \buchil{u=0} \slashgl{u=0}}
\\
&=
\bdd{\buchil{u=0}} 
- \frac{1}{2} \bdd{\btr \buchil{u=0}} r_0^2 \circg 
- \frac{1}{2} \tr \uchi|_{\Sigma_{0,0}} \bdd{\bslashgl{u=0}}
=
-2 \widehat{\circnabla^2} \bdd{\ufl{s=0}},
\end{aligned}
\right.
\\
&
\left\{
\begin{aligned}
\bdd{\bchil{u=0}'}
&= 
\bdd{\bal{u=0}^{-1}\, \ddchi'}|_{\Sigma_{0,0}}
=
\bdd{\ddchi'}|_{\Sigma_{0,0}} + \bdd{\bal{u=0}^{-1}} \cdot \chi'|_{\Sigma_{0,0}}
\\
&=
\bdd{\fl{u=0}} \circg - 2\circnabla^2 \bdd{\fl{u=0}},
\\
\bdd{\btr \bchil{u=0}'}
&=
\bdd{\bal{u=0}^{-1}\, \ddtr \ddchi'}|_{\Sigma_{0,0}} 
=
\frac{2}{r_0^2} \bdd{\fl{u=0}} - \frac{2}{r_0^2} \circDelta \bdd{\fl{u=0}},
\\
\bdd{\bhatchil{u=0}'}
&=
\bdd{\bchil{u=0}' - \frac{1}{2} \btr \bchil{u=0}' \slashgl{u=0}}
\\
&=
\bdd{\bchil{u=0}'} 
- \frac{1}{2} \bdd{\btr \bchil{u=0}'} r_0^2 \circg 
- \frac{1}{2} \tr \chi'|_{\Sigma_{0,0}} \bdd{\bslashgl{u=0}}
=
-2 \widehat{\circnabla^2} \bdd{\fl{u=0}},
\end{aligned}
\right.
\\
&\phantom{\Big\{}
\begin{aligned}
\bdd{\btal{u=0}}
&= 
\slashd \bdd{\log \bal{u=0}}
+ \bdd{\ddeta}|_{\Sigma_{0,0}}
= 
\slashd \bdd{\bal{u=0}} - \frac{1}{r_0} \slashd \bdd{f},
\end{aligned}
\end{aligned}
\label{eqn 8.2}
\end{align}
where $\widehat{\circnabla^2}$ is the trace-free part of the Hessian operator that
$\widehat{\circnabla^2} = \circnabla^2 - \frac{1}{2} \circg \circDelta$.

For $\buomegal{u=0}$, there are two kinds of linearised perturbation by two different strategies. By linearising equation \eqref{eqn 7.1}, the first strategy gives
\begin{align}
\begin{aligned}
\bddsub{1}{\buomegal{u=0}}
&= 
\frac{1}{2}\frac{\ed}{\ed u} \bdd{\bal{u}}|_{u=0} + \bdd{\duomega}|_{\Sigma_{0,0}}
=
\frac{1}{2}\frac{\ed}{\ed u} \bdd{\bal{u}}|_{u=0} - \frac{1}{r_0^2} \bdd{\ufl{s=0}}.
\end{aligned}
\label{eqn 8.3}
\end{align}
The second strategy implies that $\bddsub{2}{\buomegal{u=0}}$ satisfies the following elliptic equation
\begin{align}
\left\{
\begin{aligned}
\circDelta \bddsub{2}{\buomegal{u=0}}
&=
- \frac{3 r_0^2}{4} \mu|_{\Sigma_{0,0}} \big( \bdd{\btr \buchil{u=0}} - \overline{ \bdd{\btr \buchil{u=0}} } \big) - \circdiv \bdd{\bubetal{u=0}}
\\
&=
\frac{3}{2 r_0^2} \big(  \bdd{\fl{u=0}} -\overline{\bdd{\fl{u=0}}} \big) 
- \frac{3}{2r_0} \bdd{\bal{u=0}}
\\
&\phantom{=}
- \frac{3}{2r_0^2} \big( \bdd{\ufl{s=0}} - \overline{ \bdd{\ufl{s=0}}} \big) 
- \frac{3}{2r_0^2} \circDelta \bdd{\ufl{s=0}},
\\
\overline{\bddsub{2}{\buomegal{u=0}}}
&=
0,
\end{aligned}
\right.
\label{eqn 8.4}
\end{align}
which follows from linearising equations \eqref{eqn 2.11} \eqref{eqn 2.12}.

\item For the Gauss curvature $\bKl{u}$, there are also two ways to construct its linearised perturbation. The first strategy is to obtain the linearised perturbation of $\bKl{u}$ from the linearised perturbation of the metric $\bslashgl{u}$:
\begin{align}
\begin{aligned}
\bddsub{1}{\bKl{u=0}}
&=
\frac{1}{2r_0^4} \circdiv \circdiv \bdd{\bslashgl{u=0}}
-\frac{1}{2r_0^2} \circDelta \tr (\bdd{\bslashgl{u=0}})
- \frac{1}{2} K|_{\Sigma_{0,0}} \cdot \tr (\bdd{\bslashgl{u=0}})
\\
&=
-\frac{1}{r_0^3} \circDelta \bdd{\fl{u=0}} 
-\frac{2}{r_0^3} \bdd{\fl{u=0}}.
\end{aligned}
\label{eqn 8.5}
\end{align}
The second strategy to obtain the linearised perturbation of $\bKl{u}$ is to linearise the Gauss equation \eqref{eqn 2.13}:
\begin{align}
\begin{aligned}
\bddsub{2}{\bKl{u=0}}
&=
-\bdd{\brhol{u=0}} 
+ \frac{1}{4} \tr \chi'|_{\Sigma_{0,0}} \cdot \bdd{\btr \buchil{u=0}} 
+ \frac{1}{4} \tr \uchi|_{\Sigma_{0,0}} \cdot \bdd{\btr \bchil{u=0}'}
\\
&=
-\frac{1}{r_0^3} \circDelta \bdd{\fl{u=0}} -\frac{2}{r_0^3} \bdd{\fl{u=0}}.
\end{aligned}
\label{eqn 8.6}
\end{align}
Thus we see that two strategies give the same linearised perturbation of $\bKl{u=0}$.

\item $\bdd{\bmul{u=0}} = -\overline{\bdd{\brhol{u=0}}} = - \frac{3}{r_0^3} \overline{\bdd{\fl{u=0}}}$.

\item Linearised perturbations of curvature components at the initial leaf,
\begin{align}
\begin{aligned}
&
\begin{aligned}
\bdd{\bualphal{u=0}}
=
\bdd{\bal{u=0}^2 \cdot \ddualpha}|_{\Sigma_{0,0}}
=
\bdd{\ddualpha}|_{\Sigma_{0,0}}
=
0,
\end{aligned}
\\
&
\begin{aligned}
\bdd{\bubetal{u=0}}
=
\bdd{\bal{u=0} \cdot \ddubeta}|_{\Sigma_{0,0}}
=
\bdd{\ddubeta}|_{\Sigma_{0,0}}
=
\frac{3}{r_0^2} \slashd \bdd{\ufl{s=0}},
\end{aligned}
\\
&
\begin{aligned}
\bdd{\brhol{u=0}}
=
\bdd{\ddrho}|_{\Sigma_{0,0}}
=
\frac{3}{r_0^3} \bdd{\fl{u=0}},
\end{aligned}
\\
&
\begin{aligned}
\bdd{\bsigmal{u=0}}
=
\bdd{\ddsigma}|_{\Sigma_{0,0}}
=
0,
\end{aligned}
\\
&
\begin{aligned}
\bdd{\bbetal{u=0}}
=
\bdd{\bal{u=0}^{-1} \cdot \ddbeta}|_{\Sigma_{0,0}}
=
\bdd{\ddbeta}|_{\Sigma_{0,0}}
=
\frac{3}{r_0^2} \slashd \bdd{\fl{u=0}},
\end{aligned}
\\
&
\begin{aligned}
\bdd{\balphal{u=0}}
=
\bdd{\bal{u=0}^{-2} \cdot \ddalpha}|_{\Sigma_{0,0}}
=
\bdd{\ddalpha}|_{\Sigma_{0,0}}
=
0.
\end{aligned}
\end{aligned}
\label{eqn 8.7}
\end{align}
\end{enumerate}

\subsection{Explicit calculation of linearised perturbation at initial leaf with spherical harmonics}\label{subsec 8.2}

We already constructed linearised perturbation of the constant mass aspect function at the initial leaf in subsection \ref{subsec 8.1}. In this subsection, we explicitly calculate the linearised perturbation at the initial leaf with the help of spherical harmonics. 

Let $Y_l$ be a spherical harmonic of degree $l$, that $Y_l$ satisfies the eigenvalue equation
\begin{align*}
\circDelta Y_l = - \lambda_l Y_l,
\quad
\lambda_l = l(l+1),
\quad l\in \mathbb{N}.
\end{align*}
Note that $Y_{l=0}$ is simply a nonzero constant function. We shall use the convention that $Y_{l=0}\equiv 1$.

By the nature of the linear map, we decompose the derivation into two cases:
\begin{enumerate}[label=\emph{\roman*}.]
\item $\bdd{\fl{u=0}}=0$, $\bdd{\ufl{s=0}} = Y_l r_0$,
\item $\bdd{\fl{u=0}} = Y_l r_0$, $\bdd{\ufl{s=0}}=0$.
\end{enumerate}
Then the general case is simply the sum of two cases.

\subsubsection{Explicit calculation at initial leaf $\Sigma_{s=0,\us=0}$: case \emph{i}.}\label{subsubsec 8.2.1}

Substituting $\bdd{\fl{u=0}}=0$, $\bdd{\ufl{s=0}} = Y_l r_0$  to formulae in subsection \ref{subsec 8.1}, we obtain the explicit formulae of the linearised perturbation at the initial leaf $\Sigma_{s=0,\us=0}$ in case \emph{i.}. We consider two subcases: $l=0$ and $l \in \mathbb{Z}_{\geq 1}$.

In the subcase that $l=0$, the result is extremely simply that
\begin{align*}
\bdd{\ufl{u=0}} = \bdd{\ufl{s=0}} = Y_{l=0} r_0,
\quad
\bdd{\uh} = \bdd{\ufl{s=0}} = Y_{l=0} r_0,
\quad
\bdd{\bal{u=0}} = - Y_{l=0},
\end{align*}
and all the other linearised perturbations, including 
\begin{align*}
&
\bdd{\bslashgl{u=0}}, 
\quad 
\bdd{\br_{u=0}},
\\
&
\bdd{\buchil{u=0}}, 
\quad 
\bdd{\btr \buchil{u=0}}, 
\quad 
\bdd{\bchil{u=0}'}, 
\quad 
\bdd{\btr \bchil{u=0}'}, 
\quad 
\bdd{\btal{u=0}}, 
\quad 
\bddsub{2}{\buomegal{u=0}},
\\
&
\bdd{\bKl{u=0}},
\quad
\bdd{\bmul{u=0}},
\\
&
\bdd{\bualphal{u=0}}, 
\quad 
\bdd{\bubetal{u=0}}, 
\quad 
\bdd{\brhol{u=0}}, 
\quad 
\bdd{\bsigmal{u=0}}, 
\quad 
\bdd{\bbetal{u=0}}, 
\quad 
\bdd{\balphal{u=0}}
\end{align*}
all vanish. From equation \eqref{eqn 8.3}, we cannot determine $\bddsub{1}{\buomegal{u=0}}$ yet.

In the subcase $l \in \mathbb{Z}_{\geq 1}$, we have that
\begin{enumerate}[label=\ref{subsubsec 8.2.1}.\alph*.,leftmargin=.55in]
\item $\bdd{\ufl{u=0}} = \bdd{\ufl{s=0}} = Y_l r_0$, $\bdd{\uh} = \bdd{\ufl{s=0}} = Y_l r_0$.
\item $\bdd{\bal{u=0}} = 0$ by equation \eqref{eqn 8.1}.
\item $\bdd{\bslashgl{u=0}} = 0$, $\bdd{\br_{u=0}} =0$.
\item From equations \eqref{eqn 8.2},
\begin{align*}
&
\left\{
\begin{aligned}
\bdd{\buchil{u=0}}
&=
Y_l r_0 \circg 
- 2 r_0 \circnabla^2 Y_l,
\\
\bdd{\btr \buchil{u=0}}
&=
\bdd{\bal{u=0}\, \ddtr \dduchi}|_{\Sigma_{0,0}} 
=
\frac{2+2\lambda_l}{r_0} Y_l,
\\
\bdd{\bhatuchil{u=0}}
&=
-2 r_0 \widehat{\circnabla^2} Y_l,
\end{aligned}
\right.
\\
&
\left\{
\begin{aligned}
\bdd{\bchil{u=0}'}
&= 
0,
\\
\bdd{\btr \bchil{u=0}'}
&=
0,
\\
\bdd{\bhatchil{u=0}'}
&=
0,
\end{aligned}
\right.
\\
&\phantom{\Big\{}
\begin{aligned}
\bdd{\btal{u=0}}
&= 
0.
\end{aligned}
\end{align*}
By equation \eqref{eqn 8.4}, we have
\begin{align*}
\left\{
\begin{aligned}
\circDelta \bddsub{2}{\buomegal{u=0}}
&=
- \frac{3-3\lambda_l}{2r_0} Y_l,
\\
\overline{\bddsub{2}{\buomegal{u=0}}}
&=
0,
\end{aligned}
\right.
\end{align*}
thus
\begin{align*}
\bddsub{2}{\buomegal{u=0}} = \frac{3-3\lambda_l}{2\lambda_l r_0} Y_l.
\end{align*}

\item From equations \eqref{eqn 8.5} \eqref{eqn 8.6}, we obtain that $\bddsub{1}{\bKl{u=0}} = \bddsub{2}{\bKl{u=0}}=0$.

\item $\bdd{\bmul{u=0}} =0$.

\item By equations \eqref{eqn 8.7},
\begin{align*}
\bdd{\bubetal{u=0}} = \frac{3}{r_0} \slashd Y_l,
\end{align*}
and other linearised perturbations of curvature components $\bdd{\bualphal{u=0}}$, $\bdd{\brhol{u=0}}$, $\bdd{\bsigmal{u=0}}$, $\bdd{\bbetal{u=0}}$, $\bdd{\balphal{u=0}}$ all vanish.
\end{enumerate}

\subsubsection{Explicit calculation at initial leaf $\Sigma_{s=0,\us=0}$: case \emph{ii}.}\label{subsubsec 8.2.2}

Substituting $\bdd{\fl{u=0}} = Y_l r_0$, $\bdd{\ufl{s=0}} =0$ to formulae in subsection \ref{subsec 8.1}, we obtain the explicit formulae of the linearised perturbation at the initial leaf $\Sigma_{s=0, \us=0}$ in case \emph{ii.}. Similarly as in subsubsection \ref{subsubsec 8.2.1}, we consider two subcases: $l=0$ and $l\in \mathbb{Z}_{\geq 1}$.

In the subcase $l=0$, the result is that
\begin{align*}
&
\phantom{\Big\{}
\begin{aligned}
\bdd{\ufl{u=0}} = \bdd{\uh}= 0,
\end{aligned}
\\
&
\phantom{\Big\{}
\begin{aligned}
\bdd{\bal{u=0}} =0,
\end{aligned}
\\
&
\phantom{\Big\{}
\begin{aligned}
\bdd{\bslashgl{u=0}} = 2 r_0^2 Y_{l=0} \circg,
\quad
\bdd{\br_{u=0}} = Y_{l=0} r_0,
\end{aligned}
\\
&
\left\{
\begin{aligned}
&
\bdd{\buchil{u=0}} = Y_{l=0} r_0 \circg,
\quad
\bdd{\btr \buchil{u=0}} = -\frac{2}{r_0} Y_{l=0},
\\
&
\bdd{\bchil{u=0}'} = Y_{l=0} r_0 \circg,
\quad
\bdd{\btr \bchil{u=0}'} = \frac{2}{r_0} Y_{l=0},
\\
&
\bdd{\btal{u=0}} = 0,
\quad
\bddsub{2}{\buomegal{u=0}} = 0,
\end{aligned}
\right.
\\
&
\phantom{\Big\{}
\begin{aligned}
\bdd{\bKl{u=0}} = - \frac{2}{r_0^2} Y_{l=0},
\end{aligned}
\\
&
\phantom{\Big\{}
\begin{aligned}
\bdd{\bmul{u=0}} = - \frac{3}{r_0^2} Y_{l=0},
\end{aligned}
\\
&
\phantom{\Big\{}
\begin{aligned}
\bdd{\brhol{u=0}} = \frac{3}{r_0^2} Y_{l=0},
\end{aligned}
\end{align*}
and linearised perturbations of other curvature components vanish.

In the subcase $l\in \mathbb{Z}_{\geq 1}$, we have that
\begin{enumerate}[label=\ref{subsubsec 8.2.2}.\alph*., leftmargin=.55in]
\item $\bdd{\ufl{u=0}} = \bdd{\uh}=0$.
\item By equation \eqref{eqn 8.1}, we have that
\begin{align*}
\left\{
\begin{aligned}
\circDelta \bdd{\bal{u=0}}
&=
- (3+\lambda_l) Y_l ,
\\
\overline{\bdd{\bal{u=0}}}
&=
0,
\end{aligned}
\right.
\end{align*}
which implies that
\begin{align*}
\bdd{\bal{u=0}} = \frac{3+\lambda_l}{\lambda_l} Y_l.
\end{align*}

\item $\bdd{\bslashgl{u=0}} = 2 Y_l r_0^2 \circg$, $\bdd{\br_{u=0}} = 0$.
\item By equation \eqref{eqn 8.2},
\begin{align*}
&
\left\{
\begin{aligned}
\bdd{\buchil{u=0}}
&= 
\frac{3+2\lambda_l}{\lambda_l} Y_l r_0 \circg,
\\
\bdd{\btr \buchil{u=0}}
&=
\frac{6}{\lambda_l r_0}Y_l,
\\
\bdd{\bhatuchil{u=0}}
&= 
0,
\end{aligned}
\right.
\\
&
\left\{
\begin{aligned}
\bdd{\bchil{u=0}'}
&= 
Y_l r_0 \circg -2 r_0 \circnabla^2 Y_l,
\\
\bdd{\btr \bchil{u=0}'}
&=
\frac{2+2\lambda_l}{r_0} Y_l,
\\
\bdd{\bhatchil{u=0}'}
&= 
-2 r_0 \widehat{\circnabla^2} Y_l,
\end{aligned}
\right.
\\
&\phantom{\Big\{}
\begin{aligned}
\bdd{\btal{u=0}}
&=
\frac{3}{\lambda_l } \slashd Y_l.
\end{aligned}
\end{align*}
By equation \eqref{eqn 8.4},
\begin{align*}
\left\{
\begin{aligned}
\circDelta \bddsub{2}{\buomegal{u=0}}
&=
-\frac{9}{2 \lambda_l r_0} Y_l,
\\
\overline{\bddsub{2}{\buomegal{u=0}}}
&=
0,
\end{aligned}
\right.
\end{align*}
which implies that
\begin{align*}
\bddsub{2}{\buomegal{u=0}} = \frac{9}{2\lambda_l^2 r_0} Y_l.
\end{align*}

\item By equations \eqref{eqn 8.5} \eqref{eqn 8.6},
\begin{align*}
\begin{aligned}
\bddsub{1}{\bKl{u=0}} = \bddsub{2}{\bKl{u=0}}
=
\frac{\lambda_l -2}{r_0^2} Y_l.
\end{aligned}
\end{align*}

\item $\bdd{\bmul{u=0}} = 0$.
\item By equations \eqref{eqn 8.7},
\begin{align*}
&
\bdd{\brhol{u=0}} = \frac{3}{r_0^2} Y_l,
\\
&
\bdd{\bbetal{u=0}} = \frac{3}{r_0} \slashd Y_l,
\end{align*}
and other linearised perturbations of curvature components $\bdd{\bualphal{u=0}}$, $\bdd{\bubetal{u=0}}$, $\bdd{\bsigmal{u=0}}$, $\bdd{\balphal{u=0}}$ vanish.

\end{enumerate}

In both cases, equation \eqref{eqn 8.3} isnot sufficient to determine $\bddsub{1}{\buomegal{u=0}}$ yet. $\bddsub{1}{\buomegal{u=0}}$ will be determined in next section, see subsubsections \ref{subsubsec 9.2.1} and \ref{subsubsec 9.2.2}.

\section{Linearised perturbation of foliation: $\bddsub{1}{\cdot}$ by first strategy}\label{sec 9}
In this section, we shall linearise the inverse lapse problem of the constant mass aspect function foliation to construct the linearised perturbation of the foliation. This is the first strategy to construct linearised perturbations described in subsection \ref{subsec 7.2}. We shall use $\bddsub{1}{\cdot}$ to denote it.

\subsection{Construction of linearised perturbation $\bddsub{1}{\cdot}$ of foliation}\label{subsec 9.1}

As explained in subsections \ref{subsec 7.2} and \ref{subsec 8.1}, we assume that the most basic linearised perturbations $\bdd{\fl{u=0}}$ and $\bdd{\ufl{s=0}}$ are known. The other linearised perturbations are obtained by the following formulae.
\begin{enumerate}[label=\ref{subsec 9.1}.\alph*., leftmargin=.45in]
\item $\bddsub{1}{\ufl{u=0}} = \bddsub{1}{\uh} = \bdd{\ufl{s=0}}$. $\bddsub{1}{\fl{u}}$ satisfies the following equation
\begin{align}
\frac{\ed}{\ed u} \bddsub{1}{\fl{u}} = \bddsub{1}{\bal{u}},
\label{eqn 9.1}
\end{align}
from linearising equation \eqref{eqn 2.3}.

\item $\bddsub{1}{\bal{u}}$ satisfies the similar elliptic equation as \eqref{eqn 8.1}, obtained from linearising equation \eqref{eqn 2.2},
\begin{align}
\left\{
\begin{aligned}
\frac{1}{r^2} \circDelta \bddsub{1}{\bal{u}}
&=
- \bddsub{1}{\ddrho}|_{\Sigma_{s=u,0}} 
+ \overline{\bddsub{1}{\ddrho}|_{\Sigma_{s=u,0}}} 
- \frac{1}{r^2} \circdiv \bddsub{1}{\ddeta}|_{\Sigma_{s=u,0}},
\\
\overline{\bddsub{1}{\bal{u=0}}}
&=
- \frac{1}{r} \bddsub{1}{\br_{u}}
- \frac{r}{2} \overline{\bddsub{1}{\ddtr \dduchi}|_{\Sigma_{s,0}}},
\end{aligned}
\right.
\label{eqn 9.2}
\end{align}
which is equivalent to
\begin{align}
\left\{
\begin{aligned}
\frac{1}{r^2} \circDelta \bddsub{1}{\bal{u}}
&=
- \frac{3r_0}{r^4} \big( \bddsub{1}{\fl{u}} -\overline{\bddsub{1}{\fl{u}}} \big)
+ \frac{1}{r^3} \circDelta \bddsub{1}{\fl{u}}
\\
&\phantom{=}
- \frac{3r_0 u}{r^5} \big( \bdd{\ufl{s=0}} - \overline{\bdd{\ufl{s=0}}} \big)
- \frac{u}{r^4} \circDelta \bdd{\ufl{s=0}} ,
\\
\overline{\bddsub{1}{\bal{u}}}
&=
- \frac{1}{r} \bddsub{1}{\br_{u}}
+ \frac{1}{r} \overline{\bddsub{1}{\fl{u}}} + \frac{r-2r_0}{r^2} \overline{\bdd{\ufl{s=0}}},
\end{aligned}
\right.
\tag{\ref{eqn 9.2}}
\end{align}

\item $\bddsub{1}{\bslashgl{u}} = 2 r\bddsub{1}{\fl{u}} \circg + 2 u \bdd{\ufl{s=0}} \circg$, $\bddsub{1}{\br_u} = \overline{\bddsub{1}{\fl{u}}} + \frac{u}{r} \overline{\bdd{\ufl{s=0}}}$.

\item Linearised perturbations of connection coefficients are obtained by the following equations
\begin{align}
\begin{aligned}
&
\left\{
\begin{aligned}
\bddsub{1}{\buchil{u}}
&= 
\bddsub{1}{\bal{u}\, \dduchi}|_{\Sigma_{s=u,0}}
=
\bddsub{1}{\dduchi}|_{\Sigma_{s=u,0}} + \bddsub{1}{\bal{u}} \cdot \uchi|_{\Sigma_{s=u,0}},
\\
\bddsub{1}{\btr \buchil{u}}
&=
\bddsub{1}{\bal{u}\, \ddtr \dduchi}|_{\Sigma_{s=u,0}}
=
\bddsub{1}{\ddtr \dduchi}|_{\Sigma_{s=u,0}} + \bddsub{1}{\bal{u}} \tr \uchi|_{\Sigma_{s=u,0}},
\end{aligned}
\right.
\\
&
\left\{
\begin{aligned}
\bddsub{1}{\bchil{u}'}
&= 
\bddsub{1}{\bal{u}^{-1}\, \ddchi'}|_{\Sigma_{s=u,0}}
=
\bddsub{1}{\ddchi'}|_{\Sigma_{s=u,0}} - \bddsub{1}{\bal{u}} \cdot \chi'|_{\Sigma_{s=u,0}},
\\
\bddsub{1}{\btr \bchil{u}'}
&=
\bddsub{1}{\bal{u}^{-1}\, \ddtr \ddchi'}|_{\Sigma_{s=u,0}}
=
\bddsub{1}{\ddtr \ddchi'}|_{\Sigma_{s=u,0}} - \bddsub{1}{\bal{u}} \cdot \tr \chi'|_{\Sigma_{s=u,0}},
\end{aligned}
\right.
\\
&\phantom{\Big\{}
\begin{aligned}
\bddsub{1}{\btal{u}}
&= 
\slashd \bddsub{1}{\bal{u}}
+ \bddsub{1}{\ddeta}|_{\Sigma_{s=u,0}},
\end{aligned}
\\
&\phantom{\Big\{}
\begin{aligned}
\bddsub{1}{\buomegal{u}}
&= 
\frac{1}{2}\frac{\ed}{\ed u} \bddsub{1}{\bal{u}} + \bddsub{1}{\duomega}|_{\Sigma_{s=u,0}},
\end{aligned}
\end{aligned}
\label{eqn 9.3}
\end{align}
which are equivalent to
\begin{align}
\begin{aligned}
&
\left\{
\begin{aligned}
\bddsub{1}{\buchil{u}}
&= 
\bddsub{1}{\fl{u}} \circg 
+ \bdd{\ufl{s=0}} \circg 
- 2 \circnabla^2 \bdd{\ufl{s=0}} 
+ r \bddsub{1}{\bal{u}} \circg,
\\
\bddsub{1}{\btr \buchil{u}}
&=
-\frac{2}{r^2} \bddsub{1}{\fl{u}} 
- \frac{2r-4r_0}{r^3} \bdd{\ufl{s=0}} 
- \frac{2}{r^2} \circDelta \bdd{\ufl{s=0}}
+ \frac{2}{r} \bddsub{1}{\bal{u}},
\\
\bddsub{1}{\bhatuchil{u}}
&= 
- 2 \widehat{\circnabla^2} \bdd{\ufl{s=0}},
\end{aligned}
\right.
\\
&
\left\{
\begin{aligned}
\bddsub{1}{\bchil{u}'}
&= 
\bddsub{1}{\fl{u}} \circg 
+ \frac{u^2}{r^2} \bdd{\ufl{s=0}} \circg 
- 2 \circnabla^2 \bddsub{1}{\fl{u}} 
- u \bddsub{1}{\bal{u}} \circg,
\\
\bddsub{1}{\btr \bchil{u}'}
&=
\frac{2(r_0-u)}{r^3}\bddsub{1}{\fl{u}} 
- \frac{2u^2}{r^4} \bdd{\ufl{s=0}}
- 2 \circDelta \bddsub{1}{\fl{u}} 
- \frac{2u}{r^2} \bddsub{1}{\bal{u}},
\\
\bddsub{1}{\bhatchil{u}'}
&= 
- 2 \widehat{\circnabla^2} \bddsub{1}{\fl{u}} 
\end{aligned}
\right.
\\
&\phantom{\Big\{}
\begin{aligned}
\bddsub{1}{\btal{u}}
&= 
\slashd \bddsub{1}{\bal{u}}
+ \frac{u}{r^2} \slashd \bdd{\ufl{s=0}}
- \frac{1}{r} \slashd \bddsub{1}{\fl{u}},
\end{aligned}
\\
&\phantom{\Big\{}
\begin{aligned}
\bddsub{1}{\buomegal{u}}
&= 
\frac{1}{2}\frac{\ed}{\ed u} \bddsub{1}{\bal{u}} - \frac{r_0}{r^3} \bdd{\ufl{s=0}},
\end{aligned}
\end{aligned}
\tag{\ref{eqn 9.3}}
\end{align}

\item Linearised perturbation of the Gauss curvature is obtained by the following equation
\begin{align}
\begin{aligned}
\bddsub{1}{\bKl{u}}
&=
\frac{1}{2r^4} \circdiv \circdiv \bdd{\bslashgl{u}}
-\frac{1}{2r^2} \circDelta \tr (\bdd{\bslashgl{u}})
- \frac{1}{2} K|_{\Sigma_{s=u,0}} \cdot \tr (\bdd{\bslashgl{u}})
\\
&=
-\frac{1}{r^3} \circDelta \bddsub{1}{\fl{u}} 
- \frac{2}{r^3} \bddsub{1}{\fl{u}}
-\frac{u}{r^4} \circDelta \bdd{\ufl{s=0}}
- \frac{2u}{r^4} \bdd{\ufl{s=0}}.
\end{aligned}
\label{eqn 9.4}
\end{align}

\item $\bddsub{1}{\bmul{u}} = - \overline{\bddsub{1}{\brhol{u}}}$.

\item Linearised perturbations of curvature components are obtained by the following equations,
\begin{align}
\begin{aligned}
&
\begin{aligned}
\bddsub{1}{\bualphal{u}}
=
\bddsub{1}{\bal{u}^2 \cdot \ddualpha}|_{\Sigma_{s=u,0}}
=
\bddsub{1}{\ddualpha}|_{\Sigma_{s=u,0}}
=
0,
\end{aligned}
\\
&
\begin{aligned}
\bddsub{1}{\bubetal{u}}
=
\bddsub{1}{\bal{u} \cdot \ddubeta}|_{\Sigma_{s=u,0}}
=
\bddsub{1}{\ddubeta}|_{\Sigma_{s=u,0}}
=
\frac{3r_0}{r^3} \slashd \bdd{\ufl{s=0}},
\end{aligned}
\\
&
\begin{aligned}
\bddsub{1}{\brhol{u}}
=
\bddsub{1}{\ddrho}|_{\Sigma_{s=u,0}}
=
\frac{3r_0}{r^4} \bddsub{1}{\fl{u}}
+ \frac{3 r_0 u}{r^5} \bdd{\ufl{s=0}},
\end{aligned}
\\
&
\begin{aligned}
\bddsub{1}{\bsigmal{u}}
=
\bddsub{1}{\ddsigma}|_{\Sigma_{s=u,0}}
=
0,
\end{aligned}
\\
&
\begin{aligned}
\bddsub{1}{\bbetal{u}}
=
\bddsub{1}{\bal{u}^{-1} \cdot \ddbeta}|_{\Sigma_{s=u,0}}
=
\bddsub{1}{\ddbeta}|_{\Sigma_{s=u,0}}
=
\frac{3r_0}{r^3} \slashd \bddsub{1}{\fl{u}},
\end{aligned}
\\
&
\begin{aligned}
\bddsub{1}{\balphal{u}}
=
\bddsub{1}{\bal{u}^{-2} \cdot \ddalpha}|_{\Sigma_{s=u,0}}
=
\bddsub{1}{\ddalpha}|_{\Sigma_{s=u,0}}
=
0.
\end{aligned}
\end{aligned}
\label{eqn 9.5}
\end{align}

\end{enumerate}

\subsection{Explicit calculation of linearised perturbation $\bddsub{1}{\cdot}$ of foliation}\label{subsec 9.2}
We calculate explicitly the linearised perturbation $\bddsub{1}{\cdot}$ of the constant mass aspect function foliation at $\{\Sigma_{s,\us=0}\}$ with the help of spherical harmonics. Same as in subsection \ref{subsec 8.2}, we decompose the calculation into two cases:
\begin{enumerate}[label=\emph{\roman*}.]
\item $\bdd{\fl{u=0}}=0$, $\bdd{\ufl{s=0}} = Y_l r_0$,
\item $\bdd{\fl{u=0}} = Y_l r_0$, $\bdd{\ufl{s=0}}=0$.
\end{enumerate}

From equations in subsection \ref{subsec 9.1}, we can make the following ansatz on the linearised perturbation $\bddsub{1}{\cdot}$.
\begin{enumerate}[label=\ref{subsec 9.2}.\alph*., leftmargin=.45in]
\item $\bddsub{1}{\fl{u}} = \deltasub{1,\fl{u}}(u) Y_l$.

\item $\bddsub{1}{\bal{u}} = \deltasub{1,\bal{u}} (u) Y_l$.

\item $\bddsub{1}{\bslashgl{u}} = \deltasub{1,\subbslashgl{u}}(u) Y_l \circg$, $\bddsub{1}{\br_u} = \deltasub{1,\br_u}(u)$.

\item
\begin{align*}
&
\left\{
\begin{aligned}
\bddsub{1}{\bhatuchil{u}}
&= 
\deltasub{1,\bhatuchil{u}}(u) \widehat{\circnabla^2} Y_l,
\\
\bddsub{1}{\btr \buchil{u}}
&=
\deltasub{1,\btr \buchil{u}}(u) Y_l,
\end{aligned}
\right.
\\
&
\left\{
\begin{aligned}
\bddsub{1}{\bhatchil{u}'}
&= 
\deltasub{1,\bhatchil{u}'}(u) \widehat{\circnabla^2} Y_l,
\\
\bddsub{1}{\btr \bchil{u}'}
&=
\deltasub{1,\btr \bchil{u}'}(u) Y_l,
\end{aligned}
\right.
\\
&\phantom{\Big\{}
\begin{aligned}
\bddsub{1}{\btal{u}}
&= 
\deltasub{1,\btal{u}}(u) \slashd Y_l,
\end{aligned}
\\
&\phantom{\Big\{}
\begin{aligned}
\bddsub{1}{\buomegal{u}}
&= 
\deltasub{1,\buomegal{u}} (u) Y_l,
\end{aligned}
\end{align*}

\item $\bddsub{1}{\bKl{u}} = \deltasub{1,\bKl{u}}(u) Y_l$.

\item $\bddsub{1}{\bmul{u}} = \deltasub{1,\bmul{u}}(u)$.

\item
\begin{align*}
&
\begin{aligned}
\bddsub{1}{\bubetal{u}}
=
\deltasub{1,\bubetal{u}}(u) \slashd Y_l,
\end{aligned}
\\
&
\begin{aligned}
\bddsub{1}{\brhol{u}}
=
\deltasub{1,\brhol{u}} (u) Y_l,
\end{aligned}
\\
&
\begin{aligned}
\bddsub{1}{\bbetal{u}}
=
\deltasub{1,\bbetal{u}}(u) \dslashd Y_l.
\end{aligned}
\end{align*}

\end{enumerate}

\subsubsection{Explicit calculation of linearised perturbation $\bddsub{1}{\cdot}$: case \emph{i}.}\label{subsubsec 9.2.1}
Substituting the ansatz of the linearised perturbation $\bddsub{1}{\cdot}$ into equations in subsection \ref{subsec 9.1}, we obtain the following system of equations.
\begin{enumerate}[label=\ref{subsubsec 9.2.1}.\alph*., leftmargin=.55in]
\item By equation \eqref{eqn 9.1},
\begin{align*}
\frac{\ed}{\ed u} \deltasub{1,\fl{u}} = \deltasub{1,\bal{u}}.
\end{align*}

\item By equations \eqref{eqn 9.2},
\begin{align*}
l\geq 1:
&\quad
\begin{aligned}
- \frac{\lambda_l}{r^2} \deltasub{1,\bal{u}}
&=
-\frac{3r_0}{r^4} \deltasub{1,\fl{u}} - \frac{\lambda_l}{r^3} \deltasub{1,\fl{u}}
- \frac{3 r_0^2 u}{r^5} + \frac{\lambda_l r_0 u}{r^4},
\end{aligned}
\\
l =0:
&\quad
\begin{aligned}
\deltasub{1,\bal{u}}
&=
- \frac{r_0^2}{r^2} ,
\end{aligned}
\end{align*}

\item
$\deltasub{1,\subbslashgl{u}} = 2r \deltasub{1,\fl{u}} + 2 u r_0$,
$\deltasub{1,\br_u} = 
\left\{ 
\begin{matrix}
0, & l \geq 1, \\
\deltasub{1,\fl{u}} + \frac{u r_0}{r}, & l=0. 
\end{matrix} 
\right.$

\item By equations \eqref{eqn 9.3},
\begin{align*}
\begin{aligned}
&
\left\{
\begin{aligned}
\deltasub{1,\btr \buchil{u}}
&=
-\frac{2}{r^2} \deltasub{1,\fl{u}}
- \frac{2r r_0-4r_0^2}{r^3}
+ \frac{2r_0 \lambda_l}{r^2}
+ \frac{2}{r} \deltasub{1,\bal{u}},
\\
\bddsub{1}{\bhatuchil{u}}
&= 
\left\{
\begin{aligned}
&
- 2 r_0,
&&
l\geq 1,
\\
&
0,
&&
l=0,
\end{aligned}
\right.
\end{aligned}
\right.
\\
&
\left\{
\begin{aligned}
\deltasub{1,\btr \bchil{u}'}
&=
\frac{2(r_0-u)}{r^3}\deltasub{1,\fl{u}} 
- \frac{2r_0 u^2}{r^4}
+ \frac{2 \lambda_l}{r^2} \deltasub{1,\fl{u}} 
- \frac{2u}{r^2} \deltasub{1,\bal{u}},
\\
\bddsub{1}{\bhatchil{u}'}
&=
\left\{
\begin{aligned} 
&
- 2 \deltasub{1,\fl{u}},
&&
l\geq 1,
\\
&
0,
&&
l=0,
\end{aligned}
\right.
\end{aligned}
\right.
\\
&\phantom{\Big\{}
\begin{aligned}
\deltasub{1, \btal{u}}
&= 
\left\{
\begin{aligned}
&
\deltasub{1,\bal{u}}
+ \frac{u r_0}{r^2}
- \frac{1}{r} \deltasub{1,\fl{u}},
&&
l\geq 1,
\\
&
0,
&&
l=0,
\end{aligned}
\right.
\end{aligned}
\\
&\phantom{\Big\{}
\begin{aligned}
\deltasub{1,\buomegal{u}}
&= 
\frac{1}{2}\frac{\ed}{\ed u} \deltasub{1,\bal{u}} - \frac{r_0^2}{r^3},
\end{aligned}
\end{aligned}
\end{align*}

\item By equation \eqref{eqn 9.4},
\begin{align*}
\deltasub{1,\bKl{u}} 
=
\frac{\lambda_l-2}{r^3} \deltasub{1,\fl{u}} + \frac{u r_0}{r^4} (\lambda_l -2)
=
\frac{\lambda_l-2}{r^3} \deltasub{1,\fl{u}} + \frac{r_0}{r^3} (\lambda_l -2) - \frac{r_0^2}{r^4} (\lambda_l -2).
\end{align*}

\item $\deltasub{1,\bmul{u}} = 
\left\{ 
\begin{matrix}
0, & l \geq 1, \\
-\deltasub{1,\brhol{u}}, & l=0. 
\end{matrix} 
\right.$

\item By equations \eqref{eqn 9.5},
\begin{align*}
\deltasub{1,\bubetal{u}}
&=
\left\{
\begin{aligned}
&
\frac{3r_0^2}{r^3},&&
l\geq 1,
\\
&
0,
&&
l=0,
\end{aligned}
\right.
\\
\deltasub{1,\brhol{u}}
&=
\frac{3r_0}{r^4} \deltasub{1,\fl{u}} + \frac{3 r_0^2 u}{r^5},
\\
\deltasub{1,\bbetal{u}}
&=
\left\{
\begin{aligned}
&
\frac{3r_0}{r^3} \deltasub{1,\fl{u}}&&
l\geq 1,
\\
&
0,
&&
l=0.
\end{aligned}
\right.
\end{align*}

\end{enumerate}

We solve the above system of equations. From \ref{subsubsec 9.2.1}.a\&b, we derive that
\begin{align*}
l\geq 1:&
\quad
\frac{\ed}{\ed u} \deltasub{1,\fl{u}}
=
\Big(\frac{3r_0}{\lambda_l r^2} + \frac{1}{r} \Big) \deltasub{1,\fl{u}}
+\frac{3 r_0^2 u}{\lambda_l r^3} 
- \frac{r_0 u}{r^2},
\\
l=0:&
\quad
\frac{\ed}{\ed u} \deltasub{1,\fl{u}}
=
-\frac{r_0^2}{r^2}.
\end{align*}
Solving the above equation with the initial data $\deltasub{1,\fl{u}}(u=0) =0$, we obtain that
\begin{align*}
l\geq 1:&
\quad
\deltasub{1,\fl{u}}(u) 
= 
\frac{r}{3} \lambda_l(\lambda_l +2)
- r_0(\lambda_l +1) 
+ \frac{r_0^2}{r}
- \frac{r}{3} \lambda_l(\lambda_l-1) \exp \big[ \frac{3}{\lambda_l}(1-\frac{r_0}{r}) \big],
\\
l=0:&
\quad
\deltasub{1,\fl{u}}(u) = - r_0 (1- \frac{r_0}{r}).
\end{align*}
Substituting $\deltasub{1,\fl{u}}(u)$ to formulae of other linearised perturbations, we obtained that
\begin{enumerate}[leftmargin=.45in]
\item[$l\geq 1$:] $\deltasub{1,\fl{u}}(u) 
= 
\frac{r}{3} \lambda_l(\lambda_l +2)
- r_0(\lambda_l +1) 
+ \frac{r_0^2}{r}
- \frac{r}{3} \lambda_l(\lambda_l-1) \exp \big[ \frac{3}{\lambda_l}(1-\frac{r_0}{r}) \big]$,
\begin{align*}
\deltasub{1,\bal{u}}(u)
&=
\frac{\lambda_l(\lambda_l+2)}{3} 
-\frac{r_0^2}{r^2}
- (\lambda_l-1) \big(\frac{\lambda_l}{3} + \frac{r_0}{r} \big) \exp \big[ \frac{3}{\lambda_l}(1-\frac{r_0}{r}) \big],
\\
\deltasub{1,\subbslashgl{u}}(u)
&=
\frac{2r^2}{3} \lambda_l(\lambda_l +2)
- 2rr_0\lambda_l
- \frac{2r^2}{3} \lambda_l(\lambda_l-1) \exp \big[ \frac{3}{\lambda_l}(1-\frac{r_0}{r}) \big] ,
\\
\deltasub{1,\br_u}(u)
&=
0,
\\
\deltasub{1,\btr \buchil{u}}(u)
&=
\frac{4\lambda_l r_0}{r^2} 
+ \frac{2r_0}{r^2}(1-\lambda_l) \exp \big[ \frac{3}{\lambda_l}(1-\frac{r_0}{r}) \big],
\\
\deltasub{1,\bhatuchil{u}}
&=
-2r_0,
\\
\deltasub{1,\btr \bchil{u}'}(u)
&=
- \frac{2 \lambda_l r_0^2}{r^3} 
+ \frac{4\lambda_l r_0}{r^2} 
+ \frac{2\lambda_l(\lambda_l-2)(\lambda_l+2)}{r}
\\
&\phantom{=}
- \frac{\lambda_l-1}{r}\Big(\frac{2r_0^2}{r^2} +2(\lambda_l-1) \frac{r_0}{r} + \frac{2}{3}\lambda_l(\lambda_l-2)\Big) \exp \big[ \frac{3}{\lambda_l}(1-\frac{r_0}{r}) \big],
\\
\deltasub{1,\bhatchil{u}'}(u)
&=
-\frac{2r}{3} \lambda_l(\lambda_l +2) 
+2 r_0(\lambda_l +1) 
- \frac{2 r_0^2}{r}
+ \frac{2r}{3} \lambda_l(\lambda_l-1) \exp \big[ \frac{3}{\lambda_l}(1-\frac{r_0}{r}) \big],
\\
\deltasub{1,\btal{u}}(u)
&=
-\frac{3r_0^2}{r^2}
+\frac{r_0(\lambda_l+2)}{r}
-\frac{r_0}{r} (\lambda_l-1) \exp \big[ \frac{3}{\lambda_l}(1-\frac{r_0}{r}) \big],
\\
\deltasub{1,\buomegal{u}}(u)
&=
-\frac{3r_0^2(\lambda_l-1)}{2\lambda_l r^3} \exp \big[ \frac{3}{\lambda_l}(1-\frac{r_0}{r}) \big],
\\
\deltasub{1,\bKl{u}}(u)
&=
\frac{1}{3r^2} \lambda_l(\lambda_l +2)(\lambda_l-2)
- \frac{r_0}{r^3}\lambda_l (\lambda_l-2)
\\
&\phantom{=\;}
- \frac{1}{3r^2} \lambda_l(\lambda_l-1)(\lambda_l-2) \exp \big[ \frac{3}{\lambda_l}(1-\frac{r_0}{r}) \big].
\\
\deltasub{1,\bmul{u}}(u)
&=
0,
\\
\deltasub{1,\bubetal{u}}(u)
&=
\frac{3r_0^2}{r^3},
\\
\deltasub{1,\brhol{u}}(u)
&=
-\frac{3\lambda_l r_0^2}{r^4}
+\frac{r_0\lambda_l (\lambda_l+2)}{r^3}
-\frac{r_0}{r^3} \lambda_l  (\lambda_l-1) \exp \big[ \frac{3}{\lambda_l}(1-\frac{r_0}{r}) \big],
\\
\deltasub{1,\bbetal{u}}(u)
&=
\frac{r_0}{r^2} \lambda_l(\lambda_l +2)
- \frac{3r_0^2}{r^3}(\lambda_l +1) 
+ \frac{3r_0^3}{r^4}
- \frac{r_0}{r^2} \lambda_l(\lambda_l-1) \exp \big[ \frac{3}{\lambda_l}(1-\frac{r_0}{r}) \big].
\end{align*}

\item[$l =0$:] $\deltasub{1,\fl{u}}(u) 
= -r_0 \big( 1-\frac{r_0}{r} \big)$,
\begin{align*}
\deltasub{1,\bal{u}}(u)
&=
-\frac{r_0^2}{r^2},
\\
\deltasub{1,\subbslashgl{u}}(u)
&=
0,
\\
\deltasub{1,\br_u}(u)
&=
0,
\end{align*}
and other linearised perturbations $\deltasub{1,\btr \buchil{u}}(u)$, $\deltasub{1,\bhatuchil{u}}(u)$, $\deltasub{1,\btr \bchil{u}'}(u)$, $\deltasub{1,\bhatchil{u}'}(u)$, $\deltasub{1,\btal{u}}(u)$, $\deltasub{1,\buomegal{u}}(u)$, $\deltasub{1,\bmul{u}}$, $\deltasub{1,\bubetal{u}}(u)$, $\deltasub{1,\brhol{u}}(u)$, $\deltasub{1,\bbetal{u}}(u)$ all vanish.

\end{enumerate}

\subsubsection{Explicit calculation of linearised perturbation $\bddsub{1}{\cdot}$: case \emph{ii}.}\label{subsubsec 9.2.2}

Substituting the ansatz of the linearised perturbation $\bddsub{1}{\cdot}$ into equations in subsection \ref{subsec 9.1}, we obtain the following system of equations.
\begin{enumerate}[label=\ref{subsubsec 9.2.2}.\alph*., leftmargin=.55in]
\item By equation \eqref{eqn 9.1},
\begin{align*}
\frac{\ed}{\ed u} \deltasub{1,\fl{u}} = \deltasub{1,\bal{u}}.
\end{align*}

\item By equations \eqref{eqn 9.2},
\begin{align*}
l\geq 1:
&\quad
\begin{aligned}
- \frac{\lambda_l}{r^2} \deltasub{1,\bal{u}}
&=
-\frac{3r_0}{r^4} \deltasub{1,\fl{u}} - \frac{\lambda_l}{r^3} \deltasub{1,\fl{u}},
\end{aligned}
\\
l =0:
&\quad
\begin{aligned}
\deltasub{1,\bal{u}}
&=
0,
\end{aligned}
\end{align*}

\item
$\deltasub{1,\subbslashgl{u}} = 2r \deltasub{1,\fl{u}}$,
$\deltasub{1,\br_u} = 
\left\{ 
\begin{matrix}
0, & l \geq 1, \\
\deltasub{1,\fl{u}}, & l=0. 
\end{matrix} 
\right.$

\item By equations \eqref{eqn 9.3},
\begin{align*}
\begin{aligned}
&
\left\{
\begin{aligned}
\deltasub{1,\btr \buchil{u}}
&=
-\frac{2}{r^2} \deltasub{1,\fl{u}}
+ \frac{2}{r} \deltasub{1,\bal{u}},
\\
\bddsub{1}{\bhatuchil{u}}
&= 
0,
\end{aligned}
\right.
\\
&
\left\{
\begin{aligned}
\deltasub{1,\btr \bchil{u}'}
&=
\frac{2(r_0-u)}{r^3}\deltasub{1,\fl{u}} 
+ \frac{2 \lambda_l}{r^2} \deltasub{1,\fl{u}} 
- \frac{2u}{r^2} \deltasub{1,\bal{u}},
\\
\bddsub{1}{\bhatchil{u}'}
&=
\left\{
\begin{aligned} 
&
- 2 \deltasub{1,\fl{u}},
&&
l\geq 1,
\\
&
0,
&&
l=0,
\end{aligned}
\right.
\end{aligned}
\right.
\\
&\phantom{\Big\{}
\begin{aligned}
\deltasub{1, \btal{u}}
&= 
\left\{
\begin{aligned}
&
\deltasub{1,\bal{u}}
- \frac{1}{r} \deltasub{1,\fl{u}},
&&
l\geq 1,
\\
&
0,
&&
l=0,
\end{aligned}
\right.
\end{aligned}
\\
&\phantom{\Big\{}
\begin{aligned}
\deltasub{1,\buomegal{u}}
&= 
\frac{1}{2}\frac{\ed}{\ed u} \deltasub{1,\bal{u}},
\end{aligned}
\end{aligned}
\end{align*}

\item By equation \eqref{eqn 9.4},
\begin{align*}
\deltasub{1,\bKl{u}} 
=
\frac{\lambda_l-2}{r^3} \deltasub{1,\fl{u}}.
\end{align*}

\item $\deltasub{1,\bmul{u}} = 
\left\{ 
\begin{matrix}
0, & l \geq 1, \\
-\deltasub{1,\brhol{u}}, & l=0. 
\end{matrix} 
\right.$

\item By equations \eqref{eqn 9.5},
\begin{align*}
\deltasub{1,\bubetal{u}}
&=
0,
\\
\deltasub{1,\brhol{u}}
&=
\frac{3r_0}{r^4} \deltasub{1,\fl{u}},
\\
\deltasub{1,\bbetal{u}}
&=
\left\{
\begin{aligned}
&
\frac{3r_0}{r^3} \deltasub{1,\fl{u}}&&
l\geq 1,
\\
&
0,
&&
l=0.
\end{aligned}
\right.
\end{align*}

\end{enumerate}

We solve the above system of equations. From \ref{subsubsec 9.2.2}.a\&b, we derive that
\begin{align*}
l\geq 1:&
\quad
\frac{\ed}{\ed u} \deltasub{1,\fl{u}}
=
\Big(\frac{3r_0}{\lambda_l r^2} + \frac{1}{r} \Big) \deltasub{1,\fl{u}},
\\
l=0:&
\quad
\frac{\ed}{\ed u} \deltasub{1,\fl{u}}
=
0.
\end{align*}
Solving the above equation with the initial data $\deltasub{1,\fl{u}}(u=0) =r_0$, we obtain that
\begin{align*}
l\geq 1:&
\quad
\deltasub{1,\fl{u}}(u) 
= 
r \exp \big[ \frac{3}{\lambda_l}(1-\frac{r_0}{r}) \big],
\\
l=0:&
\quad
\deltasub{1,\fl{u}}(u) = r_0.
\end{align*}
Substituting $\deltasub{1,\fl{u}}(u)$ to formulae of other linearised perturbations, we obtained that
\begin{enumerate}[leftmargin=.45in]
\item[$l\geq 1$:] $\deltasub{1,\fl{u}}(u) 
= 
r \exp \big[ \frac{3}{\lambda_l}(1-\frac{r_0}{r}) \big]$,
\begin{align*}
\deltasub{1,\bal{u}}(u)
&=
\big( \frac{3r_0}{\lambda_l r} +1 \big) \exp \big[ \frac{3}{\lambda_l}(1-\frac{r_0}{r}) \big],
\\
\deltasub{1,\subbslashgl{u}}(u)
&=
2 r^2 \exp \big[ \frac{3}{\lambda_l}(1-\frac{r_0}{r}) \big] ,
\\
\deltasub{1,\br_u}(u)
&=
0,
\\
\deltasub{1,\btr \buchil{u}}(u)
&=
\frac{6r_0}{\lambda_l r^2}\exp \big[ \frac{3}{\lambda_l}(1-\frac{r_0}{r}) \big],
\\
\deltasub{1,\bhatuchil{u}}
&=
0,
\\
\deltasub{1,\btr \bchil{u}'}(u)
&=
\Big(\frac{6r_0^2}{\lambda_l r^3} +  \frac{6(\lambda_l-1) r_0}{\lambda_l r^2} + \frac{2\lambda_l -4}{r}\Big) \exp \big[ \frac{3}{\lambda_l}(1-\frac{r_0}{r}) \big],
\\
\deltasub{1,\bhatchil{u}'}(u)
&=
-2 r \exp \big[ \frac{3}{\lambda_l}(1-\frac{r_0}{r}) \big],
\\
\deltasub{1,\btal{u}}(u)
&=
\frac{3 r_0}{\lambda_l r} \exp \big[ \frac{3}{\lambda_l}(1-\frac{r_0}{r}) \big],
\\
\deltasub{1,\buomegal{u}}(u)
&=
\frac{9 r_0^2}{2\lambda_l^2 r^3} \exp \big[ \frac{3}{\lambda_l}(1-\frac{r_0}{r}) \big],
\\
\deltasub{1,\bKl{u}}(u)
&=
\frac{\lambda_l-2}{r^2} \exp \big[ \frac{3}{\lambda_l}(1-\frac{r_0}{r}) \big],
\\
\deltasub{1,\bmul{u}}(u)
&=
0,
\\
\deltasub{1,\bubetal{u}}(u)
&=
0,
\\
\deltasub{1,\brhol{u}}(u)
&=
\frac{3r_0}{r^3} \exp \big[ \frac{3}{\lambda_l}(1-\frac{r_0}{r}) \big],
\\
\deltasub{1,\bbetal{u}}(u)
&=
\frac{3r_0}{r^2} \exp \big[ \frac{3}{\lambda_l}(1-\frac{r_0}{r}) \big].
\end{align*}

\item[$l =0$:] $\deltasub{1,\fl{u}}(u) 
= r_0$,
\begin{align*}
\deltasub{1,\bal{u}}(u)
&=
0,
\\
\deltasub{1,\subbslashgl{u}}(u)
&=
2r r_0,
\\
\deltasub{1,\br_u}(u)
&=
r_0,
\\
\deltasub{1,\btr \buchil{u}}(u)
&=
-\frac{2 r_0}{r^2},
\\
\deltasub{1,\btr \bchil{u}'}(u)
&=
\frac{4r_0^2}{r^3}
- \frac{2 r_0}{r^2},
\\
\deltasub{1,\bKl{u}}(u)
&=
-\frac{2r_0}{r^3},
\\
\deltasub{1,\bmul{u}}(u)
&=
- \frac{3r_0^2}{r^4},
\\
\deltasub{1,\brhol{u}}(u)
&=
\frac{3r_0^2}{r^4},
\end{align*}
and other linearised perturbations $\deltasub{1,\bhatuchil{u}}(u)$, $\deltasub{1,\bhatchil{u}'}(u)$, $\deltasub{1,\btal{u}}(u)$, $\deltasub{1,\buomegal{u}}(u)$, $\deltasub{1,\bubetal{u}}(u)$, $\deltasub{1,\bbetal{u}}(u)$ all vanish.
\end{enumerate}

\section{Linearised perturbation of foliation: $\bddsub{2}{\cdot}$ by second strategy}\label{sec 10}
In this section, we shall linearise the basic equations for a constant mass aspect function foliation to obtain the linearised perturbation of the foliation, which is the second strategy described in subsection \ref{subsec 7.2}. We use $\bddsub{2}{\cdot}$ to denote it.

\subsection{Construction of linearised perturbation $\bddsub{2}{\cdot}$ of foliation}\label{subsec 10.1}
As in subsection \ref{subsec 9.1}, we assume that the most basic linearised perturbations $\bdd{\fl{u=0}}$ and $\bdd{\ufl{s=0}}$ are known. We linearised the basic equations of a constant mass aspect function foliation to derive equations of other linearised perturbations.
\begin{enumerate}[label=\ref{subsec 10.1}.\alph*., leftmargin=.5in]
\item $\bddsub{2}{\ufl{u=0}} = \bddsub{2}{\uh} = \bdd{\ufl{s=0}}$. $\bddsub{2}{\fl{u}}$ is obtained by solving the following elliptic equation
\begin{align*}
\left\{
\begin{aligned}
&
2 \circDelta \bddsub{2}{\fl{u}} 
-\frac{2(r_0-u)}{r^3} \big( \bddsub{2}{\fl{u}} -\overline{\bddsub{2}{\fl{u}}} \big)
\\
&\phantom{2 \circDelta } =
\big( \bddsub{2}{\btr \bchil{u}'} -\overline{\bddsub{2}{\btr \bchil{u}'}} \big)
- \frac{2u^2}{r^4} \big( \bdd{\ufl{s=0}} -\overline{\bddsub{2}{\ufl{s=0}}} \big)
- \frac{2u}{r^2} \big( \bddsub{2}{\bal{u}} -\overline{\bddsub{2}{\bal{u}}} \big),
\\
&
\frac{\ed}{\ed u} \overline{\bddsub{2}{\fl{u}}} = \overline{\bddsub{2}{\bal{u}}}.
\end{aligned}
\right.
\end{align*}
The first equation is similar to the equation of $\bddsub{1}{\btr \bchil{u}'}$ in \eqref{eqn 9.3}, and the second equation comes from linearising equation \eqref{eqn 2.3}.

\item $\bddsub{2}{\bal{u}}$ satisfies the following propagation equation
\begin{align*}
\frac{\ed}{\ed u} \bddsub{2}{\bal{u}}
=
2\bddsub{2}{\buomegal{u}} + \frac{2 r_0}{r^3} \bdd{\ufl{s=0}},
\end{align*}
which follows from linearising equation \eqref{eqn 7.1}.

\item $\bddsub{2}{\bslashgl{u}} = 2 r\bddsub{2}{\fl{u}} \circg + 2 u \bdd{\ufl{s=0}} \circg$, $\bddsub{1}{\br_u} = \overline{\bddsub{2}{\fl{u}}} + \frac{u}{r} \overline{\bdd{\ufl{s=0}}}$.

\item Linearised perturbations of connection coefficients satisfy the following equations
\begin{align*}
&
\begin{aligned}
\frac{\ed}{\ed u} \bddsub{2}{\btr \buchil{u}}
&=
2 \tr \uchi|_{\Sigma_{s=u,0}} \cdot \bddsub{2}{\buomegal{u}}
-
\tr \uchi|_{\Sigma_{s=u,0}} \cdot \bddsub{2}{\btr \buchil{u}},
\end{aligned}
\\
&
\begin{aligned}
\frac{\ed}{\ed u} \bddsub{2}{\btr \bchil{u}'}
&=
-2 \tr \chi'|_{\Sigma_{s=u,0}} \cdot \bddsub{2}{\buomegal{u}}
- \frac{1}{2} \tr \uchi|_{\Sigma_{s=u,0}} \cdot \bddsub{2}{\btr \bchil{u}'}
\\
&\phantom{=\;\;}
- \frac{1}{2} \tr \chi'|_{\Sigma_{s=u,0}} \cdot \bddsub{2}{\btr \buchil{u}}
+ 2 \bddsub{2}{\bmul{u}},
\end{aligned}
\\
&
\begin{aligned}
\frac{1}{r^2} \circdiv \bddsub{2}{\bhatuchil{u}} 
&=
\frac{1}{2} \slashd \bddsub{2}{\btr \buchil{u}} 
- \frac{1}{2} \tr \uchi|_{\Sigma_{s=u,0}} \cdot \bddsub{2}{\btal{u}} 
- \bddsub{2}{\bubetal{u}},
\end{aligned}
\\
&
\begin{aligned}
\frac{1}{r^2} \circdiv \bddsub{2}{\bhatchil{u}'} 
&=
\frac{1}{2} \slashd \bddsub{2}{\btr \bchil{u}'} 
+ \frac{1}{2} \tr \chi'|_{\Sigma_{s=u,0}} \cdot \bddsub{2}{\btal{u}} 
- \bddsub{2}{\bbetal{u}},
\end{aligned}
\\
&
\left\{
\begin{aligned}
&
\frac{1}{r^2} \circcurl \bddsub{2}{\btal{u}} = 0,
\\
&
\frac{1}{r^2} \circdiv \bddsub{2}{\btal{u}} = - \bddsub{2}{\brhol{u}} - \bddsub{2}{\bmul{u}},
\end{aligned}
\right.
\\
&
\left\{
\begin{aligned}
&
\frac{1}{r^2} \circDelta \bddsub{2}{\buomegal{u}}
=
- \frac{3}{4} \mu|_{\Sigma_{s=u,0}} \big( \bddsub{2}{\btr \buchil{u}} - \overline{ \bddsub{2}{\btr \buchil{u}} } \big) 
- \frac{1}{r^2} \circdiv \bddsub{2}{\bubetal{u}},
\\
&
\overline{\bddsub{2}{\buomegal{u}}}
=
0,
\end{aligned}
\right.
\end{align*}
which are equivalent to
\begin{align*}
&
\begin{aligned}
\frac{\ed}{\ed u} \bddsub{2}{\btr \buchil{u}}
&=
\frac{4}{r} \bddsub{2}{\buomegal{u}}
- \frac{2}{r} \bddsub{2}{\btr \buchil{u}},
\end{aligned}
\\
&
\begin{aligned}
\frac{\ed}{\ed u} \bddsub{2}{\btr \bchil{u}'}
&=
-\frac{2(r-r_0)}{r^2} \bddsub{2}{\buomegal{u}}
- \frac{1}{r} \bddsub{2}{\btr \bchil{u}'}
- \frac{r-r_0}{r^2} \bddsub{2}{\btr \buchil{u}}
+ 2 \bddsub{2}{\bmul{u}},
\end{aligned}
\\
&
\begin{aligned}
\circdiv \bddsub{2}{\bhatuchil{u}} 
&=
\frac{r^2}{2} \slashd \bddsub{2}{\btr \buchil{u}} 
- r \bddsub{2}{\btal{u}} 
- r^2 \bddsub{2}{\bubetal{u}},
\end{aligned}
\\
&
\begin{aligned}
\circdiv \bddsub{2}{\bhatchil{u}'} 
&=
\frac{r^2}{2} \slashd \bddsub{2}{\btr \bchil{u}'} 
+ (r-r_0) \bddsub{2}{\btal{u}} 
- r^2 \bddsub{2}{\bbetal{u}},
\end{aligned}
\\
&
\left\{
\begin{aligned}
&
\circcurl \bddsub{2}{\btal{u}} = 0,
\\
&
\circdiv \bddsub{2}{\btal{u}} = -r^2 \bddsub{2}{\brhol{u}} -r^2 \bddsub{2}{\bmul{u}},
\end{aligned}
\right.
\\
&\left\{
\begin{aligned}
&
\circDelta \bddsub{2}{\buomegal{u}}
=
- \frac{3r_0}{4 r} \mu|_{\Sigma_{s=u,0}} \big( \bddsub{2}{\btr \buchil{u}} - \overline{ \bddsub{2}{\btr \buchil{u}} } \big) - \circdiv \bddsub{2}{\bubetal{u}},
\\
&
\overline{\bddsub{2}{\buomegal{u}}}
=
0.
\end{aligned}
\right.
\end{align*}

\item Linearised perturbation of the Gauss curvature satisfies the following equation
\begin{align*}
\bddsub{2}{\bKl{u}}
&=
-\bddsub{2}{\brhol{u}} 
+ \frac{1}{4} \tr \chi'|_{\Sigma_{s=u,0}} \cdot \bddsub{2}{\btr \buchil{u}}
+ \frac{1}{4} \tr \uchi|_{\Sigma_{s=u,0}} \cdot \bddsub{2}{\btr \bchil{u}'}
\\
&=
-\bddsub{2}{\brhol{u}} 
+ \frac{r-r_0}{2r^2} \bddsub{2}{\btr \buchil{u}}
+ \frac{1}{2r} \bddsub{2}{\btr \bchil{u}'},
\end{align*}
by linearising the Gauss equation \eqref{eqn 2.13}.

\item Linearised perturbation of the mass aspect function satisfies the following equation
\begin{align*}
\frac{\ed}{\ed u} \bddsub{2}{\bmul{u}}
&=
-\frac{3}{2} \tr \uchi|_{\Sigma_{s=u,0}} \cdot \bddsub{2}{\bmul{u}}
- \frac{3}{2} \mu|_{\Sigma_{s=u,0}} \cdot \overline{\bddsub{2}{\btr \buchil{u}}},
\end{align*}
which is equivalent to
\begin{align*}
&
\begin{aligned}
\frac{\ed}{\ed u} \bddsub{2}{\bmul{u}}
&=
-\frac{3}{r}  \bddsub{2}{\bmul{u}}
- \frac{3r_0}{2r^3} \overline{\bddsub{2}{\btr \buchil{u}}}.
\end{aligned}
\end{align*}

\item Linearised perturbations of curvature components are given by the following equations
\begin{align*}
\begin{aligned}
&
\begin{aligned}
\bddsub{2}{\bualphal{u}}
=
\bddsub{2}{\bal{u}^2 \cdot \ddualpha}|_{\Sigma_{s=u,0}}
=
\bddsub{2}{\ddualpha}|_{\Sigma_{s=u,0}}
=
0,
\end{aligned}
\\
&
\begin{aligned}
\bddsub{2}{\bubetal{u}}
=
\bddsub{2}{\bal{u} \cdot \ddubeta}|_{\Sigma_{s=u,0}}
=
\bddsub{2}{\ddubeta}|_{\Sigma_{s=u,0}}
=
\frac{3r_0}{r^3} \slashd \bdd{\ufl{s=0}},
\end{aligned}
\\
&
\begin{aligned}
\bddsub{2}{\brhol{u}}
=
\bddsub{2}{\ddrho}|_{\Sigma_{s=u,0}}
=
\frac{3r_0}{r^4} \bddsub{2}{\fl{u}}
+ \frac{3 r_0 u}{r^5} \bdd{\ufl{s=0}},
\end{aligned}
\\
&
\begin{aligned}
\bddsub{2}{\bsigmal{u}}
=
\bddsub{2}{\ddsigma}|_{\Sigma_{s=u,0}}
=
0,
\end{aligned}
\\
&
\begin{aligned}
\bddsub{2}{\bbetal{u}}
=
\bddsub{2}{\bal{u}^{-1} \cdot \ddbeta}|_{\Sigma_{s=u,0}}
=
\bddsub{2}{\ddbeta}|_{\Sigma_{s=u,0}}
=
\frac{3r_0}{r^3} \slashd \bddsub{2}{\fl{u}},
\end{aligned}
\\
&
\begin{aligned}
\bddsub{2}{\balphal{u}}
=
\bddsub{2}{\bal{u}^{-2} \cdot \ddalpha}|_{\Sigma_{s=u,0}}
=
\bddsub{2}{\ddalpha}|_{\Sigma_{s=u,0}}
=
0.
\end{aligned}
\end{aligned}
\end{align*}

\end{enumerate}

\subsection{Explicit calculation of linearised perturbation $\bddsub{2}{\cdot}$ of foliation}\label{subsec 10.2}

We can adopt the same method in subsection \ref{subsec 9.2} to calculate the linearised perturbation $\bddsub{2}{\cdot}$. However, since we already obtain the linearised perturbation $\bddsub{1}{\cdot}$, and claim that two linearised perturbations $\bddsub{1}{\cdot}$ and $\bddsub{2}{\cdot}$ coincide, it is sufficient to prove the claim by checking that $\bddsub{1}{\cdot}$ satisfies equations of $\bddsub{2}{\cdot}$ in subsection \ref{subsec 10.1}. The verification is straightforward, thus we omit the details. We just mention the following formula used in the verification,
\begin{align*}
\circdiv \widehat{\circnabla^2} Y_l = \big( 1-\frac{\lambda_l}{2} \big)\slashd Y_l.
\end{align*}

Since $\bddsub{1}{\cdot}$ and $\bddsub{2}{\cdot}$ are the same, we use $\bdd{\cdot}$ to denote both for the sake of brevity.

\section{Linearised perturbation of asymptotic geometry of constant mass aspect function foliation}\label{sec 11}
We have constructed the linearised perturbation of constant mass aspect function foliation at $\{\Sigma_{s,\us=0}\}$. In this section, we shall apply it to construct linearised perturbation of the asymptotic geometry of the foliation at null infinity.

\subsection{Construction of linearised perturbation of asymptotic geometry}\label{subsec 11.1}
We assume that the most basic linearised perturbations $\bdd{\fl{u=0}}$ and $\bdd{\ufl{s=0}}$ are known. Recall the definitions of renormalised metric $\slashgl{u,r}$, renormalised Gauss curvature $\Kl{u,r}$ and their limits $\slashgl{\infty,r}$, $\Kl{\infty,r}$ in formulae \eqref{eqn 2.14}, \eqref{eqn 2.15}, \eqref{eqn 2.16}, then we construct their linearised perturbations by
\begin{align*}
\bdd{\bslashgl{u,r}}
&=
\frac{1}{r^2} \bdd{\bslashgl{u}} - \frac{2\bdd{\br_u}}{r} \circg,
\\
\bdd{\bKl{u,r}}
&=
r^2 \bdd{\bKl{u}} + 2 r \bdd{\br_u}\cdot K|_{\Sigma_{s=u,0}}
=
r^2 \bdd{\bKl{u}} + \frac{2}{r} \bdd{\br_u},
\end{align*}
and
\begin{align*}
\bdd{\bslashgl{\infty,r}}
=
\lim_{u \rightarrow + \infty} \bdd{\bslashgl{u,r}},
\quad
\bdd{\bKl{\infty,r}}
=
\lim_{u \rightarrow + \infty} \bdd{\bKl{u,r}}.
\end{align*}
Substituting $\bdd{\bslashgl{u}}$, $\bdd{\br_u}$ in \ref{subsec 9.1}.c. and $\bdd{\bKl{u}}$ in equation \eqref{eqn 9.4}, we obtain that
\begin{align*}
\bdd{\bslashgl{u,r}}
&=
\frac{2}{r} \big( \bdd{\fl{u}} -  \overline{\bdd{\fl{u}}} \big) \circg 
+ \frac{2 (r-r_0)}{r^2} \big( \bdd{\ufl{s=0}} - \overline{\bdd{\ufl{s=0}}} \big) \circg,
\\
\bdd{\bKl{u,r}}
&=
-\frac{1}{r} \circDelta \bdd{\fl{u}} 
- \frac{2}{r}  \big( \bdd{\fl{u}} -  \overline{\bdd{\fl{u}}} \big)
\\
&\phantom{=\;}
-\frac{r-r_0}{r^2} \circDelta \bdd{\ufl{s=0}}
- \frac{2(r-r_0)}{r^2}  \big( \bdd{\ufl{s=0}} - \overline{\bdd{\ufl{s=0}}} \big),
\end{align*}
and
\begin{align*}
\bdd{\bslashgl{\infty,r}}
&=
2 \lim_{u\rightarrow +\infty} \Big( \frac{\bdd{\ufl{u}}}{r} - \overline{\frac{\bdd{\ufl{u}}}{r}} \Big) \circg,
\\
\bdd{\bKl{\infty,r}}
&=
- \circDelta \lim_{u \rightarrow +\infty} \frac{\bdd{\fl{u}}}{r} 
- 2 \lim_{u \rightarrow +\infty}\Big( \frac{\bdd{\fl{u}}}{r} -  \overline{\frac{\bdd{\fl{u}}}{r}} \Big)
\end{align*}

\subsection{Explicit calculation of linearised perturbation of asymptotic geometry}\label{subsec 11.2}
We shall use the explicit calculations in subsection \ref{subsec 9.2} to calculate the linearised perturbations of asymptotic geometry at null infinity. Same as in subsection \ref{subsec 9.2}, we decompose the calculation into two cases:
\begin{enumerate}[label=\emph{\roman*}.]
\item $\bdd{\fl{u=0}}=0$, $\bdd{\ufl{s=0}} = Y_l r_0$,
\item $\bdd{\fl{u=0}} = Y_l r_0$, $\bdd{\ufl{s=0}}=0$.
\end{enumerate}

Following the ansatz in subsection \ref{subsec 9.2}, we introduce the following
\begin{align*}
&
\bdd{\bslashgl{u,r}} =
\deltasub{\subbslashgl{u,r}}(u) Y_l \circg,
\quad
\bdd{\bslashgl{\infty,r}} =
\deltasub{\subbslashgl{\infty,r}} Y_l \circg,
\\
&
\bdd{\bKl{u,r}} =
\deltasub{\bKl{u,r}}(u) Y_l,
\quad
\bdd{\bKl{\infty,r}} =
\deltasub{\bKl{\infty,r}} Y_l.
\end{align*}

\subsubsection{Explicit calculation of linearised perturbation: case \emph{i.}}\label{subsubsec 11.2.1}
Substituting $\bdd{\ufl{s=0}}$ and $\bdd{\fl{u}}$ in subsubsection \ref{subsubsec 9.2.1} into $\bdd{\bslashgl{u,r}}$, $\bdd{\bKl{u,r}}$, $\bdd{\bslashgl{\infty,r}}$, $\bdd{\bKl{\infty,r}}$, we obtain the following formulae of case \emph{i.}.
\begin{enumerate}[label=\ref{subsubsec 11.2.1}.\alph*., leftmargin=.60in]
\item 
\begin{align*}
\begin{aligned}
\deltasub{\subbslashgl{u,r}}(u)
&=
\left\{
\begin{aligned}
&
\frac{2}{3} \lambda_l(\lambda_l +2)
- \frac{2r_0}{r} \lambda_l 
- \frac{2}{3} \lambda_l(\lambda_l-1) \exp \big[ \frac{3}{\lambda_l}(1-\frac{r_0}{r}) \big],
\quad l\geq 1,
\\
&
0, \quad l=0,
\end{aligned}
\right.
\\
\deltasub{\subbslashgl{\infty,r}}
&=
\left\{
\begin{aligned}
&
\frac{2}{3} \lambda_l(\lambda_l +2)
- \frac{2}{3} \lambda_l(\lambda_l-1) \exp \big( \frac{3}{\lambda_l} \big),
\quad l\geq 1,
\\
&
0, \quad l=0,
\end{aligned}
\right.
\end{aligned}
\end{align*}

\item
\begin{align*}
\deltasub{\bKl{u,r}}(u)
&=
\left\{
\begin{aligned}
&
\frac{1}{3} \lambda_l(\lambda_l +2)(\lambda_l-2)
- \frac{r_0}{r} \lambda_l (\lambda_l-2)
- \frac{1}{3} \lambda_l(\lambda_l-1)(\lambda_l-2) \exp \big[ \frac{3}{\lambda_l}(1-\frac{r_0}{r}) \big],
\quad l\geq 1,
\\
&
0, \quad l=0,
\end{aligned}
\right.
\\
\deltasub{\bKl{\infty,r}}
&=
\left\{
\begin{aligned}
&
\frac{1}{3} \lambda_l(\lambda_l +2)(\lambda_l-2)
- \frac{1}{3} \lambda_l(\lambda_l-1)(\lambda_l-2) \exp \big( \frac{3}{\lambda_l} \big),
\quad l\geq 1,
\\
&
0, \quad l=0.
\end{aligned}
\right.
\end{align*}

\end{enumerate}

\subsubsection{Explicit calculation of linearised perturbation: case \emph{ii.}}\label{subsubsec 11.2.2}
Substituting $\bdd{\ufl{s=0}}$ and $\bdd{\fl{u}}$ in subsubsection \ref{subsubsec 9.2.2} into $\bdd{\bslashgl{u,r}}$, $\bdd{\bKl{u,r}}$, $\bdd{\bslashgl{\infty,r}}$, $\bdd{\bKl{\infty,r}}$, we obtain the following formulae of case \emph{ii.}.
\begin{enumerate}[label=\ref{subsubsec 11.2.2}.\alph*., leftmargin=.60in]
\item 
\begin{align*}
\begin{aligned}
\deltasub{\subbslashgl{u,r}}(u)
&=
\left\{
\begin{aligned}
&
2 \exp \big[ \frac{3}{\lambda_l}(1-\frac{r_0}{r}) \big],\\
&
0, \quad l=0.
\end{aligned}
\right.
\\
\deltasub{\subbslashgl{\infty,r}}
&=
\left\{
\begin{aligned}
&
2 \exp \big( \frac{3}{\lambda_l} \big),\\
&
0, \quad l=0.
\end{aligned}
\right.
\end{aligned}
\end{align*}

\item
\begin{align*}
\deltasub{\bKl{u,r}}(u)
&=
\left\{
\begin{aligned}
&
(\lambda_l-2) \exp \big[ \frac{3}{\lambda_l}(1-\frac{r_0}{r}) \big],\\
&
0, \quad l=0.
\end{aligned}
\right.
\\
\deltasub{\bKl{\infty,r}}
&=
\left\{
\begin{aligned}
&
(\lambda_l-2) \exp \big( \frac{3}{\lambda_l} \big),\\
&
0, \quad l=0.
\end{aligned}
\right.
\end{align*}

\end{enumerate}

\subsection{Geometric interpretation of linearised perturbation of asymptotic geometry in case \emph{i.}}\label{subsec 11.3}

The linearised perturbation in case \emph{i.} $\bdd{\fl{u=0}} =0, \bdd{\ufl{s=0}}=Y_l r_0$ has a particular geometric meaning, which we explain in this subsection.

\begin{figure}[H]
\begin{center}
\begin{tikzpicture}
\draw (-1,1.5) -- (-1,-1) 
(1,1.5) node[right] {$r=r_0$} -- (1,-1); 
\draw[thin] (-1,0) to [out=-90,in=180]
(0,-0.6) to [out=0,in=-90]
(1,0); 
\draw[dashed,thin] (1,0) to [out=90,in=0]
(0,0.6) to [out=180,in=90]
(-1,0); 
\draw[thin] (-1,0.1) -- (-4,-3.4)
(1,0.1) -- (4,-3.4); 
\draw[red,thick] (-1,0.6+0.2) to [out=-90,in=180]
(0,0.9+0.2) node[above] {\footnotesize $\bSigma_{u=0}$} to [out=0,in=-90]
(1,0.6+0.2); 
\draw[red,dashed,thick] (1,0.6+0.2) to [out=90,in=0]
(0,0.4+0.2) to [out=180,in=90]
(-1,0.6+0.2); 
\draw[red,thick] (-1,0.7+0.2) -- (-4,-2.8+0.2)
(1,0.7+0.2) -- (4,-2.8+0.2); 
\draw[red,thick]  (-3.55,-2.3+0.2) to [out=-125,in=180]
(0,-1.6+0.2) node[below] {\footnotesize $\bSigma_u$} to [out=0,in=-55]
(3.55,-2.3+0.2); 
\draw[red,dashed,thick] (3,-1.7+0.2) to [out=147,in=0]
(0,-2.6+0.2) to [out=180,in=33]
(-3,-1.7+0.2); 
\draw[->,blue] (0,-0.6) -- (0,1.1);
\node[right,blue] at (-0.1,0) {\tiny $\bdd{\ufl{s=0}}$};
\end{tikzpicture}
\end{center}
\caption{Linearised perturbation case \emph{i.}}
\label{fig 9}
\end{figure}
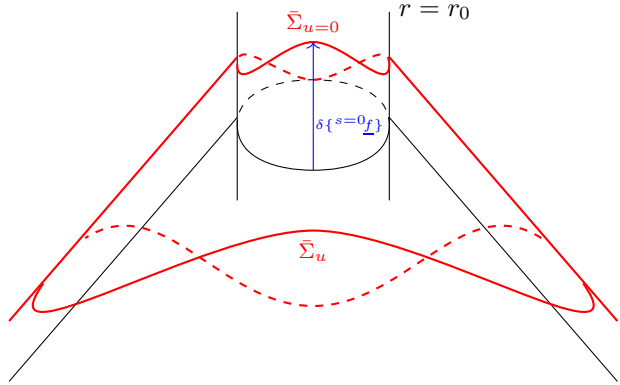

As illustrated in figure \ref{fig 9}, case \emph{i.} describes the linearised perturbation of the initial leaf $\bSigma_{u=0}$ within the event horizon $C_{s=0}$ where $r=r_0$. In other words, case \emph{i.} corresponds to the linearised perturbation through marginally trapped surfaces in the Schwarzschild spacetime. 

Such linearised perturbation through marginally trapped surfaces is of central importance when we concern the application to the null Penrose inequality. Recalling the project proposed by Christodoulou and Sauter, mentioned in subsection \ref{subsec 2.3}, it concerns the asymptotic geometry of the constant mass aspect function foliation starting from a marginally trapped surface, which is exactly the scenario covered by case \emph{i.}. We shall explain this in precise words in the following.

\subsubsection{Maps from marginally trapped surface to asymptotic geometry: $\bfg$ and $\bfk$}\label{subsec 11.3.1}
Any closed spacelike surface in the event horizon $C_{s=0}$ where $r=r_0$ is marginally trapped. Let $\calM$ denote the set of closed marginally trapped surfaces in the event horizon. Since any marginally trapped surface in $\calM$ can be parameterised by a function $\ufl{s=0}$ as its graph of $\us$ over $\vartheta$ domain in the coordinate system $\{\us, \vartheta\}$ in the event horizon $C_{s=0}$, as described in \ref{subsec 4.1}.b, we can parameterise the set $\calM$ by the function space on the sphere.

We introduce the following two maps from $\calM$ to the asymptotic geometry at past null infinity.
\begin{align}
\begin{aligned}
&
\mathbf{g}:
\quad \ufl{s=0}
\quad \stackrel{a}{\mapsto}
\quad \bSigma_{u=0} \in \calM
\quad \stackrel{b}{\mapsto}
\quad \{\bSigma_{u}\}
\quad \stackrel{c}{\mapsto}
\quad \bslashgl{\infty,r},
\\
&
\mathbf{k}:
\quad \ufl{s=0}
\quad \stackrel{a}{\mapsto}
\quad \bSigma_{u=0} \in \calM
\quad \stackrel{b}{\mapsto}
\quad \{\bSigma_{u}\}
\quad \stackrel{d}{\mapsto}
\quad \bKl{\infty,r}.
\end{aligned}
\label{eqn 11.1}
\end{align}

\begin{enumerate}[label=\alph*.]
\item $\bSigma_{u=0}$ is parameterised by $\ufl{s=0}$ as its graph of $\us$ over $\vartheta$ domain in the coordinate system $\{\us, \vartheta\}$ in the event horizon $C_{s=0}$,

\item $\{\bSigma_u\}$ is the constant mass aspect function foliation of $\ucalH$ starting from $\bSigma_{u=0}$,

\item if $\{ \bSigma_u \}$ extends to the past null infinity, then $\bslashgl{\infty,r}$ is the limit renormalised metric,

\item if $\{ \bSigma_u \}$ extends to the past null infinity, then $\bKl{\infty,r}$ is the limit renormalised Gauss curvature.

\end{enumerate}

In steps c.\&d., we need to assume that the global existence of the foliation $\{\bSigma_u\}$ extending to the past null infinity. Although the global existence assumption is not true in general, while for nearly spherically symmetric $\ucalH$, the assumption follows from the global existence results in \cite{S2008} and \cite{L2022}. Thus the maps $\bfg$ and $\bfk$ can be defined for smooth $\ufl{s=0}$ sufficiently small.\footnote{We donot specify the precise meaning of $\ufl{s=0}$ being small. An interpretation is the following: let $h$ be a smooth function and consider a one-parameter family of functions $\{ t h \}_{t\in \mathbb{R}}$, then $\bfg$ and $\bfk$ are defined for $th$ where $t$ is sufficiently small.}

\subsubsection{Linearised maps $\delta \bfg$ and $\delta \bfk$}
Denote the linearisations of $\bfg$, $\bfk$ at $\ufl{s=0}=0$ by $\delta \bfg$, $\delta \bfk$. Calculations in case \emph{i.} imply that the linearised maps $\delta \bfg$ and $\delta \bfk$ are the followings:
\begin{align}
\delta \bfg (Y_l r_0)
&=
\left\{
\begin{aligned}
&
\Big[ \frac{2}{3} \lambda_l(\lambda_l +2)
- \frac{2}{3} \lambda_l(\lambda_l-1) \exp \big( \frac{3}{\lambda_l} \big) \Big] Y_l \circg,
\quad l\geq 1,
\\
&
0, \quad l=0,
\end{aligned}
\right.
\label{eqn 11.2}
\\
\delta \bfk (Y_l r_0)
&=
\left\{
\begin{aligned}
&
\Big[ \frac{1}{3} \lambda_l(\lambda_l +2)(\lambda_l-2)
- \frac{1}{3} \lambda_l(\lambda_l-1)(\lambda_l-2) \exp \big( \frac{3}{\lambda_l} \big) \Big] Y_l,
\quad l\geq 1,
\\
&
0, \quad l=0.
\end{aligned}
\right.
\label{eqn 11.3}
\end{align}

Note that the image of $\delta \bfg$ is a pure multiplication of $\circg$, thus it implies that $\delta \bfg$ is the linear conformal deformation. It is easy to check that
\begin{align*}
\delta\bfk
=
\frac{1}{2} \circdiv \circdiv \delta \bfg 
- \frac{1}{2} \circDelta \tr (\delta \bfg) 
- \frac{1}{2} K_{\circg} \cdot \tr (\delta \bfg),
\end{align*}
which means that $\delta\bfk$ is the linearised perturbation of the Gauss curvature corresponding to the linear conformal deformation $\delta \bfg$.

We summarise the linearised perturbation of asymptotic geometry of the constant mass aspect function foliation which starts from a marginally trapped surface in the event horizon in the following theorem.

\begin{theorem}\label{thm 11.1}
Let $\bfg$ be the map characterising the limit renormalised metric, and $\bfk$ be the map characterise the limit renormalised Gauss curvature of the constant mass aspect function foliation starting from a marginally trapped surface in the event horizon, see formulae \eqref{eqn 11.1}. Let $\delta \bfg$ and $\delta \bfk$ be the corresponding linearisations at the marginally trapped surface $\Sigma_{0,0}$. The linearised map $\delta \bfg$ is a linear conformal deformation, and $\delta \bfk$ is the corresponding linearised perturbation of the Gauss curvature relative to $\delta \bfg$. The explicit formulae of $\delta \bfg$ and $\delta \bfk$ are given in \eqref{eqn 11.2} \eqref{eqn 11.3}.

Moreover $\delta \bfk$ is a bounded self-adjoint linear map from $H^2(\mathbb{S}^2,\circg)$ to $L^2(\mathbb{S}^2,\circg)$. Let $V$ be the linear space spanned by the spherical harmonics of degrees $0$ and $1$, and $V^{\perp}$ be the $L^2$ orthogonal complement of $V$. Then the kernel of $\delta \bfk$ is $V$ and the image of $\delta \bfk$ is $V^{\perp}$. $\delta \bfk$ is a bounded self-adjoint bijection from $H^2(\mathbb{S}^2,\circg) \cap V^{\perp}$ to $V^{\perp}$.

The four dimensional kernel $V$ of $\delta \bfk$ represents the set of the linearised perturbation of $\Sigma_{0,0}$, corresponding to which the linearised perturbation of the constant mass aspect function foliation stays being an asymptotic reference frame at null infinity.
\end{theorem}
\begin{proof}
The first paragraph of the theorem has been shown in the above. From formula \eqref{eqn 11.3} of $\delta \bfk$, we see that the spherical harmonic $Y_l r_0$ is the eigenfunction of $\delta \bfk$ corresponding to the eigenvalue $\Big[ \frac{1}{3} \lambda_l(\lambda_l +2)(\lambda_l-2)
- \frac{1}{3} \lambda_l(\lambda_l-1)(\lambda_l-2) \exp \big( \frac{3}{\lambda_l} \big) \Big] r_0^{-1}$ when $l\geq 1$ and the eigenvalue $0$ when $l=0$. This diagonalisation of $\delta \bfk$ implies that $\delta \bfk$ is self-adjoint.

In order to prove the second paragraph of the theorem, we examine the asymptotic of the eigenvalues of $\delta \bfk$ as $l \rightarrow +\infty$. Introduce the notation $k_l$ to denote the eigenvalues of $\delta \bfk$,
\begin{align*}
k_l 
=
\left\{
\begin{aligned}
&
\Big[ \frac{1}{3} \lambda_l(\lambda_l +2)(\lambda_l-2)
- \frac{1}{3} \lambda_l(\lambda_l-1)(\lambda_l-2) \exp \big( \frac{3}{\lambda_l} \big) \Big] r_0^{-1},
\quad
l \geq 1,
\\
&
0, \quad l=0.
\end{aligned}
\right. 
\end{align*}
Applying the Taylor expansion of $\exp \big( \frac{3}{\lambda_l} \big)$, we obtain that as $l \rightarrow +\infty$, $\lambda_l \rightarrow +\infty$
\begin{align*}
\frac{k_l}{\lambda_l} 
&=
\frac{\lambda_l^2}{3r_0} - \frac{4}{3r_0}  - \frac{1}{3r_0}(\lambda_l^2 -3\lambda_l +2)  \exp \big( \frac{3}{\lambda_l} \big)
\\
&=
-\frac{1}{2 r_0} + \frac{1}{ \lambda_l r_0} + \frac{3}{8 \lambda_l^2 r_0} + O\Big( \frac{1}{\lambda_l^3 r_0} \Big).
\end{align*}
Therefore 
\begin{align*}
\lim_{l\rightarrow +\infty} \frac{k_l}{\lambda_l}  = -\frac{1}{2r_0},
\end{align*}
which implies that $\delta \bfk$ is a bounded self-adjoint linear map from $H^2(\mathbb{S}^2,\circg)$ to $L^2(\mathbb{S}^2,\circg)$.

We show that $k_l=0$ if and only if $l=0,1$, which implies that $\ker (\delta \bfk) = V$ and $\Ima (\delta \bfk) = V^{\perp}$. The if part is obvious since $\lambda_{l=1} =2$. In order to prove the only if part, we show the sequence $\{ \frac{k_l}{\lambda_l} \}_{l\geq 2}$ is a negative sequence:
\begin{align*}
\frac{k_l}{\lambda_l}
& \leq
\frac{\lambda_l^2}{3r_0} - \frac{4}{3r_0}  - \frac{1}{3r_0}(\lambda_l^2 -3\lambda_l +2)  \big( 1+\frac{3}{\lambda_l} + \frac{9}{2\lambda_l^2} \big)
\\
& =
-\frac{(\lambda_l-2) (\lambda_l -3)}{2\lambda_l^2 r_0} 
\\
& <
0, \quad \lambda_l = l(l+1), \quad l=2,3,\cdots 
\end{align*}
Then we prove the only if part. The second paragraph of the theorem follows.

The third paragraph of the theorem follows from the second paragraph and the concept of an asymptotic reference frame introduced in subsection \ref{subsec 2.7}.
\end{proof}

\begin{remark}
We can carry out a similar analysis for the linearised perturbation of asymptotic geometry in case ii.. We briefly state the result in case ii. similar to theorem \ref{thm 11.1}. Define the linear maps $\delta \bfg_{ii.}$ and $\delta \bfk_{ii.}$ by
\begin{align*}
&
\delta \bfg_{ii.}:
\quad h 
\quad \mapsto 
\quad \bdd{\fl{u=0}} = h , \bdd{\ufl{u=0}}=0 
\quad \mapsto \bdd{\bslashgl{\infty,r}},
\\
&
\delta \bfk_{ii.}:
\quad h 
\quad \mapsto 
\quad \bdd{\fl{u=0}} = h , \bdd{\ufl{u=0}}=0 
\quad \mapsto \bdd{\bKl{\infty,r}}.
\end{align*}
The linearised perturbation $\delta \bfg_{ii.}$ of the renormalised metric at null infinity in case ii. is a linear conformal deformation and $\delta \bfk_{ii.}$ is the linearised perturbation of the Gauss curvature corresponding to $\delta \bfg_{ii.}$. The linearised perturbation $\delta \bfk_{ii.}$ of the renormalised Gauss curvature has the kernel $V$ and the image $V^{\perp}$. $\delta \bfk_{ii.}$ is a bounded self-adjoint bijection from $H^2(\mathbb{S}^2,\circg) \cap V^{\perp}$ to $V^{\perp}$. The following linearised perturbation of the initial leaf $\Sigma_{0,0}$
\begin{align*}
\bdd{\fl{u=0}} \in V, \quad \bdd{\ufl{u=0}} =0
\end{align*}
preserves the constant mass aspect function foliation to be an asymptotic reference frame at null infinity on the linearised level.

\end{remark}

\section{Linearised perturbation of energy-momentum vector at null infinity}\label{sec 12}

By theorem \ref{thm 11.1}, we already characterised the set of linearised perturbations of $\Sigma_{0,0}$, corresponding to which the perturbed constant mass aspect function foliation stays being an asymptotic reference frame at null infinity on the linearised level. This set of linearised perturbations of $\Sigma_{0,0}$ is the four dimensional kernel $V$ of the linearised map $\delta \bfk$. For such linearised perturbations of $\Sigma_{0,0}$, it is possible to define the corresponding linearised perturbation of the energy-momentum vector at null infinity. We discuss it in this section.

\subsection{Linearised perturbation of function $N$ at null infinity}\label{subsec 12.1}
In order to calculate the linearised perturbation of the energy-momentum vector at null infinity, we shall calculate the linearised perturbation of the function $N$ at null infinity, see definitions \ref{def 2.8} and \ref{def 2.9} of the energy-momentum vector. Recall that we use $\bdd{a}$ to denote the linearised perturbation of some quantity $a$. Note that for the constant mass aspect function foliation $\{\Sigma_{s,\us=0}\}$ of the null hypersurface $\uC_{\us=0}$, the limit function $\varSigma$, $\Xi$ both vanish, thus we have that for the linearised perturbation at this foliation,
\begin{align*}
\bdd{N} = - \bdd{P} = -\lim_{u\rightarrow \infty} \big( 3\bdd{\br_u} r_u^2\cdot \rhol{u} + r_u^3 \cdot \bdd{\brhol{u}} \big) .
\end{align*}
Directly applying the results in section \ref{sec 9}, we have that for the following two cases of the linearised perturbation of the initial leaf $\Sigma_{0,0}$:
\begin{enumerate}[label=\emph{\roman*}.]
\item $\bdd{\fl{u=0}}=0$, $\bdd{\ufl{s=0}} = Y_l r_0$,
\item $\bdd{\fl{u=0}} = Y_l r_0$, $\bdd{\ufl{s=0}}=0$,
\end{enumerate}
the corresponding linearised perturbation of the function $N$ at null infinity is given by
\begin{enumerate}[label=\emph{\roman*}.]
\item $\bdd{N} = 
\left\{
\begin{array}{l}
\big[- \lambda_l (\lambda_l+2)
+ \lambda_l  (\lambda_l-1) \exp \big( \frac{3}{\lambda_l} \big) \big] Y_l r_0, \quad l \geq 1, \\
0, \quad l=0. 
\end{array} 
\right.$
\item $\bdd{N} = 
\left\{
\begin{array}{l}
- 3 \exp \big( \frac{3}{\lambda_l} \big) Y_l r_0, \quad l \geq 1, \\
0, \quad l=0. 
\end{array} 
\right.$
\end{enumerate}

\subsection{Linearised perturbation of energy-momentum vector at null infinity}\label{subsec 12.2}

We calculate the linearised perturbation of the energy-momentum vector for the linearised perturbation of $\Sigma_{0,0}$ in $V$. 

\subsubsection{Linearised perturbation of Bondi energy at null infinity}\label{subsubsec 12.2.1}
We already calculate the linearised perturbation of the function $N$ above, thus we can easily calculate the linearised perturbation of the Bondi energy $E^{\gamma_{\infty}}$. For the following linearised perturbation of the initial leaf $\Sigma_{0,0}$,
\begin{align*}
\bdd{\fl{u=0}}=0,
\quad
\bdd{\ufl{s=0}} = Y_{l=1} r_0,
\end{align*}
we have that
\begin{enumerate}[label=\alph*.]
\item
$\bdd{N} = (2 e^{\frac{3}{2}} - 8) Y_{l=1} r_0$.

\item 
$\bdd{\bslashgl{\infty,r}} = \big[ \frac{16}{3} - \frac{4}{3} \exp \big( \frac{3}{2} \big) \big] Y_{l=1} \circg$, $\bdd{\dvol_{\subbslashgl{\infty,r}}} = \big[ \frac{16}{3} - \frac{4}{3} \exp \big( \frac{3}{2} \big) \big] Y_{l=1} \dvol_{\circg}$.
\item From equation \eqref{eqn 2.18},
\begin{align*}
\bdd{E^{\gamma_{\infty}}} = \frac{1}{8\pi} \Big[ \int \bdd{N} \cdot \dvol_{\circg} + \int N \cdot \bdd{\dvol_{\subbslashgl{\infty,r}}} \Big] =0.
\end{align*}
\end{enumerate}

\subsubsection{Linearised perturbation of linear momentum at null infinity}\label{subsubsec 12.2.2}
In order to calculate the linearised perturbation of the linear momentum, we shall first settle one more problem, explained in the following. Note that in definition \ref{def 2.9} equation \eqref{eqn 2.19}, the linear momentum is defined relative to a set of functions $\{x^1, x^2, x^3\}$ associated with an asymptotic reference frame $\gamma_{\infty}$ at null infinity. Since there is no canonical way to determine this set of functions for $\gamma_{\infty}$, we shall explain the choice of the set of functions when calculating the linearised perturbation of the linear momentum. We describe the choice of the set of functions $\{x^1, x^2, x^3\}$ and the linearised perturbation of these functions in the following. 

Firstly note that the set $\{x^1, x^2, x^3\}$ can be defined equivalently as an orthogonal basis $\{ Y_{l=1}^1, Y_{l=1}^2, Y_{l=1}^3\}$ of the first eigenspace of the Laplacian $\slashDelta_{\subbslashgl{\infty,r}}$ with the normalisation condition that
\begin{align*}
\int |Y_{l=1}^i|^2 \dvol_{\subbslashgl{\infty,r}} = \frac{4\pi}{3}.
\end{align*}
Then the linear momentum $\vec{P}^{\gamma_{\infty}}$ relative to the set of functions $\{ Y_{l=1}^1, Y_{l=1}^2, Y_{l=1}^3\}$ with respect to the reference frame $\gamma_{\infty}$ is defined by 
\begin{align*}
P^{\gamma_{\infty},i}(\ucalH) 
=
\frac{1}{8\pi} \int Y_{l=1}^i \cdot N \dvol_{\subslashgl{\infty,r}}.
\tag{\ref{eqn 2.19}}
\end{align*}
In order to calculate the linearised perturbation of the linear momentum, we shall describe the linearised perturbation $\bdd{Y_{l=1}^i}$. In the following, we give the equation satisfied by $\bdd{Y_{l=1}^i}$ for the linearised perturbation of the limit renormalised metric $\bdd{\bslashgl{\infty,r}} = \big[ \frac{16}{3} - \frac{4}{3} \exp \big( \frac{3}{2} \big) \big] Y_{l=1} \circg$.
\begin{enumerate}[label=\alph*.]
\item
$\bdd{\slashDelta_{\subbslashgl{\infty,r}}} = \big[ \frac{4}{3} \exp \big( \frac{3}{2} \big) - \frac{16}{3} \big] Y_{l=1} \circDelta$.

\item
Linearise the eigenvalue equation $\slashDelta_{\subbslashgl{\infty,r}} Y_{l=1}^i = - \lambda_{l=1} Y_{l=1}^i$, we obtain that
\begin{align}
\circDelta \bdd{Y_{l=1}^i} = -2 \bdd{Y_{l=1}^i} + 2 \big[ \frac{4}{3} \exp \big( \frac{3}{2} \big) - \frac{16}{3} \big] Y_{l=1} Y_{l=1}^i,
\label{eqn 12.1}
\end{align}
and as an corollary that
\begin{align*}
\int \Big\{  \bdd{Y_{l=1}^i} - \big[ \frac{4}{3} \exp \big( \frac{3}{2} \big) - \frac{16}{3} \big] Y_{l=1} Y_{l=1}^i \Big\} \dvol_{\circg} =0,
\tag{\ref{eqn 12.1}.co}
\label{eqn 12.1.co}
\end{align*}
which can be also derived easily from linearising the equation $\int Y_{l=1}^i \dvol_{\subbslashgl{\infty,r}} =0$.

\item
Linearise the normalisation condition $\int |Y_{l=1}^i|^2 \dvol_{\subbslashgl{\infty,r}} = \frac{4\pi}{3}$, we obtain that
\begin{align}
\int Y_{l=1}^i \bdd{Y_{l=1}^i} \dvol_{\circg} =0.
\label{eqn 12.2}
\end{align}

\item
Linearise the orthogonal condition $\int Y_{l=1}^i Y_{l=1}^j \dvol_{\subbslashgl{\infty,r}} = \frac{4\pi}{3} \delta_{ij}$, we obtain that
\begin{align}
\int Y_{l=1}^i \cdot \bdd{Y_{l=1}^j} \dvol_{\circg}  + \int Y_{l=1}^j \cdot \bdd{Y_{l=1}^i} \dvol_{\circg} =0.
\label{eqn 12.3}
\end{align}
\end{enumerate}
Then any set of solutions $\{ \bdd{Y_{l=1}^i} \}_{i=1,2,3}$ of equations \eqref{eqn 12.1} - \eqref{eqn 12.3} is an admissible linearised perturbation of the set of functions $\{ Y_{l=1}^i \}_{i=1,2,3}$ corresponding to the linearised perturbation $\bdd{\bslashgl{\infty,r}} =  \big[ \frac{16}{3} - \frac{4}{3} \exp \big( \frac{3}{2} \big) \big] Y_{l=1} \circg$. One can solve equations \eqref{eqn 12.1} - \eqref{eqn 12.3} through a geometric method via looking at the Lorentzian rotations in the Minkowski spacetime, but we donot demonstrate this solving procedure here. What we really need to calculate the linearised perturbation of the linear momentum is equation \eqref{eqn 12.1.co}.

Suppose that for the following linearised perturbation of the initial leaf $\Sigma_{0,0}$\begin{align*}
\bdd{\fl{u=0}}=0,
\quad
\bdd{\ufl{s=0}} = Y_{l=1} r_0,
\quad
Y_{l=1} = c_1 Y_{l=1}^1 + c_2 Y_{l=1}^2 + c_3 Y_{l=1}^3,
\end{align*}
we choose a corresponding admissible linearised perturbation $\{ \bdd{Y_{l=1}^i} \}_{i=1,2,3}$ of the set of functions $\{ Y_{l=1}^i \}_{i=1,2,3}$. We calculate the relative linearised perturbation $\bdd{P^{\gamma_{\infty},i}}$.
\begin{enumerate}[label=\alph*.]
\item
$\bdd{N} = (2 e^{\frac{3}{2}} - 8) Y_{l=1} r_0$.

\item 
$\bdd{\bslashgl{\infty,r}} = \big[ \frac{16}{3} - \frac{4}{3} \exp \big( \frac{3}{2} \big) \big] Y_{l=1} \circg$, $\bdd{\dvol_{\subbslashgl{\infty,r}}} = \big[ \frac{16}{3} - \frac{4}{3} \exp \big( \frac{3}{2} \big) \big] Y_{l=1} \dvol_{\circg}$.

\item From equation \eqref{eqn 2.19},
\begin{align*}
\bdd{P^{\gamma_{\infty},i}} 
&= 
\frac{1}{8\pi} \int N \cdot \big\{ \bdd{Y_{l=1}^i} \dvol_{\circg} + Y_{l=1}^i \bdd{\dvol_{\subbslashgl{\infty,r}}} \big\}
+
\frac{1}{8\pi} \int  Y_{l=1}^i \cdot \bdd{N}\dvol_{\circg}
\\
&=
\frac{1}{8\pi} \int N \cdot \Big\{ \bdd{Y_{l=1}^i} - \big[ \frac{4}{3} \exp \big( \frac{3}{2} \big) - \frac{16}{3} \big] Y_{l=1} Y_{l=1}^i \Big\} \dvol_{\circg}
+
\frac{1}{8\pi} \int  Y_{l=1}^i \cdot \bdd{N}\dvol_{\circg}
\\
&=
\frac{1}{8\pi} \int  (2 e^{\frac{3}{2}} - 8) Y_{l=1}  Y_{l=1}^i r_0 \dvol_{\circg}
\\
&=
c_i \big( \frac{1}{3} e^{\frac{3}{2}} - \frac{4}{3} \big) r_0,
\end{align*}
where the third equality follows from formula \eqref{eqn 12.1.co}.
\end{enumerate}
Note that the final result $\bdd{P^{\gamma_{\infty},i}}$ is independent of the choice of admissible linearised perturbation $\{ \bdd{Y_{l=1}^i} \}_{i=1,2,3}$.

\subsubsection{Linearised perturbation of Bondi mass at null infinity}\label{subsubsec 12.2.3}
We can easily calculate the linearised perturbation of the Bondi mass $\bdd{m_B}$ by $\bdd{E^{\gamma_{\infty}}}$ and $\bdd{P^{\gamma_{\infty},i}}$. For the following linearised perturbation of the initial leaf $\Sigma_{0,0}$
\begin{align*}
\bdd{\fl{u=0}}=0,
\quad
\bdd{\ufl{s=0}} = Y_{l=1} r_0,
\end{align*}
the corresponding linearised perturbation of the Bondi mass is
\begin{align*}
\bdd{m_{B}} = \bdd{\sqrt{(E^{\gamma_{\infty}})^2 - | \vec{P}^{\gamma_{\infty}} |^2}} = 0.
\end{align*}

\subsubsection{Summary of linearised perturbation of energy-momentum vector at null infinity}\label{subsubsec 12.2.4}

We summarise the above result of linearised perturbation of energy-momentum vector at null infinity in the following theorem.
\begin{theorem}\label{thm 12.1}
Consider an arbitrary linearised perturbation of the initial leaf $\Sigma_{0,0}$ in the kernel $V$ of the linear map $\delta \bfk$, which is assumed taking the form
\begin{align*}
\bdd{\fl{u=0}}=0,
\quad
\bdd{\ufl{s=0}} = c_0 + Y_{l=1} r_0,
\quad
Y_{l=1} = c_1 Y_{l=1}^1 + c_2 Y_{l=1}^2 + c_3 Y_{l=1}^3.
\end{align*}
By theorem \ref{thm 11.1}, the resulting linearised perturbation of the constant mass aspect function foliation is preserved to be an asymptotic reference frame at null infinity. The resulting linearised perturbation of the energy-momentum vector and the Bondi mass at null infinity are given by the following formulae
\begin{align*}
\bdd{E^{\gamma_{\infty}}} = 0,
\quad
\bdd{P^{\gamma_{\infty},i}} 
=
c_i \big( \frac{1}{3} e^{\frac{3}{2}} - \frac{4}{3} \big) r_0,
\quad
\bdd{m_B} = 0.
\end{align*}
The linear map from $V$ to the $3$-dimensional space of the linearised perturbation of the momentum vector
\begin{align*}
\bdd{\vec{P}^{\gamma_{\infty}}}
=
\big(\bdd{P^{\gamma_{\infty},1}}, \bdd{P^{\gamma_{\infty},2}}, \bdd{P^{\gamma_{\infty},3}}    \big)
\end{align*}
is surjective, and the kernel is $V_0$ which corresponds to the linearised perturbation of $\Sigma_{0,0}$ that $\bdd{\fl{u=0}}=0$, $\bdd{\ufl{s=0}} = c_0$.
\end{theorem}
\begin{proof}
We already prove the theorem in the case $\bdd{\ufl{s=0}} = Y_{l=1} r_0$, thus the case left is when $\bdd{\ufl{s=0}} = c_0$. This follows easily from the time translation invariance of the Schwarzschild spacetime.
\end{proof}

\begin{remark}
We can carry out the calculation in this section similarly for the linearised perturbation of the initial leaf $\Sigma_{0,0}$
\begin{align*}
\bdd{\fl{u=0}} = c_0 + Y_{l=1} r_0 \in V, 
\quad
\bdd{\ufl{u=0}} =0,
\quad
Y_{l=1} = c_1 Y_{l=1}^1 + c_2 Y_{l=1}^2 + c_3 Y_{l=1}^3.
\end{align*}
We briefly state the results of the corresponding linearised perturbations of the energy-momentum vector and the Bondi mass without derivations:
\begin{align*}
\bdd{E^{\gamma_{\infty}}}=0,
\quad
\bdd{P^{\gamma_{\infty},i}} = -\frac{c_i}{2} \exp \big( \frac{3}{2} \big) r_0,
\quad
\bdd{m_B} =0.
\end{align*}
\end{remark}

\section{Outlook on linearised perturbation of constant mass aspect function foliation in perturbed Schwarzschild spacetime}\label{sec 13}

In previous sections, we obtain the precise result on the linearised perturbation of the constant mass aspect function foliation in a Schwarzschild spacetime. As mentioned in the introduction section \ref{sec 1} and subsection \ref{subsec 2.3}, this result serves as the first step to carry out the project of proving the null Penrose inequality in a perturbed Schwarzschild spacetime using the constant mass aspect function foliation and its perturbation.

Aside from this paper, references \cite{L2020} \cite{L2022} addressed different aspects of a perturbed Schwarzschild spacetime which are also essential for the above project: \cite{L2020} dealt with the set of closed marginally trapped surfaces in a perturbed Schwarzschild spacetime, and \cite{L2022} concerned the global existence and regularity of null hypersurfaces in a perturbed Schwarzschild exterior.

The follow-up step should be to study the same linearised perturbation problem in a perturbed Schwarzschild spacetime. In this section, we give an outlook on this problem. We will briefly discuss the appropriate class of perturbed Schwarzschild spacetimes, formulate the linearised perturbation problem for the above class, and discuss the possible solution of the problem. We shall see that the precise result obtained in a Schwarzschild spacetime serves as the model, which contains the key structure, for the result in the general case of a perturbed Schwarzschild spacetime.

\subsection{Vacuum perturbed Schwarzschild metric}\label{subsec 13.1}
We shall describe an appropriate class of vacuum perturbed Schwarzschild metrics in this subsection. Here we adopt the method using a double null coordinate system as in \cite{L2022}. Let $\{ s, \us, \theta^1, \theta^2\}$ be the double null coordinate system introduced in subsection \ref{subsec 3.1}. We recall a definition introduced in \cite{L2022}.
\begin{definition}[$\kappa$-neighbourhood $M_{\kappa}$]\label{def 13.1}
Let $\{\us,s\}$ be the double null foliation of the Schwarzschild spacetime $\left(\calS, g \right)$ (see subsection \ref{subsec 3.1}), then the $\kappa$-neighbourhood $M_{\kappa}$ of the null hypersurface $\uC_{\us=0}$ is defined by
\begin{align*}
M_{\kappa} = \left\{ p\in \mathcal{S}: s(p)>-\kappa r_0, \vert \us \vert < \kappa r_0 \right\}.
\end{align*}
See figure \ref{fig 10}.
\begin{figure}[H]
\centering
\begin{tikzpicture}
\draw[->] (0,0) -- (3,-3) node[right] {$s$};
\draw[->] (0,0) -- (0.5,0.5) node[right] {$\us$};
\draw[dotted,thick] (3.25, -2.75) -- (0,0.5) -- (-0.5,0) -- (2.75,-3.25);
\draw[fill] (0,0)  circle [radius=0.03];
\draw[->] (-0.5,-0.5) node[left] {\footnotesize $\Sigma_{0,0}$} to [out=30,in=-120] (-0.06,-0.06);
\draw[->] (2,-1) node[right] {\footnotesize $\uC_{\us=0}$} to [out=-150,in=60]   (1.55,-1.45);
\draw[->] (1.4,-2.6) node[left] {\footnotesize$M_{\kappa}$: $\kappa$-neighbourhood }to [out=30,in=-110] (1.85,-2.15);
\end{tikzpicture}
\caption{$\kappa$-neighbourhood $M_{\kappa}$}
\label{fig 10}
\end{figure}
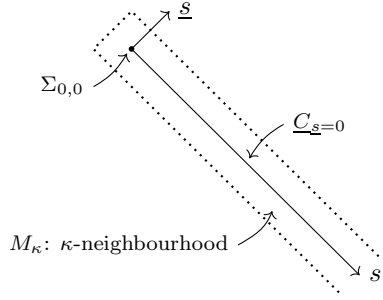
\end{definition}

Then we shall consider the perturbation of the Schwarzschild metric in the $\kappa$-neighbourhood $M_{\kappa}$. Use the double null coordinate system $\{ s, \us, \theta^1, \theta^2\}$ on $M_{\kappa}$ inherited from the Schwarzschild spacetime $(\calS, g)$. Up to coordinate transformations of $\{\theta^1, \theta^2\}$ on each $\Sigma_{\us,s}$, a general Lorentzian metric $g_{\kappa}$ on $M_{\kappa}$, with $\{ s, \us\}$ preserved being double null, takes the form
\begin{align*}
g_{\kappa}=2\Omega_{\kappa}^2  ( \ed s  \otimes \ed \us  +\ed \us \otimes \ed s ) + (\slashg_{\kappa})_{ab} \big( \ed \theta^a - b_{\kappa}^{a} \ed s\big) \otimes \big( \ed \theta^b - b_{\kappa}^{b} \ed s\big).
\end{align*}
In the following, assume that $g_{\kappa}$ is Ricci flat, i.e. $(M_{\kappa}, g_{\kappa})$ is vacuum. Following the construction in subsection \ref{subsec 3.2}, we use the metric components of $g_{\kappa}$ together with the connection coefficients and curvature components to describe the geometry of $(M_{\kappa}, g_{\kappa})$. Then as in \cite{L2022}, we describe the perturbed Schwarzschild metric by comparing these geometric quantities with their values in the Schwarzschild spacetime. We roughly explain the class of vacuum perturbed Schwarzschild metrics in the following definition.
\begin{definition}[Rough definition of the appropriate vacuum perturbed Schwarzschild metric]\label{def 13.2}
Let $\epsilon$ be a positive real number and $g_{\epsilon}$ be a Ricci-flat Lorentzian metric on $M_{\kappa}$ that in coordinates $\{\us,s,\theta^1,\theta^2\}$
\begin{align}
g_{\epsilon}=2\Omega_{\epsilon}^2  \left( \ed s  \otimes \ed \us  +\ed \us \otimes \ed s \right) + \left( \slashg_{\epsilon}\right)_{ab} \left( \ed \theta^a - b_{\epsilon}^{a} \ed s\right) \otimes \left( \ed \theta^b - b_{\epsilon}^{b} \ed s\right).
\end{align}
We define the area radius function $r_{\epsilon}(\us,s)$ by
\begin{align}
4\pi r_{\epsilon}^2(\us,s) =\int_{\Sigma_{\us,s}} 1 \cdot \dvol_{\subslashg_{\epsilon}}.
\end{align}
$g_{\epsilon}$ is called $\epsilon$-close to the Schwarzschild metric $g_{\calS}$ on $M_{\kappa}$, if certain assumptions on the difference between the metric components, the connection coefficients, the curvature components of $g_{\epsilon}$ and $g_{\calS}$ hold. We illustrate the form of assumptions by demonstrating for some components and coefficients in the following. Let $n$ be a positive integer.
\begin{enumerate}[label=\alph*.]
\item Assumptions on the metric components: demonstration for $\slashg_{\epsilon}$. For $k\geq 0$, $l\geq 1$, $m\geq 1$, $k +l +m \leq n+2$,
\begin{align*}
&
\begin{aligned}
&  1-\epsilon < \vert \frac{r_{\epsilon}}{r_{\calS}} \vert < 1+ \epsilon, 
\\
&  \vert \slashg_{\epsilon} - \slashg_{\calS} \vert_{\subcircg} < \epsilon r_{\epsilon}^2, 
&&  \vert \circnabla^{k} \left( \slashg_{\epsilon} - \slashg_{\calS} \right) \vert_{\subcircg} < \epsilon r_{\epsilon}^2, 
\\
&  \vert \circnabla^{k} \partial_{s}^l \left( \slashg_{\epsilon} - \slashg_{\calS} \right) \vert_{\subcircg} < \frac{\epsilon r_0 }{r_{\epsilon}^{l-1}},
&& \vert \circnabla^k \partial_{\us}^m \left( \slashg_{\epsilon} - \slashg_{\calS} \right) \vert_{\subcircg} < \frac{\epsilon r_{\epsilon}}{r_0^{m-1}},
\end{aligned}
\\
& \vert \circnabla^k \partial_{s}^l \partial_{\us}^m \left( \slashg_{\epsilon} - \slashg_{\calS} \right) \vert_{\subcircg} < \frac{\epsilon}{ r_0^{m-1} r_{\epsilon}^{l-1}}.
\end{align*}
\item Assumptions on the connection coefficients : demonstration for $\tr \chi_{\epsilon}$. For $k\geq 0$, $l\geq 1$, $m\geq 1$, $k+l +m \leq n+1$,
\begin{align*}
&
\begin{aligned}
& \vert \tr\chi_{\epsilon} - \tr\chi_{\calS} \vert_{\subcircg} < \frac{\epsilon}{ r_{\epsilon} },
&& \vert \circnabla^{k}  \left(  \tr\chi_{\epsilon} - \tr\chi_{\calS} \right) \vert_{\subcircg} < \frac{\epsilon }{ r_{\epsilon} } , 
\\
& \vert  \circnabla^k \partial_{s}^l  \left(  \tr\chi_{\epsilon} - \tr\chi_{\calS} \right) \vert_{\subcircg} < \frac{\epsilon }{r_{\epsilon}^{1+l}},
&& \vert  \circnabla^k \partial_{\us}^m  \left(  \tr\chi_{\epsilon} - \tr\chi_{\calS} \right)  \vert_{\subcircg} < \frac{\epsilon }{ r_{\epsilon}^{2}r_0^{m-1}},
\end{aligned}
\\
& \vert \circnabla^k \partial_{s}^l \partial_{\us}^m   \left(  \tr\chi_{\epsilon} - \tr\chi_{\calS} \right) \vert_{\circg} < \frac{\epsilon }{r_{\epsilon}^{2+l}r_0^{m-1}}.
\end{align*}
\item Assumptions on the curvature components: demonstration for $\ualpha_{\epsilon}$. For $k\geq 0$, $l\geq 1$, $m\geq 2$, $k+ l+ m \leq n$,
\begin{align*}
&
\begin{aligned}
& \vert \ualpha_{\epsilon} \vert_{\subcircg} < \frac{\epsilon r_0^{\frac{3}{2}}}{r_{\epsilon}^{\frac{3}{2}}},
&& \vert \circnabla^k \ualpha_{\epsilon} \vert_{\subcircg} < \frac{\epsilon r_0^{\frac{3}{2}}}{r_{\epsilon}^{\frac{3}{2}}},
&& \vert \circnabla^k \partial_{s}^l \ualpha_{\epsilon}  \vert_{\subcircg} < \frac{\epsilon r_0^{\frac{3}{2}}}{r_{\epsilon}^{\frac{3}{2}+l}},
\\
& \vert  \circnabla^k \partial_{\us}  \ualpha_{\epsilon}  \vert_{\subcircg} < \frac{\epsilon r_0^{\frac{3}{2}}}{r_{\epsilon}^{\frac{5}{2}}},
&& \vert \circnabla^k \partial_{s}^l \partial_{\us} \ualpha_{\epsilon}   \vert_{\subcircg} < \frac{\epsilon r_0^{\frac{3}{2}}}{r_{\epsilon}^{\frac{5}{2}+l}},
\end{aligned}
\\
&
\begin{aligned}
& \vert  \circnabla^k \partial_{\us}^{m}  \ualpha_{\epsilon}  \vert_{\subcircg} < \frac{\epsilon r_0}{r_{\epsilon}^{3} r_0^{m-2}},
&& \vert \circnabla^k \partial_{s}^l \partial_{\us}^{m} \ualpha_{\epsilon}   \vert_{\subcircg} < \frac{\epsilon r_0}{r_{\epsilon}^{3+l} r_0^{m-2}}.
\end{aligned}
\end{align*}
\end{enumerate}
\end{definition} 

\begin{remark}
Note that the quantity $\epsilon$ in above definition is dimensionless. We see that the above perturbation of the Schwarzschild metric is only local near a null hypersurface $\uC_{\us=0}$, rather than a global perturbation of the whole Schwarzschild exterior.

The decay assumptions in the above definition of $g_{\epsilon}$ are taken from the results in the global nonlinear stability of Minkowski spacetime \cite{CK1993} by Christodoulou and Klainerman. See also the extension in \cite{Bi2009}. Note that the Kerr black hole spacetime is strongly asymptotically flat, thus it satisfies the decay behaviour proved in \cite{CK1993} in a $\kappa$-neighbourhood $M_{\kappa}$. Therefore the Kerr black hole with sufficiently small angular momentum lies in the class of vacuum perturbed Schwarzschild spacetime in above definition. An explicit axisymmetric double null coordinate system of the Kerr black hole is given in \cite{PI1998}.
\end{remark}

\subsection{Formulation of linearised perturbation problem}\label{subsec 13.2}
\begin{figure}[H]
\begin{center}
\begin{tikzpicture}
\draw[dashed] (-1,0) to [out=70,in=-110] (-0.7,0.9);
\draw (-1,0) to [out=-110,in=70] (-1.85+0.5,-2.4+1.4);
\draw[dashed] (1,0) to [out=110,in=-70] (0.7,0.9);
\draw[->] (1,0) to [out=-70,in=110] (1.85-0.5,-2.4+1.4) node[right] {\small $s$}; 
\node[above right] at (1.3,-1.6) {\tiny $\uC_{\us=0}$};
\draw[dashed] (-1,0) to [out=-45,in=135] (-0.5,-0.5);
\draw (-1,0) to [out=135, in= -45] (-2,1);
\draw[dashed] (1,0) to [out=-135,in=45] (0.5,-0.5);
\draw[->] (1,0) to [out=45, in= -135] (1.9,0.9) node[right] {\tiny $C_{s=0}$}to [out=45,in=-135] (2.3,1.3) node[right] {\small $\us$}; 
\draw[dashed] (-1.5,0.5) to [out=135,in=45] (1.5,0.5);
\draw (1.5,0.5) to [out=-135,in=0] (0,0) to [out=180,in=-45] (-1.5,0.5); 
\node[below] at (0,0.1) {\scriptsize $\Sigma_0$}; 
\draw[dashed] (-1.4,1.3) to [out=-110,in=70] (-1.6,0.7);
\draw (-1.6,0.7) to [out=-110,in=70] (-2.75+0.5,-2.4+1.4);
\draw[dashed] (1.4,1.3) to [out=-70,in=110] (1.6,0.7);
\draw[->] (1.6,0.7) to [out=-70,in=110] (2.75-0.5,-2.4+1.4) node[right]{\small $s$}; 
\node[right] at (2.4-0.4,-1.5+1) {\scriptsize $\ucalH$}; 
\draw[dashed] (-2.3+0.4,-2.2+2.1) to [out=70,in=180] (-0.5,-2.5+2) node[below] {\scriptsize $\Sigma$} to [out=0,in=110] (2.3-0.4,-2.2+2.1);
\draw (2.3-0.4,-2.2+2.1) to [out=-70,in=0] (1,-2+1.9) to [out=180,in=-110] (-2.3+0.4,-2.2+2.1); 
\draw[->] (-1.73,0.73) to [out=-110,in=70] (-2,0);
\node[left] at (-2,0) {$\bfs(\ufl{s=0})$};
\draw[->] (-1.05,-0.1) to [out=135,in=-45] (-1.65,0.5);
\node[right] at (-1.6,0.5) {\tiny $\ufl{s=0}$};
\end{tikzpicture}
\end{center}
\caption{The parameterisation map $\bfs$.}
\label{fig 11}
\end{figure}
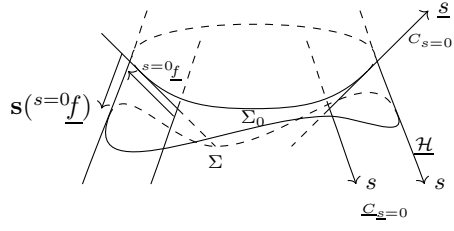

Let $(M_{\kappa}, g_{\epsilon})$ be as in definition \ref{def 13.2}. \cite{L2020} gave the characterisation of the set of closed marginally trapped surfaces in $(M_{\kappa}, g_{\epsilon})$. We present an overview of the characterisation here. Adopting the second kind of parametrisation $(f, \ufl{s=0})$ of a spacelike surface $\Sigma$ in \ref{subsec 4.1}.b of subsection \ref{subsec 4.1}, there exists a map $\bfs$ from the parameterisation function $\ufl{s=0}$ to $f$, such that the spacelike surface $\Sigma$ with the second parameterisation $(f, \ufl{s=0}) = (\bfs(\ufl{s=0}), \ufl{s=0})$ is a closed marginally trapped surface. This map $\bfs$ is called the parameterisation map of closed marginally trapped surfaces. See figure \ref{fig 11}.

Let $\bSigma_0$ be a closed marginally trapped surface and $\ucalH$ be the incoming null hypersurface where $\bSigma_0$ is embedded. By the method in subsection \ref{subsec 2.4}, we construct the constant mass aspect function foliation $\{ \bSigma_u \}$ emanating from $\bSigma_0$. 

The above construction naturally gives rise to a map from the set of closed marginally trapped surfaces to constant mass aspect function foliations in $(M_{\kappa}, g_{\epsilon})$. We can investigate, through this map, the change of the asymptotic geometry of the constant mass aspect function foliation $\{ \bSigma_u\}$ with respect to the deformation of the closed marginally trapped surface $\bSigma_0$. 

Now we can heuristically formulate the linearised perturbation problem of the asymptotic geometry of the constant mass aspect function foliation. Use the notations $\bfg_{\epsilon}$ and $\bfk_{\epsilon}$ to denote the maps in a vacuum perturbed Schwarzschild spacetime $(M_{\kappa}, g_{\epsilon})$ similar to $\bfg$, $\bfk$ in subsection \ref{subsec 11.3}. Let $\calM$ be the set of closed marginally trapped surfaces in $(M_{\kappa}, g_{\epsilon})$. The maps $\bfg_{\epsilon}$ and $\bfk_{\epsilon}$ are defined heuristically as follows.
\begin{align}
\begin{aligned}
&
\bfg_{\epsilon}:
\quad \ufl{s=0}
\quad \stackrel{a}{\mapsto} 
\quad (\bfs(\ufl{s=0}), \ufl{s=0})
\quad \stackrel{b}{\mapsto}
\quad \bSigma_0 \in \calM
\quad \stackrel{c}{\mapsto}
\quad \{\bSigma_{u}\}
\quad \stackrel{d}{\mapsto}
\quad \bslashgl{\infty,r},
\\
&
\bfk_{\epsilon}:
\quad \ufl{s=0}
\quad \stackrel{a}{\mapsto} 
\quad (\bfs(\ufl{s=0}), \ufl{s=0})
\quad \stackrel{b}{\mapsto}
\quad \bSigma_0 \in \calM
\quad \stackrel{c}{\mapsto}
\quad \{\bSigma_{u}\}
\quad \stackrel{e}{\mapsto}
\quad \bKl{\infty,r}.
\end{aligned}
\label{eqn 13.3}
\end{align}
\begin{enumerate}[label=\alph*.]
\item As explained in the beginning of this subsection, $\bfs$ is the parameterisation map of closed marginally trapped surfaces,

\item $\bSigma_0$ has the second parameterisation $(\bfs(\ufl{s=0}), \ufl{s=0})$, thus it is a closed marginally trapped surface in $(M_{\kappa}, g_{\epsilon})$,

\item $\{\bSigma_u\}$ is the constant mass aspect function foliation of $\ucalH$ emanating from $\bSigma_0$,

\item if $\{ \bSigma_u \}$ extends to the past null infinity, then $\bslashgl{\infty,r}$ is the limit renormalised metric,

\item if $\{ \bSigma_u \}$ extends to the past null infinity, then $\bKl{\infty,r}$ is the limit renormalised Gauss curvature.
\end{enumerate}
We see that the maps $\bfg_{\epsilon}$ and $\bfk_{\epsilon}$ characterise the asymptotic intrinsic geometry of leaves of the constant mass aspect function foliation emanating from a marginally trapped surface. If one is interested in the asymptotic behaviour of other geometric quantities, one can construct other corresponding maps by the above approach. Clearly in order to define the above maps rigorously, we have to solve several problems first: the global existence and regularity of the incoming null hypersurface $\ucalH$ from $\bSigma_0$ to the past null infinity, the global existence and regularity of the constant mass aspect function foliation $\{ \bSigma_u\}$ from $\bSigma_0$ to the past null infinity. These problems are already separately addressed in \cite{L2022} and \cite{S2008}.

Given the above constructions, we can formulate the problem of linearised perturbation of the asymptotic geometry of the constant mass aspect function foliation as the study of the linearised perturbations of maps $\bfg_{\epsilon}$, $\bfk_{\epsilon}$ and other maps corresponding to different geometric quantities associated with the foliation.

\subsection{Outlook on solution of linearised perturbation problem}\label{subsec 13.3}
In formulae \eqref{eqn 11.2}, \eqref{eqn 11.3}, we obtain the precise result of the linearised perturbations $\delta \bfg$ and $\delta \bfk$. Based on these formulae, we prove the main theorem \ref{thm 11.1} of this paper, summarising the property of the linearised perturbation of the asymptotic geometry of the constant mass aspect function foliation in a Schwarzschild spacetime.

Naturally the above precise result in a Schwarzschild spacetime serves as a model for the solution of the linearised perturbation problem in a vacuum perturbed Schwarzschild spacetime. It is reasonable to compare the corresponding linearisation perturbations $\delta \bfg_{\epsilon}$ with $\delta \bfg$, and $\delta \bfk_{\epsilon}$ with $\delta \bfk$. More precisely, let $\ufl{s=0}$ be the parameterisation function of a marginally trapped surface $\bSigma_0$ in $(M_{\kappa}, g_{\epsilon})$, we shall compare the linearised perturbations $\delta \bfg_{\epsilon}|_{\ufl{s=0}}, \delta \bfk_{\epsilon}|_{\ufl{s=0}}$ with $\delta \bfg, \delta \bfk$.

Note that $\delta \bfg, \delta \bfk$ in a Schwarzschild spacetime are the linearised perturbations at $\ufl{s=0}=0$ rather than at an arbitrary $\ufl{s=0}$. Thus there are at least two dimensionless quantities entering the comparison of the corresponding linearised perturbations: one is the dimensionless quantity $\epsilon$ in definition \ref{def 13.2} of $(M_{\kappa}, g_{\epsilon})$, and the other one, denoted by $\udelta$, measures the deviation of $\ufl{s=0}$ from a constant function.

A reasonable choice of the quantity $\udelta$ is a certain Sobolev norm of the differential of $\ufl{s=0}$. For example, let $n$ be a positive integer and $p > 1$, we assume that the $\rmW^{n,p}(\circg)$ Sobolev norm of $\slashd \ufl{s=0}$ is bounded by
\begin{align*}
\Vert \slashd \ufl{s=0} \Vert^{n,p} \leq \udelta r_0.
\end{align*}

Let $\Vvert \cdot \Vvert^{n_1 \rightarrow n_2,p}$ denote the operator norm from Sobolev spaces $\rmW^{n_1,p}(\circg)$ to $\rmW^{n_2,p}(\circg)$. The possible form taken by the comparison between $\delta \bfk_{\epsilon}|_{\ufl{s=0}}$ and $\delta \bfk$ could be
\begin{align*}
\Vvert \delta \bfk_{\epsilon}|_{\ufl{s=0}} - \delta \bfk \Vvert^{k \rightarrow k-2,p} \leq \varrho(k, p, \epsilon, \udelta) r_0^{-1}.
\end{align*}
where $\varrho(k, p, \epsilon, \udelta)$ is a certain function of $k,p$ and the small parameters $\epsilon, \udelta$. At present, we cannot determine the exact form of the comparison between $\delta \bfk_{\epsilon}|_{\ufl{s=0}}$ and $\delta \bfk$ without careful analyses in $(M_{\kappa}, g_{\epsilon})$. The rigorous and thoughtful analyses of this problem are left for subsequent papers. We end the discussion here by the last comment that we expect that the optimal choice of $k$ above should be $n+1$ in the regularity of $\ufl{s=0}$, and $\varrho(k, p,\epsilon, \udelta)$ could be a linear function of $\epsilon$ and $\udelta$, i.e. $\varrho(k, p, \epsilon, \udelta) = c(k,p) (\epsilon + \udelta)$ where $c(k,p)$ is a constant depending on $k,p$.

\section*{Acknowledgements}
\addcontentsline{toc}{section}{Acknowledgement}
This paper emerges from the author's thesis \cite{L2018} on the null Penrose inequality in a perturbed Schwarzschild black hole. The author is grateful to Demetrios Christodoulou for his constant encouragement and generous guidance. The author also thanks Alessandro Carlotto and Lydia Bieri for many helpful discussions.

\appendix
\section{Derivations of equations \eqref{eqn 2.5'} and \eqref{eqn 2.11'}}\label{appen A}

Note that equation \eqref{eqn 2.5'} follows from equations \eqref{eqn 2.11'} and the following lemma.
\begin{lemma}\label{lem A.1}
Let $\{\nu_u\}$ be a family of functions along the foliation $\{\Sigma_u\}$ of a null hypersurface $\ucalH$, we have
\begin{align*}
\uL \overline{\nu_u} = \overline{\uL \nu_u} + \overline{ \nu_u \tr \uchil{u}} - \overline{\nu_u}\, \overline{\tr \uchil{u}}.
\end{align*}
\end{lemma}
\begin{proof}
Since $\lie_{\uL} \dvol_u = \tr \uchil{u} \dvol_u$, 
then 
\begin{align*}
&
\textstyle
\uL \big(  \int_{\Sigma_u} \dvol_u \big) = \int_{\Sigma_u} \tr \uchil{u} \dvol_u,
\quad
\uL \big( \int_{\Sigma_u} \nu_u \dvol_u \big) =  \int_{\Sigma_u} (\uL \nu_u + \nu_u \tr \uchil{u} )\dvol_u.
\end{align*}
Substituting the above to the $\uL$ derivative of $\overline{\nu_u} = \frac{\int_{\Sigma_u} \nu_u \dvol_u}{\int_{\Sigma_u} \dvol_u}$ , the lemma follows.
\end{proof}

We introduce a coordinate system $\{u, \theta^1, \theta^2\}$ of $\ucalH$ where the $u$-level set is $\Sigma_u$ and $\uL \theta^1 = \uL \theta^2 =0$. Then $\uL = \partial_u$. The derivation of equation \eqref{eqn 2.11'} is as follows.
\begin{enumerate}[label=\alph*.]
\item Variation of the metric $\slashgl{u}$: 
$\uL \slashg_{ab} = 2 \uchi_{ab}$, 
$\uL (\slashg^{-1})^{ab} =-2 \uchi^{ab}$.

\item Variation of the Christoffel symbol $\slashGamma _{ab}^c$ of the Levi-Civita connection $\slashnabla$ of $(\Sigma_u, \slashg)$:
\begin{align*}
\uL \slashGamma_{ab}^c
=
(\slashg^{-1})^{cd}
( \slashnabla_a \uchi_{bd} + \slashnabla_b \uchi_{ad} - \slashnabla_d \uchi_{ab} ).
\end{align*}

\item Variation of the Gauss curvature $K$ of $(\Sigma_u, \slashg)$:
\begin{align*}
\uL K = \slashdiv \slashdiv \hatuchi - \frac{1}{2} \slashDelta \tr \uchi - K \tr \uchi.
\end{align*}

\item Variation of $\tr \chi'\, \tr \uchi$:
\begin{align*}
\uL ( \tr \chi'\, \tr \uchi )
=
- \tr \chi'\, ( \tr \uchi )^2 
- \tr \chi'\, |\hatuchi|^2
- 2 \tr \uchi\, | \eta |^2 
+ 2 \mu\, \tr \uchi.
\end{align*}

\item Variation of $\slashdiv \eta$: 
\begin{align*}
&
\left\{
\begin{aligned}
&
\lie_{\uL} \eta
=
- \hatuchi \cdot \eta
- \frac{1}{2} \tr \uchi\, \eta
+ \ubeta
+ 2 \slashd \uomega,
\\
&
\uL  (\slashnabla_a \eta_b )
= 
\slashnabla_a ( \lie_{\uL} \eta )_b
- 
(\slashg^{-1})^{cd}
( \slashnabla_a \uchi_{bd} + \slashnabla_b \uchi_{ad} - \slashnabla_d \uchi_{ab} ) \eta_c,
\\
&
\uL \slashdiv \eta
=
(\slashg^{-1})^{ab}\, \uL(\slashnabla_a \eta_b)
+ \uL ( \slashg^{-1} )^{ab}\, \slashnabla_a \eta_b
\end{aligned}
\right.
\\
\Rightarrow
&
\uL \slashdiv \eta 
=
-3 \slashdiv \hatuchi \cdot \eta 
- 3 ( \hatuchi, \slashnabla \eta)
- \frac{3}{2} \tr \uchi\, \slashdiv \eta
-\frac{1}{2} ( \slashd \tr \uchi, \eta)
+ \slashdiv \ubeta
+ 2 \slashDelta \uomega.
\end{align*}

\item Variation of $\mu$: equation \eqref{eqn 2.11'} follows from substituting c, d, e and equation \eqref{eqn 2.8} to $\uL \mu$.

\end{enumerate}

\Address

\end{document}